\documentclass[review,3p]{elsarticle}

\usepackage{mathtools}
\usepackage{amssymb}
\usepackage{mathrsfs}
\usepackage{graphicx}
\usepackage{algorithm}
\usepackage{algorithmic}
\usepackage{xcolor}
\usepackage{upgreek}
\usepackage{cases}
\usepackage{caption}
\usepackage{booktabs}
\usepackage[normalsize,nooneline,tight]{subfigure}
\usepackage[hidelinks]{hyperref}

\biboptions{sort&compress}

\newcommand{\Real}{\mathbb{R}}

\newcommand{\bu}{\boldsymbol{u}}

\newcommand{\bsigma}{\boldsymbol{\sigma}}

\newcommand{\bepsilon}{\boldsymbol{\epsilon}}
\newcommand{\bepsilone}{\bepsilon^\mathrm{e}}
\newcommand{\bepsilonp}{\bepsilon^\mathrm{p}}
\newcommand{\ld}{\lambda_{\alpha_1}}
\newcommand{\lp}{\lambda_{\kappa}}
\newcommand{\latt}{\lambda_{\alpha_2}}

\newcommand{\psiep}{\psi^{\mathrm{e}_+}}
\newcommand{\psiem}{\psi^{\mathrm{e}_-}}
\newcommand{\sigmad}{\boldsymbol{\sigma}_{\text{dev}}}
\newcommand{\nhat}{\hat{\boldsymbol{n}}}
\newcommand{\nnode}{n_\text{node}}
\newcommand{\ndof}{n_\text{DOF}}
\newcommand{\tp}{t^\text{prim}}
\newcommand{\td}{t^\text{dual}}

\newcommand{\sigmap}{\sigma^{\mathrm{p}}}
\newcommand{\etad}{\eta_\mathrm{d}}
\newcommand{\etap}{\eta_\mathrm{p}}

\newcommand{\buu}{\mathbf{u}}
\newcommand{\bul}{\boldsymbol{\lambda}}
\newcommand{\buL}{\boldsymbol{\Lambda}}
\newcommand{\Kuu}{\mathbf{K}_{\bu\bu}}
\newcommand{\Kuk}{\mathbf{K}_{\bu\kappa}}
\newcommand{\Kua}{\mathbf{K}_{\bu\alpha}}
\newcommand{\Kku}{\mathbf{K}_{\kappa\bu}}
\newcommand{\Kkk}{\mathbf{K}_{\kappa\kappa}}
\newcommand{\Kka}{\mathbf{K}_{\kappa\alpha}}
\newcommand{\Kau}{\mathbf{K}_{\alpha\bu}}
\newcommand{\Kaa}{\mathbf{K}_{\alpha\alpha}}
\newcommand{\Kak}{\mathbf{K}_{\alpha\kappa}}
\newcommand{\Kklp}{\mathbf{K}_{\kappa\lp}}
\newcommand{\Kald}{\mathbf{K}_{\alpha\ld}}
\newcommand{\Kalat}{\mathbf{K}_{\alpha\latt}}

\newcommand{\bld}{\boldsymbol{\lambda}_{\alpha_1}}
\newcommand{\blp}{\boldsymbol{\lambda}_{\kappa}}
\newcommand{\blat}{\boldsymbol{\lambda}_{\alpha_2}}
\newcommand{\Ld}{\boldsymbol{\Lambda}_{\alpha_1}}
\newcommand{\Lp}{\boldsymbol{\Lambda}_{\kappa}}
\newcommand{\Lat}{\boldsymbol{\Lambda}_{\alpha_2}}
\newcommand{\Sp}{\mathbf{S}_{\kappa}}
\newcommand{\Sd}{\mathbf{S}_{\alpha_1}}
\newcommand{\Sat}{\mathbf{S}_{\alpha_2}}
\newcommand{\Ns}{\mathbf{N}_{\mathrm{s}}}

\newcommand{\Bs}{\mathbf{B}_{\mathrm{s}}}
\newcommand{\Bv}{\mathbf{B}_{\mathrm{v}}}
\newcommand{\nhatu}{\hat{\mathbf{n}}}
\newcommand{\psiepu}{\uppsi^{\mathrm{e}_+}}

\begin{document}

\begin{frontmatter}

\title{Interior-point methods for the phase-field approach \\ to brittle and ductile fracture\tnoteref{t1,t2}}
\author[label1]{J.~Wambacq\corref{cor1}} 
\ead{jef.wambacq@kuleuven.be}
\author[label1]{J.~Ulloa}
\ead{jacintoisrael.ulloa@kuleuven.be}
\author[label1]{G.~Lombaert}
\ead{geert.lombaert@kuleuven.be}
\author[label1]{S.~Fran\c{c}ois}
\ead{stijn.francois@kuleuven.be}
\address[label1]{KU Leuven, Department of Civil Engineering, Leuven, Belgium}
\cortext[cor1]{Corresponding author: Jef Wambacq}

\tnotetext[t1]{{\itshape Postprint version.}}
\tnotetext[t2]{{\itshape Published version:}
        J.~Wambacq, J.~Ulloa, G.~Lombaert, and S.~Fran\c{c}ois.
       \newblock Interior-point methods for the phase-field approach to brittle and ductile fracture.
       \newblock {\em Computer Methods in Applied Mechanics and Engineering}, 375:113612, 2021.\\
       {\itshape DOI: \tt\url{https://doi.org/10.1016/j.cma.2020.113612}} \vspace{3mm}}

\begin{abstract}

The governing equations of the variational approach to brittle and ductile fracture emerge from the minimization of a non-convex energy functional subject to irreversibility constraints. This results in a multifield problem governed by a mechanical balance equation and evolution equations for the internal variables. While the balance equation is subject to kinematic admissibility of the displacement field, the evolution equations for the internal variables are subject to irreversibility conditions, and take the form of variational inequalities, which are typically solved in a relaxed or penalized way that can lead to deviations of the actual solution. This paper presents an interior-point method that allows to rigorously solve the system of variational inequalities. With this method, a sequence of perturbed constraints is considered, which, in the limit, recovers the original constrained problem. As such, no penalty parameters or modifications of the governing equations are involved. The interior-point method is applied in both a staggered and a monolithic scheme for both brittle and ductile fracture models. In order to stabilize the monolithic scheme, a perturbation is applied to the Hessian matrix of the energy functional. The presented algorithms are applied to three benchmark problems and compared to conventional methods, where irreversibility of the crack phase-field is imposed using a history field or an augmented Lagrangian.
\end{abstract}


\end{frontmatter}

\clearpage

\section{Introduction}
Many solids and structures exhibit complex material behavior where fracture combined with substantial ductility plays a dominant role. Examples include mild steels that show significant plastic deformations before the onset of fracture in various engineering applications. A variational approach to the modeling of brittle and ductile fracture provides a rigorous mathematical setting that allows to describe these phenomena in a physically sound manner.

Such a variational approach to quasi-static brittle fracture was introduced by Francfort and Marigo~\cite{fran98a}. This approach overcomes the limitations of the classical Griffith theory, and allows modeling complex crack topologies as well as the initiation, propagation, coalescence and branching of cracks without additional ad-hoc criteria. The method leads to a free discontinuity problem, which, in most cases, is intractable to solve. Therefore, Bourdin~et~al.~\cite{bour00a,bour08a} proposed to approximate discontinuities in the displacement field using a phase field, leading to a non-convex energy functional, which gave rise to the regularized variational approach to brittle fracture. Recently, several authors proposed extensions of this approach to model ductile fracture~\cite{ales14a,duda15a,amba15a,amba16a,bord16a,kuhn16a,mieh16a,mieh16b,mieh16c,rodr18a,ales18b} and ductile fatigue~\cite{ullo21a}, where plastic deformations can occur together with the fracture process. A review of phase-field models of fracture coupled with plasticity was conducted by Alessi~et~al.~\cite{ales18a}.

In the literature, various algorithms have been presented to minimize the non-convex energy functional and to solve the governing equations of the phase-field approach to brittle and ductile fracture. Monolithic schemes minimize the energy functional for all arguments simultaneously. However, due to the non-convexity, convergence issues are often encountered~\cite{heis15a,gera16a}. Therefore, staggered schemes, based on alternate minimization, which exploit the convexity of the energy functional with respect to each individual argument, are often adopted~\cite{bour08a,mieh10b,rodr18a}. This results in an algorithmic decoupling of the governing equations, which can be independently solved in a robust way. Convergence proofs of the alternate minimization algorithm are presented in~\cite{bour07a,burk10a,brun20a}. The efficiency of the staggered and monolithic schemes, however, is problem-dependent. Alternatively, Farrell and Maurini~\cite{farr17a} suggest a hybrid approach, where a staggered scheme is used to obtain a state close to a local minimum of the energy functional, followed by a minimization with a monolithic scheme.

Based on thermodynamical grounds, it is required to impose irreversibility conditions, leading to a system of variational inequalities~\cite{pham11a,gera19a}. These conditions can also be shown to be a consequence of the stronger Drucker-Illyushin postulate, as thoroughly discussed in Marigo~\cite{mari89a,mari00a} and Pham and Marigo~\cite{pham10a}. To enforce irreversibility in the case of brittle fracture, Bourdin~\cite{bour07a} proposed setting Dirichlet-type conditions on the phase field to prevent crack healing in fully developed cracks. Similar formulations \cite{arti15a,gera16a} followed later, where the Dirichlet condition is imposed in an approximated way. These approaches, however, are not sufficient to avoid crack healing in partially damaged regions, leading to violations of the irreversibility constraints. To prevent crack healing in these regions, Miehe~et~al.~\cite{mieh10a} proposed the use of a smooth penalty term. This was later extended by Wheeler~et~al.~\cite{whee14a} by adopting an augmented Lagrangian which uses the Moreau-Yosida approximation of an indicator function to prevent ill-conditioning caused by the penalization. The augmented Lagrangian approach has since been used in many follow-up works, e.g.~\cite{wick17a,wick17b,geel19a,brun20a}. However, both the penalty approach of Miehe et al.~\cite{mieh10a} and the augmented Lagrangian approach of Wheeler~et~al.~\cite{whee14a} require a parameter to penalize violations of the irreversibility constraints. A too small penalty parameter is insufficient to impose the constraint and will not prevent crack healing, while a penalty parameter that is too large often leads to ill-conditioning. Therefore, Gerasimov and De~Lorenzis~\cite{gera19a} analytically derived a lower bound for the penalty parameter, which is sufficient to enforce the constraint for brittle fracture. Heister~et~al.~\cite{heis15a} depart from the penalization approach and use a primal-dual active set strategy, based on the works of Hinterm{\"u}ller~et~al.~\cite{hint02a}. This approach can be identified as a semi-smooth Newton method, which, under certain assumptions, leads to locally super-linear convergence rates. Similarly, several authors~\cite{leon14a,farr17a,tann18b} solve the bound-constrained problem using the reduced-space active-set method of Benson and Munson~\cite{bens06a} or the bounded trust region method provided in PETSc~\cite{petsc-user-ref}. Alternatively, Lancioni and Royer-Carfagni~\cite{lanc09a} apply an a posteriori projection technique, similar to the way irreversibility is treated in the level set method of Allaire~et~al.~\cite{alla07a}. This technique leads to results which do not necessarily satisfy the system of variational inequalities~\cite{amor09a}. Miehe~et~al.~\cite{mieh10b} impose the condition of irreversibility implicitly by introducing a history field, which replaces the tensile elastic energy density in the damage evolution equation by its maximal value attained within the loading history. This, however, destroys the variational structure of the problem, and is not proven to rigorously impose the irreversibility conditions.

In the case of ductile fracture, several authors~\cite{amba15a,duda15a,amba16a,kuhn16a,bord16a} adopt standard local return-mapping algorithms to solve the plastic evolution problem. This approach can only be used for local plasticity and, when combined with a crack phase-field, can lead to mesh-sensitive results or unrealistic localization patterns~\cite{mieh16c}. In the case of ductile fracture with non-local plasticity, Miehe~et~al.~\cite{mieh16a,mieh16b} consider a length scale for both an isotropic hardening variable and the crack phase-field. Thermodynamic consistency is then taken care of by relying on a viscous regularization of the problem, leading to an unconstrained formulation. Alternatively, Ulloa~et~al.~\cite{ullo16a,ullo21a} and Rodr{\'i}guez~et~al.~\cite{rodr18a} apply the a posteriori projection technique of Lancioni and Royer-Carfagni~\cite{lanc09a} to both the crack phase-field and an isotropic hardening variable. In addition, similar to brittle fracture, a smooth penalty term can be adopted to solve the plastic evolution problem. However, the downside of this approach is that a lower bound for the required penalty is not known, and two independent penalty parameters for damage and plasticity are required.

Modifications of the governing equations or inadequate penalization can lead to deviations from the solution of the original rate-independent problem. Therefore, the use of interior-point methods is particularly attractive since they allow to rigorously solve the constrained equations in both staggered and monolithic schemes without requiring a penalty parameter or modifications of the governing equations. In addition, interior-point methods can be used for large-scale optimization problems, since, generally, the number of iterations that is required to obtain an optimal solution depends very little on the problem dimension~\cite{gond12a}.

In the context of computational mechanics, interior-point methods have been adopted mainly to solve problems with plasticity or contact mechanics. Therein, the system of governing equations is typically cast in a conic optimization problem, which is then solved using an interior-point method. Applications include limit analysis~\cite{makr06a,makr07a}, shakedown analysis~\cite{vu04a,makr06b,nguy08a}, granular flow~\cite{krab12a,zhan14a}, contact problems~\cite{mang18a,elbo20b}, viscoplasticity~\cite{bley15a,bley18a}, and elastoplasticity~\cite{yone11a,sche20a,elbo20a}. Krabbenh\o ft~et al.~\cite{krab07b} present a primal-dual interior-point method for the solution of elastoplastic problems, while Vavasis et al.~\cite{vava20a} recently proposed a primal interior-point method for the solution of problems governed by cohesive fracture. Inspired by these developments, the present work presents a primal-dual interior-point method for the phase-field approach to brittle and ductile fracture.

The outline of this paper is as follows. In section~\ref{sec:notation}, the notation and the governing equations are presented. Next, monolithic and staggered schemes based on interior-point methods are elaborated for the rigorous solution of the brittle and ductile phase-field problem in section~\ref{sec:interiorpoint}. Section~\ref{sec:examples} presents three numerical examples illustrating the results and systematically compares the staggered and monolithic schemes. In addition, the results are benchmarked against results obtained using the history field of Miehe~et~al.~\cite{mieh10b} and the augmented Lagrangian approach of Wheeler et al.~\cite{whee14a}. The results show that interior-point methods can be successfully applied to solve the equations that govern the phase-field approach to brittle and ductile fracture, in both staggered and monolithic schemes. Finally, the work is concluded in section~\ref{sec:conclusions}.

\section{Notation and governing equations} \label{sec:notation}
This section summarizes the governing equations of the phase-field approach to ductile fracture presented by Rodr{\'i}guez~et~al.~\cite{rodr18a}. The governing equations are an extension of the classical equations of the phase-field approach to brittle fracture~\cite{bour00a,bour08a}, and recover the original model when plastic effects are omitted. Adopting a variational approach, the governing equations are obtained by taking the directional derivatives of an energy functional. In a second step, these equations are discretized in time, which results in an incremental system of equations that is suitable for numerical implementation. The incremental equations are discretized using finite elements, resulting in a constrained system of discrete non-linear equations.

\subsection{Energy functional}
Consider a solid domain $\Omega\subset\Real^3$, on which displacements $\bar{\bu}(\boldsymbol{x},t)$ are imposed on the Dirichlet boundary $\Gamma_\mathrm{D}\subset\Real^3$. No external forces are considered, and the deformation process is assumed quasi-static. In addition, all deformations are assumed to be infinitesimally small, except in certain localized regions. In this setting, the aim of this paper is to present an algorithm which allows to compute the displacement field $\bu:\Omega\times T\rightarrow\Real^3$, the plastic strain tensor $\bepsilonp:\Omega\times T\rightarrow\Real^{3\times3}_\text{dev}$, and the crack phase-field $\alpha:\Omega\times T\rightarrow[0,1]$ in the solid over a time interval $T\coloneqq[0,t_\text{max}]$, where $\Real^{3\times3}_\text{dev}\coloneqq\{\boldsymbol{e}\in\Real^{3\times3}_\text{sym}\,\vert\, \text{tr}\,\boldsymbol{e}=0\}$ is the set of symmetric and deviatoric second-order tensors.

The governing equations of the coupled damage-plasticity model are obtained from an energy functional by applying the energetic formulation for rate-independent systems~\cite{miel06a}. This framework postulates the following two principles: (a)~a~stability condition, leading to the weak forms of the balance equation and the yield criteria, and (b) energy balance, leading to the consistency conditions of both the damage and the plastic evolution problems. In addition, thermodynamic consistency leads to irreversibility constraints. Rodr{\'i}guez et~al.~\cite{rodr18a} propose to minimize a functional with respect to $\{\bu,\bepsilonp,\kappa,\alpha\}$, where $\kappa:\Omega\times T\rightarrow\Real_+$ is the equivalent plastic strain evolving as $\dot{\kappa}\coloneqq\sqrt{2/3}\,\Vert\dot{\bepsilon}^\mathrm{p}\Vert$. The functional reads
\begin{equation}
\Pi(\bu,\bepsilonp,\kappa,\alpha)\coloneqq\int_\Omega \bigg[g(\alpha)\psiep(\bepsilon,\bepsilonp)+\psiem(\bepsilon,\bepsilonp)+w(\alpha)+\frac{1}{2}\etad^2\Vert\nabla \alpha\Vert^2 + g(\alpha)\Big(\sigmap\kappa+\frac{1}{2}\etap^2\Vert\nabla \kappa\Vert^2\Big)\bigg] \,\mathrm{d}\boldsymbol{x},
\label{eq:functional}
\end{equation}
where $\bepsilon\coloneqq\nabla^\mathrm{s}\bu=\frac{1}{2}(\nabla\otimes\bu+\bu\otimes\nabla)$ is the second-order infinitesimal strain tensor, ${g(\alpha)\coloneqq(1-\alpha)^2}$ is the usual quadratic degradation function, $\Vert\square\Vert\coloneqq\Vert\square\Vert_{\ell^2}$ is the Euclidean norm, $\sigmap$ is the yield strength, $\etad$ is a length scale that governs the localization of damage, and $\etap$ is a length scale that governs the localization of the plastic strains before the onset of damage. For notational simplicity, explicit dependence on space and time has been omitted in Eq.~\eqref{eq:functional}, as will be done hereinafter. In the energy functional, perfect plasticity with gradient regularization is considered. Additional mechanisms, such as isotropic hardening, can be incorporated in a straightforward way~\cite{rodr18a}.

In Eq.~\eqref{eq:functional}, the elastic energy density is decomposed into so-called positive and negative parts, denoted by $\psiep(\bepsilon,\bepsilonp)$ and $\psiem(\bepsilon,\bepsilonp)$, respectively, enabling to only degrade the positive elastic energy part. Various definitions of the positive and negative parts are presented in the literature~\cite{amor09a,mieh10a,fred10a,stei19a}. Similar to the ductile fracture model presented in~\cite{rodr18a}, the volumetric-deviatoric decomposition proposed by Amor~et~al.~\cite{amor09a} is used, where the elastic energy density is decomposed as
\begin{equation}
\begin{dcases}
\psiep(\bepsilon,\bepsilonp)\coloneqq\frac{1}{2}K\big<\text{tr}(\bepsilone)\big>_+^2+G(\bepsilone_{\text{dev}}:\bepsilone_{\text{dev}}), \\
\psiem(\bepsilon,\bepsilonp)\coloneqq\frac{1}{2}K\big<\text{tr}(\bepsilone)\big>_-^2,
\end{dcases}
\label{eq:split}
\end{equation}
where $K$ is the bulk modulus, $G$ is the shear modulus, $\big<\square\big>_\pm\coloneqq\frac{1}{2}(\square\pm|\square|)$ is the ramp function, and $\bepsilone_{\text{dev}}:\Omega\times T\rightarrow\Real^{3\times3}_\text{dev}$ is the deviatoric part of the elastic strain tensor $\bepsilone\coloneqq\bepsilon-\bepsilonp$. From the decomposed energy density, the stress tensor is determined as follows:
\begin{equation}
\bsigma(\bepsilon,\bepsilonp,\alpha)= g(\alpha)\frac{\partial\psiep}{\partial\bepsilon}(\bepsilon,\bepsilonp)+\frac{\partial\psiem}{\partial\bepsilon}(\bepsilon,\bepsilonp).\label{eq:stress}
\end{equation}

In Eq.~\eqref{eq:functional}, $w(\alpha)$ represents the dissipation due to the local damage evolution, for which two options are typically adopted:
\begin{equation}
w(\alpha)\coloneqq
\begin{dcases*}
w_0\alpha & AT-1, \\
w_0\alpha^2 & AT-2,
\end{dcases*}
\end{equation}
where $w_0$ is a local damage dissipation constant. The AT-1 model leads to an elastic regime before the onset of damage or plasticity, whereas in the AT-2 model, damage starts to evolve immediately upon loading~\cite{mari16a,gera19a}.

Using these definitions, the plastic length scale $\etap=\ell_\mathrm{p}\sqrt{\sigmap}$ and damage length scale $\etad=\ell_\mathrm{d}\sqrt{2w_0}$ can be related to the plastic characteristic length $\ell_\mathrm{p}$, and the damage characteristic length $\ell_\mathrm{d}$, respectively. Note that these two parameters have the units of length.

\subsection{Governing equations in continuous form}
To derive the governing equations from the energy functional~\eqref{eq:functional}, the following function spaces are defined:
\begin{equation}
\begin{alignedat}{2}
\mathscr{U}&\coloneqq\{\boldsymbol{u}\in H^1(\Omega;\Real^3)\,\vert\,\boldsymbol{u}=\bar{\bu} \text{ on } \Gamma_{\mathrm{D}}\},
&&\quad\tilde{\mathscr{U}}\coloneqq\{\tilde{\boldsymbol{u}}\in H^1(\Omega;\Real^3)\,\vert\,\tilde{\boldsymbol{u}}=\boldsymbol{0} \text{ on } \Gamma_{\mathrm{D}}\}, \\
\mathscr{B}&\coloneqq L^2(\Omega;\Real^{3\times3}_\text{dev}),
&&\quad\tilde{\mathscr{B}}\coloneqq\mathscr{B}, \\
\mathscr{K}&\coloneqq H^1(\Omega;\Real_+),
&&\quad\tilde{\mathscr{K}}(\tilde{\boldsymbol{p}})\coloneqq\{\tilde{\xi}\in\mathscr{K}\,\vert\,\tilde{\xi}=\sqrt{2/3}\,\Vert\tilde{\boldsymbol{p}}\Vert,\,\tilde{\boldsymbol{p}}\in\tilde{\mathscr{B}}\}, \\
\mathscr{D}&\coloneqq H^1(\Omega;[0,1]),
&&\quad\tilde{\mathscr{D}}\coloneqq H^1(\Omega;\Real_+).
\end{alignedat}
\label{eq:functionspaces}
\end{equation}
The set $\{\bu,\bepsilonp,\kappa,\alpha\}\in\mathscr{U}\times\mathscr{B}\times\mathscr{K}\times\mathscr{D}$ contains admissible displacement fields and internal variables, considering kinematic admissibility of the displacement field and boundedness of the equivalent plastic strain and the crack phase-field. Similarly, the set $\{\tilde{\bu},\tilde{\bepsilon}^\mathrm{p},\tilde{\kappa},\tilde{\alpha}\}\in\tilde{\mathscr{U}}\times\tilde{\mathscr{B}}\times\tilde{\mathscr{K}}\big(\tilde{\bepsilon}^\mathrm{p}\big)\times\tilde{\mathscr{D}}$ contains admissible virtual displacement fields and internal variables, considering the monotonicity constraint on the crack phase-field, as well as a constraint linking the variation of the equivalent plastic strain to the variation of the plastic strain tensor, i.e.~$\tilde{\kappa}=\sqrt{2/3}\,\Vert\tilde{\bepsilon}^\mathrm{p} \Vert$. The weak forms of the governing equations, which hold for almost all times $t\in I$, are obtained by taking directional derivatives of the energy functional~\eqref{eq:functional} with respect to $\{\bu,\bepsilonp,\kappa,\alpha\}$, leading to
\begin{equation}
\begin{dcases}
\int_\Omega \bsigma(\bepsilon,\bepsilonp,\alpha):\nabla^\mathrm{s}\tilde{\bu}\,\mathrm{d}\boldsymbol{x}=0, \\
\int_\Omega\bigg[\bigg(-\sqrt{\frac{3}{2}}\big\Vert\sigmad(\bepsilon,\bepsilonp,\alpha)\big\Vert+g(\alpha)\sigmap\bigg)\tilde{\kappa}+g(\alpha)\etap^2\nabla\kappa\cdot\nabla\tilde{\kappa}\bigg]\,\mathrm{d}\boldsymbol{x}\geq0, \\
\int_\Omega \bigg[\bigg(g'(\alpha)\bigg(\psiep(\bepsilon,\bepsilonp)+\sigmap\kappa+\frac{1}{2}\etap^2\Vert\nabla \kappa\Vert^2\bigg)+w'(\alpha)\bigg)\tilde{\alpha}+\etad^2\nabla\alpha\cdot\nabla\tilde{\alpha}\bigg]\,\mathrm{d}\boldsymbol{x}\geq0,
\end{dcases}
\end{equation}
for all $\{\tilde{\bu},\tilde{\bepsilon}^\mathrm{p},\tilde{\kappa},\tilde{\alpha}\}\in\tilde{\mathscr{U}}\times\tilde{\mathscr{B}}\times\tilde{\mathscr{K}}\big(\tilde{\bepsilon}^\mathrm{p}\big)\times\tilde{\mathscr{D}}$, where the optimality conditions are satisfied for $\bsigma:\tilde{\bepsilon}^\mathrm{p}=\Vert\sigmad\Vert\Vert\tilde{\bepsilon}^\mathrm{p}\Vert=\sqrt{3/2}\,\Vert\sigmad\Vert\tilde{\kappa}$, and $\sigmad$ is the deviatoric part of the stress tensor $\boldsymbol{\sigma}$ given in Eq.~\eqref{eq:stress}. The evolution law of the equivalent plastic strain (i.e.~$\dot{\kappa}\coloneqq\sqrt{2/3}\,\Vert\dot{\bepsilon}^\mathrm{p}\Vert$) completes the system of variational inequalities. The system corresponds to the weak forms of the balance equation, the evolution equations for plasticity and damage, and the boundary conditions for $\nabla\kappa$ and $\nabla\alpha$, where kinematic admissibility and boundedness of the internal variables are ensured by seeking solutions in the function spaces~\eqref{eq:functionspaces} and considering irreversibility of the crack phase-field.

\subsection{Governing equations in incremental form}
In order to obtain a numerical step-by-step time integration of the governing equations, the time interval $T$ is discretized using a pseudo-time step parameter $n\in\mathbb{N}\,\cup\,\{0\}$. A quantity evaluated at the previous time step $t_n$, where all variables are known, is denoted with a subscript $n$, while a quantity evaluated at $t_{n+1}$ is written without a subscript for brevity. The energy functional~\eqref{eq:functional} can then be modified using an indicator function, explicitly considering the boundedness of the crack phase-field as follows:
\begin{equation}
\widehat{\Pi}\big(\bu,\bepsilonp,\kappa,\alpha;\alpha_n\big)\coloneqq\Pi\big(\bu,\bepsilonp,\kappa,\alpha\big)+\int_\Omega I_{[\alpha_n,1]}(\alpha)\,\mathrm{d}\boldsymbol{x},
\label{eq:functionalindicator}
\end{equation}
where $I_\mathrm{C}(\square)$ is an indicator function that equals 0 if $\square\in \mathrm{C}$, and $+\infty$ otherwise. By adding the indicator function to the functional, the admissible space for the crack phase-field is no longer constrained and its function space is relaxed to $H^1(\Omega;\Real)$.

Given the modified energy functional~\eqref{eq:functionalindicator}, the governing equations are obtained by taking the directional derivatives with respect to $\{\bu,\bepsilonp,\kappa,\alpha\}$, considering the corresponding function spaces, which results in:
\begin{subnumcases}{\label{eq:directionalderivativesincremental}}
\int_\Omega \bsigma(\bepsilon,\bepsilonp,\alpha):\nabla^\mathrm{s}\tilde{\bu}\,\mathrm{d}\boldsymbol{x}=0, \label{eq:directionalderivativesincremental1} \\
\int_\Omega\bigg[\bigg(-\sqrt{\frac{3}{2}}\big\Vert\sigmad(\bepsilon,\bepsilonp,\alpha)\big\Vert+g(\alpha)\sigmap\bigg)\tilde{\kappa}+g(\alpha)\etap^2\nabla\kappa\cdot\nabla\tilde{\kappa}\bigg]\,\mathrm{d}\boldsymbol{x}\geq0, \label{eq:directionalderivativesincremental2} \\
\int_\Omega \bigg[\bigg(g'(\alpha)\bigg(\psiep(\bepsilon,\bepsilonp)+\sigmap\kappa+\frac{1}{2}\etap^2\Vert\nabla \kappa\Vert^2\bigg)+w'(\alpha)+\partial I_{[\alpha_n,1]}(\alpha)\bigg)\tilde{\alpha}+\etad^2\nabla\alpha\cdot\nabla\tilde{\alpha}\bigg]\,\mathrm{d}\boldsymbol{x}\ni0, \label{eq:directionalderivativesincremental3}
\end{subnumcases}
for all $\{\tilde{\bu},\tilde{\bepsilon}^\mathrm{p},\tilde{\kappa},\tilde{\alpha}\}\in\tilde{\mathscr{U}}\times\tilde{\mathscr{B}}\times\tilde{\mathscr{K}}\big(\tilde{\bepsilon}^\mathrm{p}\big)\times\tilde{\mathscr{D}}$ using $\bsigma:\tilde{\bepsilon}^\mathrm{p}=\sqrt{3/2}\,\Vert\sigmad\Vert\tilde{\kappa}$, where $\partial I_\mathrm{C}(\square)$ is the subdifferential of the indicator function, which gives 0 if $\square\in \mathrm{C}\setminus\partial \mathrm{C}$, $\Real_-$ if $\square\in\partial \mathrm{C}$, and $\varnothing$ otherwise.

In the incremental setting, the plastic strain tensor can be eliminated from the system in Eq.~\eqref{eq:directionalderivativesincremental} by using an incremental form of the associative flow rule, i.e.\ $\bepsilonp=\bepsilonp_n+\sqrt{3/2}\nhat(\kappa-\kappa_n)$, where ${\nhat\coloneqq\sigmad/\left\Vert\sigmad\right\Vert}$ is the direction of the plastic flow. Similar to the classical radial return method~\cite{simo98a}, an implicit scheme is used, where the direction of the plastic flow at the stress state at $t_{n+1}$ and the trial state coincide, i.e.\ $\nhat=\sigmad^\mathrm{tr}/\left\Vert\sigmad^\mathrm{tr}\right\Vert$, where $\sigmad^\mathrm{tr}\coloneqq2g(\alpha)G(\bepsilon_\text{dev}-\bepsilonp_n)$ is the deviatoric part of the trial stress tensor. In addition, the equivalent plastic strain must be constrained to $[\kappa_n,+\infty)$. As such, the plastic evolution problem~\eqref{eq:directionalderivativesincremental2} can be stated as
\begin{equation}
\int_\Omega\bigg[\bigg(-\sqrt{\frac{3}{2}}\big\Vert\sigmad^\mathrm{tr}(\bepsilon,\bepsilonp_n,\alpha)\big\Vert+g(\alpha)(\sigmap+3G(\kappa-\kappa_n))+\partial I_{[\kappa_n,+\infty)}(\kappa)\bigg)\tilde{\kappa}+g(\alpha)\etap^2\nabla\kappa\cdot\nabla\tilde{\kappa}\bigg]\,\mathrm{d}\boldsymbol{x}\ni0.
\label{eq:plasticincremental}
\end{equation}

\subsection{Solution of the constrained equations}
Various solution techniques can be adopted in order to solve the constrained Eqs.~\eqref{eq:directionalderivativesincremental3} and~\eqref{eq:plasticincremental}. Hereafter, the history field approach of Miehe et al.~\cite{mieh10b}, the augmented Lagrangian approach of Wheeler~et~al.~\cite{whee14a} and the proposed interior-point method will be considered.

\paragraph{History field} Miehe et al.~\cite{mieh10b} propose to replace the positive elastic energy density in the weak form of the damage evolution equation by its maximal value attained within the loading history. Later, several authors (e.g.~\cite{amba15a,bord16a}) applied a similar technique for ductile fracture, as will be done hereafter. The weak form of the modified damage evolution equation reads
\begin{equation}
\int_\Omega \bigg[\bigg(g'(\alpha)\mathcal{H}(\bepsilon,\bepsilonp,\kappa)+w'(\alpha)\bigg)\tilde{\alpha}+\etad^2\nabla\alpha\cdot\nabla\tilde{\alpha}\bigg]\,\mathrm{d}\boldsymbol{x}=0,
\end{equation}
where $\mathcal{H}(\bepsilon,\bepsilonp,\kappa)\coloneqq\max\limits_{t\in[0,t_{n+1}]} \big(\psiep(\bepsilon,\bepsilonp)+\sigmap\kappa+\frac{1}{2}\etap^2\Vert\nabla \kappa\Vert^2\big)$ denotes the history field. This technique can only be used to solve the damage evolution problem~\eqref{eq:directionalderivativesincremental3}, and a different method is required for the solution of the plastic evolution problem~\eqref{eq:plasticincremental}.

\paragraph{Augmented Lagrangian}
Wheeler et al.~\cite{whee14a} suggest to replace the indicator function of the crack phase-field by its smooth Moreau-Yosida approximation, which reads
\begin{equation}
I_{[\alpha_n,1]}(\alpha)\approx\frac{1}{2\gamma}\big<\Xi+\gamma(\alpha-\alpha_n)\big>_-^2+\frac{1}{2\gamma}\big<\Xi+\gamma(\alpha-1)\big>_+^2-\frac{1}{2\gamma}\Xi^2,
\end{equation}
where $\gamma\in\Real_+$ is a sufficiently high penalty parameter and $\Xi$ is an auxiliary field obtained iteratively (see Wheeler et al.~\cite{whee14a}). Using the Moreau-Yosida approximation of the indicator function, the weak form of the damage evolution equation becomes
\begin{equation}
\int_\Omega \bigg[\bigg(g'(\alpha)\bigg(\psiep(\bepsilon,\bepsilonp)+\sigmap\kappa+\frac{1}{2}\etap^2\Vert\nabla \kappa\Vert^2\bigg)+w'(\alpha)+\big<\Xi+\gamma(\alpha-\alpha_n)\big>_-+\big<\Xi+\gamma(\alpha-1)\big>_+\bigg)\tilde{\alpha}+\etad^2\nabla\alpha\cdot\nabla\tilde{\alpha}\bigg]\,\mathrm{d}\boldsymbol{x}=0.
\end{equation}
The augmented Lagrangian approach can also be adopted for the solution of the plastic evolution problem~\eqref{eq:plasticincremental}, given that a suitable penalization constant can be determined.

\paragraph{Interior-point method} For the interior-point method, the system of variational inequalities~\eqref{eq:directionalderivativesincremental3}~and~\eqref{eq:plasticincremental} is first cast in a system of variational equalities using the definition of the subdifferential of the indicator functions, similar to the procedure in Mang et al.~\cite{mang20a}:
\begin{subnumcases}{\label{eq:directionalderivativesreduced}}
\hspace{-1mm}\int_\Omega\hspace{-0.5mm}\bigg[\bigg(-\sqrt{\frac{3}{2}}\big\Vert\sigmad^\mathrm{tr}(\bepsilon,\bepsilonp_n,\alpha)\big\Vert+g(\alpha)(\sigmap+3G(\kappa-\kappa_n))-\lp\bigg)\tilde{\kappa}+g(\alpha)\etap^2\nabla\kappa\cdot\nabla\tilde{\kappa}\bigg]\,\mathrm{d}\boldsymbol{x}=0, \label{eq:directionalderivativesplastic} \\
\hspace{-1mm}\int_\Omega\hspace{-0.5mm}\bigg[\bigg(g'(\alpha)\bigg(\psiep(\bepsilon,\bepsilonp(\kappa))\hspace{-0.5mm}+\hspace{-0.5mm}\sigmap\kappa+\frac{1}{2}\etap^2\Vert\nabla \kappa\Vert^2\bigg)\hspace{-0.5mm}+\hspace{-0.5mm}w'(\alpha)\hspace{-0.5mm}-(\ld\hspace{-0.5mm}-\hspace{-0.5mm}\latt)\bigg)\tilde{\alpha}+\etad^2\nabla\alpha\cdot\nabla\tilde{\alpha}\bigg]\,\mathrm{d}\boldsymbol{x}=0, \label{eq:directionalderivativesdamage}
\end{subnumcases}
where $\lp$, $\ld$ and $\latt$ act as Lagrange multipliers satisfying the following KKT conditions:
\begin{equation}
\begin{dcases}
\lp(\kappa-\kappa_n)=0,&\lp\geq0,\hspace{5.56mm}(\kappa-\kappa_n)\geq0, \\
\ld(\alpha-\alpha_n)=0,&\ld\geq0,\hspace{4mm}(\alpha-\alpha_n)\geq0, \\
\latt(1-\alpha)=0,&\latt\geq0,\hspace{4mm}(1-\alpha)\geq0.
\end{dcases}
\label{eq:kktconditions}
\end{equation}
The KKT conditions express (a) consistency conditions for the evolution of the plastic strains and the crack phase-field, (b) that the state of the solid must be admissible, and (c) thermodynamic consistency and boundedness of the crack phase-field. In section~\ref{sec:interiorpoint}, an interior-point method will be presented, which allows to solve the constrained system of equations. Therein, a sequence of perturbed KKT conditions is considered, which in the limit recovers the unperturbed conditions. A damped Newton method is then used to solve the perturbed problem and to stay within the feasible region given by the inequality constraints.

\subsection{Spatial discretization and system matrices}
The governing equations in weak form (Eqs.~\eqref{eq:directionalderivativesincremental1} and~\eqref{eq:directionalderivativesreduced}) are spatially discretized using finite elements with $\nnode$ nodes and $\ndof$ degrees of freedom, resulting in the following discrete system of non-linear equations in the nodal quantities $\buu\in\Real^{\ndof}$ and $\boldsymbol{\kappa}$, $\boldsymbol{\alpha}$, $\blp$, $\bld$, $\blat \in\Real^{\nnode}$:
\begin{equation}
\begin{dcases}
\mathbf{r}_{\bu}\hspace{-0.7mm}\coloneqq\hspace{-0.7mm}\int_\Omega\Bv^\mathrm{T}((1-\Ns\boldsymbol{\alpha})^2\mathbf{C}^++\mathbf{C}^-)(\Bv\buu-\bepsilonp)\,\mathrm{d}\boldsymbol{x}=\mathbf{0}, \\
\mathbf{r}_{\kappa}\hspace{-0.7mm}\coloneqq\hspace{-0.7mm}\int_\Omega \bigg[\Ns^\mathrm{T}\bigg(\hspace{-1.4mm}-\hspace{-0.5mm}\sqrt{\frac{3}{2}}\big\Vert\sigmad^\mathrm{tr}\big\Vert\hspace{-0.5mm}+\hspace{-0.5mm}(1\hspace{-0.5mm}-\hspace{-0.5mm}\Ns\boldsymbol{\alpha})^2(\sigmap\hspace{-0.5mm}+\hspace{-0.5mm}3G\Ns(\boldsymbol{\kappa}\hspace{-0.5mm}-\hspace{-0.5mm}\boldsymbol{\kappa}_n))\hspace{-0.5mm}-\hspace{-0.5mm}\Ns\blp\bigg)\hspace{-0.5mm}+\hspace{-0.5mm}\Bs^\mathrm{T}(1\hspace{-0.5mm}-\hspace{-0.5mm}\Ns\boldsymbol{\alpha})^2\etap^2\Bs\boldsymbol{\kappa}\bigg]\,\mathrm{d}\boldsymbol{x}=\mathbf{0}, \\
\mathbf{r}_{\alpha}\hspace{-0.7mm}\coloneqq\hspace{-0.7mm}\int_\Omega \bigg[\Ns^\mathrm{T}\bigg(\hspace{-1.4mm}-\hspace{-0.5mm}2(1\hspace{-0.5mm}-\hspace{-0.5mm}\Ns\boldsymbol{\alpha})\bigg(\hspace{-0.5mm}\psiepu\hspace{-0.5mm}+\hspace{-0.5mm}\sigmap\Ns\boldsymbol{\kappa}\hspace{-0.5mm}+\hspace{-0.5mm}\frac{1}{2}\etap^2\Vert\Bs \boldsymbol{\kappa}\Vert^2\hspace{-0.5mm}\bigg)\hspace{-0.8mm}+\hspace{-0.5mm}w'(\Ns\boldsymbol{\alpha})\hspace{-0.5mm}-\hspace{-0.5mm}\Ns(\bld\hspace{-0.5mm}-\hspace{-0.5mm}\blat)\hspace{-0.5mm}\bigg)\hspace{-0.6mm}+\hspace{-0.5mm}\Bs^\mathrm{T}\etad^2\Bs\boldsymbol{\alpha}\bigg]\,\mathrm{d}\boldsymbol{x}\hspace{-0.4mm}=\hspace{-0.4mm}\mathbf{0}.
\end{dcases}
\label{eq:weakformdisc}
\end{equation}
In Eq.~\eqref{eq:weakformdisc}, $\Ns\in\Real^{1\times\nnode}$ is the shape function matrix, and $\Bs\in\Real^{3\times\nnode}$ and ${\Bv\in\Real^{6\times\ndof}}$ are the discrete scalar and vector differential operators, respectively. Moreover, $\psiepu\coloneqq\frac{1}{2}\Vert\Bv\buu-\bepsilonp\Vert^2_{\mathbf{C}^+}\in\Real_+$ denotes the discrete positive part of the elastic energy density, where, in Voigt notation, $\bepsilonp\coloneqq\bepsilonp_n+\sqrt{3/2}\nhatu^\mathrm{tr}\Ns(\boldsymbol{\kappa}-\boldsymbol{\kappa}_n)\in\Real^6$ is the plastic strain vector, $\sigmad^\mathrm{tr}\in\Real^6$ is the deviatoric part of the trial stress, $\nhatu^\mathrm{tr}\in\Real^6$ is the direction of the plastic flow at the trial state, and $\mathbf{C}^+\in\Real^{6\times 6}_\text{sym}$ and $\mathbf{C}^-\in\Real^{6\times6}_\text{sym}$ are the constitutive matrices corresponding to the decomposition of the elastic energy density in Eq.~\eqref{eq:split}. The function $w'(\Ns\boldsymbol{\alpha})$ equals $w_0$ for the AT-1 model and $2w_0\Ns\boldsymbol{\alpha}$ for the AT-2 model. The Jacobian matrices $\mathbf{K}_{\square\raisebox{0.1mm}{$\bullet$}}\coloneqq\frac{\partial\mathbf{r}_\square}{\partial\bullet}$ of the system in Eq.~\eqref{eq:weakformdisc} read 
\begin{equation}
\begin{alignedat}{1}
&\Kuu=\int_\Omega \Bv^\mathrm{T}((1-\Ns\boldsymbol{\alpha})^2\mathbf{C}^++\mathbf{C}^-)\Bv \,\mathrm{d}\boldsymbol{x}, \\
&\Kkk=\int_\Omega (1-\Ns\boldsymbol{\alpha})^2(3G\Ns^\mathrm{T}\Ns+\Bs^\mathrm{T}\etap^2\Bs)\,\mathrm{d}\boldsymbol{x}, \\
&\Kaa=\int_\Omega \bigg(\Ns^\mathrm{T}\Ns\Big(2\uppsi^{\mathrm{e}_+}+2\sigmap\Ns\boldsymbol{\kappa}+\etap^2\big\Vert\Bs \boldsymbol{\kappa}\big\Vert^2+w''(\Ns\boldsymbol{\alpha})\Big)+\Bs^\mathrm{T}\etad^2\Bs\bigg)\,\mathrm{d}\boldsymbol{x}, \\
&\Kuk=\Kku ^\mathrm{T}=\int_\Omega -\sqrt{\frac{3}{2}}\Bv^\mathrm{T}((1-\Ns\boldsymbol{\alpha})^2\mathbf{C}^++\mathbf{C}^-)\nhatu^\mathrm{tr}\Ns \,\mathrm{d}\boldsymbol{x}, \\
&\Kua=\Kau^\mathrm{T}=\int_\Omega -2(1-\Ns\boldsymbol{\alpha})\Bv^\mathrm{T}\mathbf{C}^+(\Bv\buu-\bepsilonp)\Ns\,\mathrm{d}\boldsymbol{x}, \\
&\Kak=\Kka ^\mathrm{T}=\int_\Omega-2(1-\Ns\boldsymbol{\alpha})\Ns^\mathrm{T}\bigg(-\sqrt{\frac{3}{2}}\Ns(\Bv\buu-\bepsilonp)^\mathrm{T}\mathbf{C}^+\nhatu^\mathrm{tr}+\sigmap\Ns+\boldsymbol{\kappa}^\mathrm{T}\Bs^\mathrm{T}\etap^2\Bs\bigg) \,\mathrm{d}\boldsymbol{x}, \\
&\Kklp=\Kald=-\Kalat=\int_\Omega -\Ns^\mathrm{T}\Ns\,\mathrm{d}\boldsymbol{x},
\end{alignedat}
\end{equation}
where the direction of the plastic flow is considered to be independent of the displacement field, and the function $w''(\Ns\boldsymbol{\alpha})$ equals zero if the AT-1 model is used and $2w_0$ if the AT-2 model is used. In addition, the KKT conditions~\eqref{eq:kktconditions} are discretized as
\begin{equation}
\begin{dcases}
\mathbf{r}_{\lp}&\hspace*{-3mm}\coloneqq\Lp(\boldsymbol{\kappa}-\boldsymbol{\kappa}_n)=\mathbf{0},\hspace{5.63mm}\blp\geq\mathbf{0},\hspace{4.95mm}(\boldsymbol{\kappa}-\boldsymbol{\kappa}_n)\geq\mathbf{0},\\
\mathbf{r}_{\ld}&\hspace*{-3mm}\coloneqq\Ld(\boldsymbol{\alpha}-\boldsymbol{\alpha}_n)=\mathbf{0},\hspace{3.4mm}\bld\geq\mathbf{0},\hspace{3.4mm}(\boldsymbol{\alpha}-\boldsymbol{\alpha}_n)\geq\mathbf{0},\\
\mathbf{r}_{\latt}&\hspace*{-3mm}\coloneqq\Lat(\mathbf{e}-\boldsymbol{\alpha})=\mathbf{0},\hspace{6.13mm}\blat\geq\mathbf{0},\hspace{3.4mm}(\mathbf{e}-\boldsymbol{\alpha})\geq\mathbf{0},
\end{dcases}
\label{eq:kktconditionsdisc}
\end{equation}
where $\mathbf{e}\in\Real^{\nnode}$ is a vector of ones, and $\Lp$, $\Ld$, $\Lat\in\Real^{\nnode\times\nnode}$ are diagonal matrices, whose diagonal entries are given by the components of $\blp$, $\bld$, and $\blat$, respectively. As such, Eqs.~\eqref{eq:weakformdisc} and~\eqref{eq:kktconditionsdisc} represent the set of discretized governing equations from which the primal variables $\{\buu,\boldsymbol{\kappa},\boldsymbol{\alpha}\}$ and the dual variables $\{\blp,\bld,\blat\}$ can be computed.

\section{Interior-point method} \label{sec:interiorpoint}
This section presents two algorithms for the solution of the constrained system of Eqs.~\eqref{eq:weakformdisc} and~\eqref{eq:kktconditionsdisc}, where an interior-point method is used to rigorously impose the constraints on the equivalent plastic strain and the crack phase-field. First, a monolithic scheme where all equations are solved simultaneously is elaborated. Second, a staggered scheme based on alternate minimization and an algorithmic decoupling of the equations is presented.

In both schemes, a perturbation is applied to the KKT conditions~\eqref{eq:kktconditionsdisc}, which reads
\begin{subnumcases}{\label{eq:perturbedkkt}}
\mathbf{r}_{\lp}\hspace*{1.5mm}\approx\Lp(\boldsymbol{\kappa}-\boldsymbol{\kappa}_n)-\mu_\kappa^{(k)}\mathbf{e}=\mathbf{0},\hspace{5.6mm}\blp\geq\mathbf{0},\hspace{5mm}(\boldsymbol{\kappa}-\boldsymbol{\kappa}_n)\geq\mathbf{0}, \label{eq:perturbedplastickkt} \\
\mathbf{r}_{\ld}\approx\Ld(\boldsymbol{\alpha}-\boldsymbol{\alpha}_n)-\mu_\alpha^{(k)}\mathbf{e}=\mathbf{0},\hspace{3.4mm}\bld\geq\mathbf{0},\hspace{3.4mm}(\boldsymbol{\alpha}-\boldsymbol{\alpha}_n)\geq\mathbf{0}, \label{eq:perturbeddamagekkt1} \\
\mathbf{r}_{\latt}\approx\Lat(\mathbf{e}-\boldsymbol{\alpha})-\mu_\alpha^{(k)}\mathbf{e}=\mathbf{0},\hspace{6.27mm}\blat\geq\mathbf{0},\hspace{3.4mm}(\mathbf{e}-\boldsymbol{\alpha})\geq\mathbf{0}, \label{eq:perturbeddamagekkt2}
\end{subnumcases}
where $\mu_\kappa^{(k)}\in\Real_+$ and $\mu_\alpha^{(k)}\in\Real_+$ are barrier parameters, with $k\in\mathbb{N}\,\cup\,\{0\}$. The perturbed problem is then solved for a sequence of barrier parameters $\big\{\mu_\kappa^{(k)}\big\}$ and $\big\{\mu_\alpha^{(k)}\big\}$ that converges to zero, in order to recover the original KKT conditions~\eqref{eq:kktconditionsdisc}.

\subsection{Monolithic scheme}\label{sec:monolithicscheme}
First, the primal-dual system is obtained through linearization of Eqs.~\eqref{eq:weakformdisc} and~\eqref{eq:perturbedkkt}, and later condensed, allowing an efficient computation of an undamped Newton step. Second, a damped Newton step is computed, taking into account the inequality constraints in the KKT conditions~\eqref{eq:perturbedkkt}. In the present method, only one iteration is performed for a given barrier parameter, as commonly done with interior-point methods~\cite{noce06a}. As such, the counter $k$ of the sequence of barrier parameters $\big\{\mu_\kappa^{(k)}\big\}$ and $\big\{\mu_\alpha^{(k)}\big\}$ coincides with the iteration counter of the approximated primal and dual variables. After each iteration, the barrier parameter is updated and the linearization is repeated.

A convergence criterion based on three residuals will be adopted, i.e.\ $\max\{{\tt RES}_{\bu},{\tt RES}_{\kappa},{\tt RES}_{\alpha}\}\leq{\tt TOL}$, where ${\tt TOL}$ is a user-defined tolerance. In the literature, different definitions of the residuals have been presented~\cite{amor09a,amba15b,ales18b,rodr18a,gera19a}. Hereafter, the infinity norm of the optimality conditions ${\tt RES}_{\bu}\coloneqq\Vert \mathbf{r}_{\bu} \Vert_\infty$, ${\tt RES}_{\kappa}\coloneqq\max\{\Vert \mathbf{r}_{\kappa} \Vert_\infty, \Vert \mathbf{r}_{\lp} \Vert_\infty\}$, and $ {\tt RES}_{\alpha}\coloneqq\max\{\Vert \mathbf{r}_{\alpha} \Vert_\infty, \Vert \mathbf{r}_{\ld} \Vert_\infty, \Vert \mathbf{r}_{\latt} \Vert_\infty\}$ is used to check convergence. The residuals are normalized through a division by their corresponding units, rendering their value unitless.

\subsubsection*{Primal-dual system}
Linearization of Eqs.~\eqref{eq:weakformdisc}~and~\eqref{eq:perturbedkkt} leads to the following primal-dual system:
\begin{equation}
\begin{bmatrix}
\Kuu & \Kuk & \Kua & \mathbf{0} & \mathbf{0} & \mathbf{0} \\
\Kku & \Kkk & \Kka & \Kklp & \mathbf{0} & \mathbf{0} \\
\Kau & \Kak & \Kaa & \mathbf{0} & \Kald & \Kalat \\
\mathbf{0} & \buL_\kappa & \mathbf{0} & \Sp & \mathbf{0} & \mathbf{0} \\
\mathbf{0} & \mathbf{0} & \phantom{-}\buL_{\alpha_1} & \mathbf{0} & \Sd & \mathbf{0} \\
\mathbf{0} & \mathbf{0} & -\buL_{\alpha_2} & \mathbf{0} & \mathbf{0} & \Sat \\
\end{bmatrix}^{(k)}
\begin{bmatrix}
\Delta\buu \\
\Delta\boldsymbol{\kappa} \\
\Delta\boldsymbol{\alpha} \\
\Delta\bul_\kappa \\
\Delta\bul_{\alpha_1} \\
\Delta\bul_{\alpha_2} \\
\end{bmatrix}
=-
\begin{bmatrix}
\mathbf{r}_{\bu} \\
\mathbf{r}_{\kappa} \\
\mathbf{r}_{\alpha} \\
\mathbf{r}_{\lambda_\kappa} \\
\mathbf{r}_{\lambda_{\alpha_1}} \\
\mathbf{r}_{\lambda_{\alpha_2}}
\end{bmatrix}^{(k)},
\label{eq:primaldual}
\end{equation}
where $\Sp$, $\Sd$, $\Sat\in\Real^{\nnode\times\nnode}$ are diagonal matrices whose diagonal entries are given by the components of the inequality constraints $\boldsymbol{\kappa}^{(k)}-\boldsymbol{\kappa}_n$, $\boldsymbol{\alpha}^{(k)}-\boldsymbol{\alpha}_n$ and $\mathbf{e}-\boldsymbol{\alpha}^{(k)}$, respectively. This system can be condensed by eliminating the dual variables, resulting in the following reduced expression in the primal variables:
\begin{equation}
\begin{bmatrix}
\Kuu & \Kuk & \Kua \\
\Kku & \Kkk^{'} & \Kka \\
\Kau & \Kak & \Kaa^{'}
\end{bmatrix}^{(k)}
\begin{bmatrix}
\Delta\buu \\
\Delta\boldsymbol{\kappa} \\
\Delta\boldsymbol{\alpha}
\end{bmatrix}
=-
\begin{bmatrix}
\mathbf{r}_{\bu} \\
\mathbf{r}_{\kappa}^{'} \\
\mathbf{r}_{\alpha}^{'}
\end{bmatrix}^{(k)},
\label{eq:condensedprimal}
\end{equation}
with
\begin{equation}
\begin{alignedat}{1}
\Kkk^{'}&\coloneqq\Kkk^{(k)}-\Kklp^{(k)}\Lp^{(k)}\Sp^{-1}, \\
\Kaa^{'}&\coloneqq\Kaa^{(k)}-\Kald^{(k)}\Ld^{(k)}\Sd^{-1}+\Kalat^{(k)}\Lat^{(k)}\Sat^{-1}, \\
\mathbf{r}_{\kappa}^{'}&\coloneqq\mathbf{r}_{\kappa}^{(k)}-\Kklp^{(k)}(\blp^{(k)}-\mu_\kappa^{(k)}\Sp^{-1}\mathbf{e}), \\
\mathbf{r}_{\alpha}^{'}&\coloneqq\mathbf{r}_{\alpha}^{(k)}-\Kald^{(k)}(\bld^{(k)}-\mu_\alpha^{(k)}\Sd^{-1}\mathbf{e})-\Kalat^{(k)}(\blat^{(k)}-\mu_\alpha^{(k)}\Sat^{-1}\mathbf{e}).
\end{alignedat}
\end{equation}
The dual variables can then be obtained as
\begin{subnumcases}{\label{eq:condenseddual}}
\Delta\blp\hspace*{1.55mm}=-\blp^{(k)}-\Sp^{-1}(\Lp^{(k)}\Delta\boldsymbol{\kappa}-\mu_\kappa^{(k)}\mathbf{e}), \label{eq:condensedlp} \\
\Delta\bld=-\bld^{(k)}-\Sd^{-1}(\Ld^{(k)}\Delta\boldsymbol{\alpha}-\mu_\alpha^{(k)}\mathbf{e}), \label{eq:condensedld} \\
\Delta\blat=-\blat^{(k)}-\Sat^{-1}(-\Lat^{(k)}\Delta\boldsymbol{\alpha}-\mu_\alpha^{(k)}\mathbf{e}) \label{eq:condensedlat},
\end{subnumcases}
which can be efficiently computed owing to the diagonal structure of $\Sp$, $\Sd$ and $\Sat$. Similar to Krabbenh{\o}ft et~al.~\cite{krab07b}, a selective integration rule is used, where the terms involving the dual variables are integrated using quadrature points that coincide with the nodes. As such, $\Kklp$, $\Kald$ and $\Kalat$ become diagonal and the system matrix in Eq.~\eqref{eq:condensedprimal} becomes symmetric.

\subsubsection*{Damped Newton step}
A damped Newton step is performed with step lengths $\tp\in\Real_+$ and $\td\in\Real_+$ such that the inequality constraints in the KKT conditions~\eqref{eq:perturbedkkt} hold for all iterations. The values of the new iterate are then given by
\begin{equation}
\begin{alignedat}{1}
\big\{\buu^{(k+1)},\boldsymbol{\kappa}^{(k+1)},\boldsymbol{\alpha}^{(k+1)}\big\}&\coloneqq\big\{\buu^{(k)},\boldsymbol{\kappa}^{(k)},\boldsymbol{\alpha}^{(k)}\big\}+\tp\big\{\Delta\buu,\Delta\boldsymbol{\kappa},\Delta\boldsymbol{\alpha}\big\}, \\
\big\{\blp^{(k+1)},\bld^{(k+1)},\blat^{(k+1)}\big\}&\coloneqq\big\{\blp^{(k)},\bld^{(k)},\blat^{(k)}\big\}+\td\big\{\Delta\blp,\Delta\bld,\Delta\blat\big\},
\end{alignedat}
\label{eq:update}
\end{equation}
where the step lengths $\tp$ and $\td$ are computed as
\begin{equation}
\begin{alignedat}{1}
\tp&\coloneqq\text{max}\left\{t\in(0,1]\,\vert\begin{bmatrix}\boldsymbol{\kappa}^{(k)}-\boldsymbol{\kappa}_n\\ \boldsymbol{\alpha}^{(k)}-\boldsymbol{\alpha}_n \\ \mathbf{e}-\boldsymbol{\alpha}^{(k)}\end{bmatrix}+t\begin{bmatrix}\Delta\boldsymbol{\kappa} \\ \Delta\boldsymbol{\alpha} \\ -\Delta\boldsymbol{\alpha}\end{bmatrix}\geq(1-\tau)\begin{bmatrix}\boldsymbol{\kappa}^{(k)}-\boldsymbol{\kappa}_n \\ \boldsymbol{\alpha}^{(k)}-\boldsymbol{\alpha}_n \\ \mathbf{e}-\boldsymbol{\alpha}^{(k)}\end{bmatrix}\right\}, \\
\td&\coloneqq\text{max}\left\{t\in(0,1]\,\vert\begin{bmatrix}\blp^{(k)} \\ \bld^{(k)} \\ \blat^{(k)}\end{bmatrix}+t\begin{bmatrix}\Delta\blp \\ \Delta\bld \\ \Delta\blat\end{bmatrix}\geq(1-\tau)\begin{bmatrix}\blp^{(k)} \\ \bld^{(k)} \\ \blat^{(k)}\end{bmatrix}\right\}.
\end{alignedat}
\label{eq:stepsize}
\end{equation}
In Eq.~\eqref{eq:stepsize}, $\tau\in(0,1)$ is a constant which prevents the inequality constraints and the dual variables from approaching their lower bound of zero too quickly~\cite{noce06a}. A too small value of $\tau$ prevents the algorithm from taking large steps, while a too large value of $\tau$ allows steps that are too close to the lower bounds of the constrained fields. Hereinafter, a value of 0.999 is assigned to~$\tau$.

\subsubsection*{Computation of the barrier parameter} 
The sequence of barrier parameters $\big\{\mu_\kappa^{(k)}\big\}$ and $\big\{\mu_\alpha^{(k)}\big\}$ in Eq.~\eqref{eq:perturbedkkt} must converge to zero to recover the original expression for the KKT conditions~\eqref{eq:kktconditionsdisc}. Because a single system is solved in the monolithic scheme, only one barrier parameter $\mu^{(k)}\in\Real_+$ is considered in the perturbation of the KKT conditions, i.e.~$\mu_\kappa^{(k)}=\mu_\alpha^{(k)}=\mu^{(k)}$. Similar to Nocedal and Wright~\cite{noce06a}, the barrier parameter $\mu^{(k)}$ is expressed as the product of a centering parameter and a duality measure, which vanishes once the KKT conditions~\eqref{eq:kktconditionsdisc} hold:
\begin{equation}
\mu^{(k)}:=\sigma^{(k)}\frac{(\boldsymbol{\kappa}^{(k)}-\boldsymbol{\kappa}_n)^\mathrm{T}\blp^{(k)}+(\boldsymbol{\alpha}^{(k)}-\boldsymbol{\alpha}_n)^\mathrm{T}\bld^{(k)}+(\mathbf{e}-\boldsymbol{\alpha}^{(k)})^\mathrm{T}\blat^{(k)}}{n_\text{const}},
\label{eq:barrierparameter}
\end{equation}
where $\sigma^{(k)}\in\Real_+$ is the centering parameter which corresponds to the desired reduction factor of the duality measure, and $n_\text{const}$ is the total number of nodal constraints. In the literature, many strategies are presented to update the centering parameter~\cite{fiac90a,mehr92a,vand99a,noce06a,wach06a}. In this paper, an adaptive strategy based on a predictor direction is employed, similar to the primal-dual interior-point method of Mehrotra~\cite{mehr92a}. This technique requires solving the primal-dual system~\eqref{eq:primaldual} with $\mu^{(k)}=0$ to obtain the affine-scaling direction $\{\Delta\boldsymbol{\kappa}^\mathrm{aff},\Delta\boldsymbol{\alpha}^\mathrm{aff},\Delta\blp^\mathrm{aff},\Delta\bld^\mathrm{aff},\Delta\blat^\mathrm{aff}\}$. The centering parameter is then computed as
\begin{equation}
\sigma^{(k)}\coloneqq\Bigg[\frac{(\boldsymbol{\kappa}^\mathrm{aff}-\boldsymbol{\kappa}_n)^\mathrm{T}\blp^\mathrm{aff}+(\boldsymbol{\alpha}^\mathrm{aff}-\boldsymbol{\alpha}_n)^\mathrm{T}\bld^\mathrm{aff}+(\mathbf{e}-\boldsymbol{\alpha}^\mathrm{aff})^\mathrm{T}\blat^\mathrm{aff}}{(\boldsymbol{\kappa}^{(k)}-\boldsymbol{\kappa}_n)^\mathrm{T}\blp^{(k)}+(\boldsymbol{\alpha}^{(k)}-\boldsymbol{\alpha}_n)^\mathrm{T}\bld^{(k)}+(\mathbf{e}-\boldsymbol{\alpha}^{(k)})^\mathrm{T}\blat^{(k)}}\Bigg]^3,
\label{eq:centeringparameter}
\end{equation}
where $\boldsymbol{\alpha}^\mathrm{aff}$, $\boldsymbol{\kappa}^\mathrm{aff}$, $\blp^\mathrm{aff}$, $\bld^\mathrm{aff}$ and $\blat^\mathrm{aff}$ are computed with Eqs.~\eqref{eq:update} and~\eqref{eq:stepsize} using the affine-scaling direction. This choice of the centering parameter allows to reduce the barrier parameter in case large affine-scaling steps can be taken without violating the inequality constraints. Conversely, when not much progress can be made along the affine-scaling direction before reaching the boundary of the inequality constraints, the barrier parameter increases. This heuristic choice works well in practice, but has no real analytical justification~\cite{noce06a}.

\subsubsection*{Starting values}
First, it is checked if any damage value of the previous time step exceeds a threshold $\tt{CRTOL}$. If so, those variables are fixed to 1 in the current time step, and both the monotonicity constraint and the upper bound constraint for that variable are removed from the primal-dual system. Hereafter, a value of 0.999 is assigned to $\tt{CRTOL}$. Second, the internal variables which are not fixed are initialized at the values of the previous time step plus a small increment, such that all inequalities hold strictly (e.g.\ $\boldsymbol{\alpha}^{(0)}\coloneqq\boldsymbol{\alpha}_n+(1-\tt{CRTOL})\mathbf{e}$ and $\boldsymbol{\kappa}^{(0)}\coloneqq\boldsymbol{\kappa}_n+(1-\tt{CRTOL})\mathbf{e}$). Third, the initial estimate $\buu^{(0)}$ of the displacement field is computed by freezing all the internal variables and taking an elastic step (i.e.~solving the first equation in~\eqref{eq:weakformdisc} for fixed internal variables). Finally, the dual variables are initialized by computing the affine-scaling direction $\{\Delta\blp^{\mathrm{aff}(0)},\Delta\bld^{\mathrm{aff}(0)},\Delta\blat^{\mathrm{aff}(0)}\}$ and computing each component $\lambda^{(0)}\in\{\blp^{(0)},\bld^{(0)},\blat^{(0)}\}$ as $\max\{1,\vert1+\Delta\lambda^\mathrm{aff(0)}\vert\}$ as proposed by Gertz~et~al.~\cite{gert04a}.

\subsubsection*{Stabilization of the monolithic scheme}
To overcome divergence issues resulting from the non-convexity of the energy functional, a heuristic stabilization is applied to the system matrix in Eq.~\eqref{eq:condensedprimal}. Since the convergence issues are attributed to the coupled equilibrium-damage system, the stabilization is determined for fixed plastic variables. The stabilization adjusts the inertia of the Hessian matrix of the functional for fixed plastic variables by adding a constant $\delta^{(k)}\in\Real_+$ to its diagonal, similar to W{\"a}chter and Biegler~\cite{wach06a}. The inertia of a matrix $\mathbf{K}$ is the scalar triplet $\{n_+,n_-,n_0\}$ that indicates the number of positive, negative, and zero eigenvalues (i.e. $\text{inertia}(\mathbf{K}):=\{n_+,n_-,n_0\}$). As such, the modified Hessian matrix $\bar{\mathbf{K}}^{(k)}\in\Real^{(\ndof+\nnode)\times(\ndof+\nnode)}$ of the equilibrium-damage problem reads
\begin{equation}
\bar{\mathbf{K}}^{(k)}\coloneqq
\begin{bmatrix}
\Kuu & \Kua \\
\Kau & \Kaa
\end{bmatrix}^{(k)}
\hspace{-0.8mm}
+\delta^{(k)}\mathbf{I},
\end{equation}
where $\mathbf{I}$ is the identity matrix.

The stabilization factor should be large enough for all eigenvalues of $\bar{\mathbf{K}}^{(k)}$ to be strictly positive, but as small as possible to avoid large modifications of the original matrix~\cite{noce06a}. Because the magnitude of the stabilization factor is not known in advance, successively larger values are tested until the desired inertia is obtained. Algorithm~\ref{alg:inertiacorrection} shows how the stabilization factor is computed. This algorithm is adapted from Nocedal and Wright~\cite{noce06a} and aims to avoid unnecessarily large modifications of the original matrix, while trying to minimize the number of matrix factorizations. The required number of factorizations can be altered by changing the starting value of the stabilization factor or the multiplication factors. In algorithm~\ref{alg:inertiacorrection}, the values of Nocedal and Wright~\cite{noce06a} are adopted. The inertia of the modified matrix $\bar{\mathbf{K}}^{(k)}$ is obtained from its symmetric indefinite factorization computed using the MA57 routine in the Harwell Subroutine Library~\cite{hsl}.
\begin{algorithm}
\caption{Stabilization of the monolithic scheme}
\begin{algorithmic}
\label{alg:inertiacorrection}
\REQUIRE $\delta^{(k-1)}$ from the previous interior-point iteration
\vspace{-0.3em}
\STATE $\delta^{(k)} \leftarrow 0$
\vspace{-0.3em}
\IF{inertia($\bar{\mathbf{K}}^{(k)}$) $\neq$ $\{\ndof+\nnode,0,0\}$}
\vspace{-0.3em}
\IF{$\delta^{(k-1)}=0$}
\vspace{-0.3em}
\STATE $\delta^{(k)}\leftarrow 10^{-4}$
\vspace{-0.3em}
\ELSE
\vspace{-0.3em}
\STATE $\delta^{(k)} \leftarrow \delta^{(k-1)}/2$
\vspace{-0.3em}
\ENDIF
\vspace{-0.3em}
\WHILE{inertia($\bar{\mathbf{K}}^{(k)}$) $\neq$ $\{\ndof+\nnode,0,0\}$}
\vspace{-0.3em}
\STATE $\delta^{(k)} \leftarrow 10\delta^{(k)}$
\vspace{-0.3em}
\ENDWHILE
\vspace{-0.3em}
\ENDIF
\vspace{-0.3em}
\RETURN $\delta^{(k)}$
\end{algorithmic}
\end{algorithm}

\vspace{-1em}
Once the stabilization factor is determined, the undamped step of primal variables $\{\Delta\buu,\Delta\boldsymbol{\kappa},\Delta\boldsymbol{\alpha}\}$ is computed from a perturbed version of the condensed primal-dual system~\eqref{eq:condensedprimal}:
\begin{equation}
\begin{bmatrix}
\Kuu+\delta\mathbf{I} & \Kuk & \Kua \\
\Kku & \Kkk^{'} & \Kka \\
\Kau & \Kak & \Kaa^{'}+\delta\mathbf{I}
\end{bmatrix}^{(k)}
\begin{bmatrix}
\Delta\buu \\
\Delta\boldsymbol{\kappa} \\
\Delta\boldsymbol{\alpha}
\end{bmatrix}
=-
\begin{bmatrix}
\mathbf{r}_{\bu} \\
\mathbf{r}_{\kappa}^{'} \\
\mathbf{r}_{\alpha}^{'}
\end{bmatrix}^{(k)},
\label{eq:perturbedprimal}
\end{equation}
of which the system matrix is symmetric and can be factorized using a symmetric indefinite factorization.

\vspace{-0.5em}
\subsubsection*{Overview}
Algorithm~\ref{alg:interiorpoint} presents an overview of the interior-point method used in the monolithic scheme.
\begin{algorithm}[!h]
\caption{Monolithic interior-point method}
\begin{algorithmic}
\label{alg:interiorpoint}
\REQUIRE $\bu_n$, $\boldsymbol{\kappa}_n$, $\boldsymbol{\alpha}_n$.
\vspace{-0.3em}
\STATE Initialize iteration counter $k\coloneqq0$ and starting values $\{\bu^{(0)},\boldsymbol{\kappa}^{(0)},\boldsymbol{\alpha}^{(0)},\blp^{(0)},\bld^{(0)},\blat^{(0)}\}$.
\vspace{-0.3em}
\STATE Initialize the barrier parameter $\mu^{(0)}$ (Eqs.~\eqref{eq:barrierparameter} and~\eqref{eq:centeringparameter}).
\vspace{-0.3em}
\STATE Initialize the stabilization factor $\delta^{(0)}$ (algorithm~\ref{alg:inertiacorrection}).\vspace{-0.2em}
\REPEAT
\vspace{-0.3em}
\STATE Set $k\leftarrow k+1$.
\vspace{-0.3em}
\STATE 1) Compute the undamped Newton step $\{\Delta\bu,\Delta\boldsymbol{\kappa},\Delta\boldsymbol{\alpha}\}$ and $\{\Delta\blp,\Delta\bld,\Delta\blat\}$ (Eqs.~\eqref{eq:condenseddual} and~\eqref{eq:perturbedprimal}).\hspace{-1mm}
\vspace{-0.3em}
\STATE 2) Compute the primal and dual variables $\{\bu^{(k)},\boldsymbol{\kappa}^{(k)},\boldsymbol{\alpha}^{(k)}\}$ and $\{\blp^{(k)},\bld^{(k)},\blat^{(k)}\}$ (Eqs.~\eqref{eq:update} and~\eqref{eq:stepsize}).\hspace{-1mm}
\vspace{-0.3em}
\STATE 3) Compute the barrier parameter $\mu^{(k)}$ (Eqs.~\eqref{eq:barrierparameter} and~\eqref{eq:centeringparameter}).
\vspace{-0.3em}
\STATE 4) Compute the stabilization factor $\delta^{(k)}$ (algorithm~\ref{alg:inertiacorrection}).
\vspace{-0.3em}
\STATE 5) Compute ${\tt RES}_{\bu}^{(k)}$, ${\tt RES}_{\kappa}^{(k)}$ and ${\tt RES}_{\alpha}^{(k)}$.
\vspace{-0.3em}
\UNTIL $\max\{{\tt RES}^{(k)}_{\bu},{\tt RES}^{(k)}_{\kappa},{\tt RES}^{(k)}_{\alpha}\}\leq{\tt TOL}$.
\vspace{-0.3em}
\STATE Set $\{\bu,\boldsymbol{\kappa},\boldsymbol{\alpha}\} \leftarrow \{\bu^{(k)},\boldsymbol{\kappa}^{(k)},\boldsymbol{\alpha}^{(k)}\}$.
\vspace{-0.3em}
\RETURN $\bu$, $\boldsymbol{\kappa}$, $\boldsymbol{\alpha}$.
\end{algorithmic}
\end{algorithm}

\vspace{-0.5em}
\subsection{Staggered scheme}
An alternative staggered scheme, which algorithmically decouples the governing equations and does not require further stabilization will now be presented. First, the alternate minimization algorithm is presented. This algorithm solves a sequence of convex sub-problems of which the optimality conditions correspond to Eqs.~\eqref{eq:weakformdisc} and~\eqref{eq:kktconditionsdisc}. The convex sub-problems are obtained by minimizing the energy functional~\eqref{eq:functional} with respect to the displacement field, plastic strain tensor and equivalent plastic strain, and crack phase-field separately. Second, the solution of the non-linear balance equation is elaborated. In a third section, an interior-point method for the rigorous solution of the constrained equations governing the rate-independent non-local plastic evolution problem is presented. The proposed method provides an alternative to viscous regularization~\cite{mieh16a,mieh16b} or the a posteriori projection technique~\cite{ullo16a,rodr18a}. In the case of local plasticity, return-mapping algorithms~\cite{amba15a,duda15a,amba16a,kuhn16a,bord16a} can be adopted to rigorously obtain a solution. Finally, in the fourth section, an interior-point method is elaborated for the rigorous solution of the damage evolution problem without modifications of the governing equations or the need for a penalty parameter. The damage evolution problem can be alternatively solved by setting Dirichlet-type conditions on the crack phase-field~\cite{bour07a}, penalization~\cite{mieh10a,whee14a,gera19a}, the a posteriori projection technique~\cite{lanc09a,ullo16a,rodr18a}, active-set strategies~\cite{heis15a,farr17a}, or the use of a history field~\cite{mieh10b}.

\vspace{-0.6em}
\subsubsection{Alternate minimization}
The use of alternate minimization was originally proposed in the context of brittle fracture by Bourdin~et~al.~\cite{bour00a}. Later, several authors have extended this approach to ductile fracture~\cite{ales14a,amba15a,amba16a,rodr18a}. The alternate minimization algorithm solves the constrained system of Eqs.~\eqref{eq:weakformdisc} and~\eqref{eq:kktconditionsdisc} by fixing two of the three fields and solving for the other field. The free field is then updated, and the process is repeated for the other fields until convergence is achieved. Algorithm~\ref{alg:alternateminimization} provides an overview of the alternate minimization algorithm.

Contrary to the monolithic scheme presented in section~\ref{sec:monolithicscheme}, the starting values of each sub-problem are based on the results of the previous staggered iteration, i.e. $\bu^{(i-1)}$, $\boldsymbol{\kappa}^{(i-1)}+(1-\tt{CRTOL})\mathbf{e}$ and $\boldsymbol{\alpha}^{(i-1)}+(1-\tt{CRTOL})\mathbf{e}$, where $i\in\mathbb{N}\,\cup\,\{0\}$ corresponds to the $i$-th staggered iteration. In case the AT-1 model is adopted, this might lead to violations of the upper bound of the crack phase-field and a starting value of $\boldsymbol{\alpha}_n+(1-\tt{CRTOL})\mathbf{e}$ is used instead. For the initial staggered iteration, the values of the previous time step are used.

The non-linear equations in algorithm~\ref{alg:alternateminimization} are solved using the iterative procedures described in the following sections. For brevity, the superscript $\square^{(i)}$ of the staggered iteration loop is omitted, and a superscript $\square^{(k)}$ is introduced, corresponding to the $k$-th approximation of $\square^{(i)}$.

\begin{algorithm}[h!]
\caption{Alternate minimization}
\begin{algorithmic}
\label{alg:alternateminimization}
\REQUIRE $\bu_n$, $\boldsymbol{\kappa}_n$, $\boldsymbol{\alpha}_n$.

\STATE Initialize iteration counter $i\coloneqq0$ and set $\{\bu^{(0)},\boldsymbol{\kappa}^{(0)},\boldsymbol{\alpha}^{(0)}\}\coloneqq\{\bu_n,\boldsymbol{\kappa}_n,\boldsymbol{\alpha}_n\}$.

\REPEAT

\STATE Set $i\leftarrow i+1$.

\STATE 1) Solve the \underline{balance equation} (section~\ref{sec:stagbalance})
\begin{equation}
\mathbf{r}_{\bu}\big(\buu^{(i)},\bepsilonp(\boldsymbol{\kappa}^{(i-1)}),\boldsymbol{\alpha}^{(i-1)}\big)=\mathbf{0}
\label{eq:stagbalance}
\end{equation}
\hspace{4mm}for $\buu^{(i)}$ such that ${\tt RES}_{\bu}^{(i)}\leq\tt{TOL}$.

\STATE 2) Solve the \underline{plastic evolution problem} (section~\ref{sec:stagplastic})
\begin{subnumcases}{\label{eq:stagplastic}}
\mathbf{r}_{\kappa}\big(\buu^{(i)},\boldsymbol{\kappa}^{(i)},\boldsymbol{\alpha}^{(i-1)},\blp^{(i)}\big)=\mathbf{0}, \label{eq:stagplasticeq} \\
\mathbf{r}_{\lp}\big(\boldsymbol{\kappa}^{(i)},\blp^{(i)}\big)=\mathbf{0},\quad\blp^{(i)}\geq\mathbf{0},\quad\big(\boldsymbol{\kappa}^{(i)}-\boldsymbol{\kappa}_n\big)\geq\mathbf{0} \label{eq:stagplastickkt}
\end{subnumcases}
\hspace{4mm}for $\{\boldsymbol{\kappa}^{(i)},\blp^{(i)}\}$ such that ${\tt RES}_{\kappa}^{(i)}\leq\tt{TOL}$.

\STATE 3) Solve the \underline{damage evolution problem} (section~\ref{sec:stagdamage})
\begin{subnumcases}{\label{eq:stagdamage}}
\mathbf{r}_{\alpha}\big(\buu^{(i)},\bepsilonp(\boldsymbol{\kappa}^{(i)}),\boldsymbol{\kappa}^{(i)},\boldsymbol{\alpha}^{(i)},\bld^{(i)},\blat^{(i)}\big)=\mathbf{0}, \label{eq:stagdamageeq} \\
\mathbf{r}_{\ld}\big(\boldsymbol{\alpha}^{(i)},\bld^{(i)}\big)=\mathbf{0},\quad\bld^{(i)}\geq\mathbf{0},\quad\big(\boldsymbol{\alpha}^{(i)}-\boldsymbol{\alpha}_n\big)\geq\mathbf{0}, \label{eq:stagdamagekkt1} \\
\mathbf{r}_{\latt}\big(\boldsymbol{\alpha}^{(i)},\blat^{(i)}\big)=\mathbf{0},\quad\blat^{(i)}\geq\mathbf{0},\quad\big(\mathbf{e}-\boldsymbol{\alpha}^{(i)}\big)\geq\mathbf{0} \label{eq:stagdamagekkt2}
\end{subnumcases}
\hspace{4mm}for $\{\boldsymbol{\alpha}^{(i)},\bld^{(i)},\blat^{(i)}\}$ such that ${\tt RES}_{\alpha}^{(i)}\leq\tt{TOL}$.

\STATE 4) Update ${\tt RES}_{\bu}^{(i)}$ and ${\tt RES}_{\kappa}^{(i)}$ using $\{\bu^{(i)},\boldsymbol{\kappa}^{(i)},\boldsymbol{\alpha}^{(i)},\blp^{(i)}\}$.

\UNTIL ${\tt RES}_{\bu}^{(i)}\leq{\tt TOL}$ and ${\tt RES}_{\kappa}^{(i)}\leq{\tt TOL}$.

\STATE Set $\{\bu,\boldsymbol{\kappa},\boldsymbol{\alpha}\} \leftarrow \{\bu^{(i)},\boldsymbol{\kappa}^{(i)},\boldsymbol{\alpha}^{(i)}\}$.

\RETURN $\bu$, $\boldsymbol{\kappa}$, $\boldsymbol{\alpha}$.
\end{algorithmic}
\end{algorithm}

\subsubsection{Balance equation} \label{sec:stagbalance}
Due to the decomposition of the elastic energy density in Eq.~\eqref{eq:split}, the balance equation (Eq.~\eqref{eq:stagbalance}) becomes non-linear. Similar to Amor~et~al.~\cite{amor09a}, an undamped Newton-Raphson procedure is used to solve it. The linearized expression reads
\begin{equation}
\Kuu^{(k)}\Delta\buu=-\mathbf{r}_{\bu}^{(k)}.
\end{equation}
The displacement field is then updated by performing an undamped Netwon step:
\begin{equation}
\buu^{(k+1)}\coloneqq\buu^{(k)}+\Delta\buu.
\end{equation}
Alternatively, a damped Newton step can be performed. However, the use of an undamped Newton step did not lead to convergence issues in the present method. In addition, in order to fairly compare the results of this approach to the monolithic scheme presented in section~\ref{sec:monolithicscheme}, where the damped method is only used to stay in the feasible region, only the undamped method to solve the balance equation will be considered hereafter.

\subsubsection{Plastic evolution problem} \label{sec:stagplastic}
Similar to the monolithic scheme presented in section~\ref{sec:monolithicscheme}, a perturbation is applied to the KKT conditions~\eqref{eq:stagplastickkt}. Since the displacement field and the crack phase-field are kept constant, only Eq.~\eqref{eq:stagplastic} needs to be solved. As such, algorithm~\ref{alg:interiorpoint} can be used to solve the plastic evolution problem, where only the equivalent plastic strain $\boldsymbol{\kappa}$ and the plastic dual Lagrange multiplier $\blp$ are solved for and the stabilization of the Hessian matrix is omitted. In addition, the barrier parameter $\mu_\kappa^{(k)}$ is computed using Eqs.~\eqref{eq:barrierparameter} and~\eqref{eq:centeringparameter} by setting all discretized variables but $\boldsymbol{\kappa}$, $\boldsymbol{\kappa}_n$, $\boldsymbol{\kappa}^\mathrm{aff}$, $\blp$ and $\blp^\mathrm{aff}$ to zero. For completeness, the primal-dual system and the damped Newton step are elaborated next.

\subsubsection*{Primal-dual system}
Linearization of Eq.~\eqref{eq:stagplastic} leads to the following primal-dual system:
\begin{equation}
\begin{bmatrix}
\Kkk & \Kklp \\
\buL_\kappa & \Sp \\
\end{bmatrix}^{(k)}
\begin{bmatrix}
\Delta\boldsymbol{\kappa} \\
\Delta\bul_\kappa \\
\end{bmatrix}
=-
\begin{bmatrix}
\mathbf{r}_{\kappa} \\
\mathbf{r}_{\lambda_\kappa} \\
\end{bmatrix}^{(k)}.
\end{equation}
Similar to the monolithic solution procedure presented in section~\ref{sec:monolithicscheme}, the dual variables are condensed out, and the primal variables and dual variables can be computed using $\Kkk^{'}\Delta\boldsymbol{\kappa}=-\mathbf{r}_{\kappa}^{'}$ and Eq.~\eqref{eq:condensedlp}, respectively.

\subsubsection*{Damped Newton step}
A damped Newton step is performed with step lengths $t^\text{prim}_\kappa\in\Real_+$ and $t^\text{dual}_\kappa\in\Real_+$ such that the inequality constraints in the KKT conditions~\eqref{eq:stagplastickkt} hold for all iterations. As such, the values of the new iterate are given by
\begin{equation}
\begin{alignedat}{1}
\boldsymbol{\kappa}^{(k+1)}&\coloneqq\boldsymbol{\kappa}^{(k)}+t^\text{prim}_\kappa\Delta\boldsymbol{\kappa}, \\
\blp^{(k+1)}&\coloneqq\blp^{(k)}+t^\text{dual}_\kappa\Delta\blp,
\end{alignedat}
\end{equation}
where the step lengths $t^\text{prim}_\kappa$ and $t^\text{dual}_\kappa$ are computed as
\begin{equation}
\begin{alignedat}{1}
t^\text{prim}_\kappa&\coloneqq\text{max}\big\{t\in(0,1]\,\vert\,(\boldsymbol{\kappa}^{(k)}-\boldsymbol{\kappa}_n)+t\Delta\boldsymbol{\kappa} \geq(1-\tau)(\boldsymbol{\kappa}^{(k)}-\boldsymbol{\kappa}_n)\big\}, \\
t^\text{dual}_\kappa&\coloneqq\text{max}\big\{t\in(0,1]\,\vert\,\blp^{(k)}+t\Delta\blp\geq(1-\tau)\blp^{(k)}\big\}.
\end{alignedat}
\end{equation}

\clearpage

\subsubsection{Damage evolution problem} \label{sec:stagdamage}
Similar to the monolithic scheme presented in section~\ref{sec:monolithicscheme} and the interior-point method for the solution of the plastic evolution problem presented in section~\ref{sec:stagplastic}, a perturbation is applied to the KKT conditions~\eqref{eq:stagdamagekkt1} and~\eqref{eq:stagdamagekkt2}. Since the displacement field and the plastic variables are kept constant, only Eq.~\eqref{eq:stagdamage} needs to be solved. As such, algorithm~\ref{alg:interiorpoint} can be used to solve the damage evolution problem, where only the crack phase-field $\boldsymbol{\alpha}$ and the damage dual Lagrange multipliers $\bld$ and $\blat$ are solved for and the stabilization of the Hessian matrix is omitted. In addition, the barrier parameter $\mu_\alpha^{(k)}$ is computed using Eqs.~\eqref{eq:barrierparameter} and~\eqref{eq:centeringparameter} by setting all variables but $\boldsymbol{\alpha}$, $\boldsymbol{\alpha}_n$, $\boldsymbol{\alpha}^\mathrm{aff}$, $\bld$, $\bld^\mathrm{aff}$, $\blat$ and $\blat^\mathrm{aff}$ to zero. Similar to section~\ref{sec:stagplastic}, the primal-dual system and the damped Newton step are elaborated next.

\subsubsection*{Primal-dual system}
Linearization of Eq.~\eqref{eq:stagdamage} leads to the following primal-dual system:
\begin{equation}
\begin{bmatrix}
\Kaa & \Kald & \Kalat \\
\phantom{-}\buL_{\alpha_1} & \Sd & \mathbf{0} \\
-\buL_{\alpha_2} & \mathbf{0} & \Sat \\
\end{bmatrix}^{(k)}
\begin{bmatrix}
\Delta\boldsymbol{\alpha} \\
\Delta\bul_{\alpha_1} \\
\Delta\bul_{\alpha_2} \\
\end{bmatrix}
=-
\begin{bmatrix}
\mathbf{r}_{\alpha} \\
\mathbf{r}_{\lambda_{\alpha_1}} \\
\mathbf{r}_{\lambda_{\alpha_2}}
\end{bmatrix}^{(k)}.
\end{equation}
Similar to the monolithic solution procedure presented in section~\ref{sec:monolithicscheme}, the dual variables are condensed out, and the primal variables and dual variables can be computed using $\Kaa^{'}\Delta\boldsymbol{\alpha}=-\mathbf{r}_{\alpha}^{'}$ and Eqs.~\eqref{eq:condensedld} and~\eqref{eq:condensedlat}, respectively.

\subsubsection*{Damped Newton step}
A damped Newton step is performed with step lengths $t^\text{prim}_\alpha\in\Real_+$ and $t^\text{dual}_\alpha\in\Real_+$ such that the inequality constraints in the KKT conditions~\eqref{eq:stagdamagekkt1} and~\eqref{eq:stagdamagekkt2} hold for all iterations. As such, the values of the new iterate are given by
\begin{equation}
\begin{alignedat}{1}
\boldsymbol{\alpha}^{(k+1)}&\coloneqq\boldsymbol{\alpha}^{(k)}+t^\text{prim}_\alpha\Delta\boldsymbol{\alpha}, \\
\big\{\bld^{(k+1)},\blat^{(k+1)}\big\}&\coloneqq\big\{\bld^{(k)},\blat^{(k)}\big\}+t^\text{dual}_\alpha\big\{\Delta\bld,\Delta\blat\big\},
\end{alignedat}
\end{equation}
where the step lengths $t^\text{prim}_\alpha$ and $t^\text{dual}_\alpha$ are computed as
\begin{equation}
\begin{alignedat}{1}
t^\text{prim}_\alpha&\coloneqq\text{max}\left\{t\in(0,1]\,\vert\begin{bmatrix}\boldsymbol{\alpha}^{(k)}-\boldsymbol{\alpha}_n \\ \mathbf{e}-\boldsymbol{\alpha}^{(k)}\end{bmatrix}+t\begin{bmatrix}\Delta\boldsymbol{\alpha} \\ -\Delta\boldsymbol{\alpha}\end{bmatrix}\geq(1-\tau)\begin{bmatrix}\boldsymbol{\alpha}^{(k)}-\boldsymbol{\alpha}_n \\ \mathbf{e}-\boldsymbol{\alpha}^{(k)}\end{bmatrix}\right\}, \\
t^\text{dual}_\alpha&\coloneqq\text{max}\left\{t\in(0,1]\,\vert\begin{bmatrix}\bld^{(k)} \\ \blat^{(k)}\end{bmatrix}+t\begin{bmatrix}\Delta\bld \\ \Delta\blat\end{bmatrix}\geq(1-\tau)\begin{bmatrix}\bld^{(k)} \\ \blat^{(k)}\end{bmatrix}\right\}.
\end{alignedat}
\end{equation}

\clearpage

\section{Numerical examples} \label{sec:examples}
This section illustrates the performance of the monolithic and staggered interior-point methods using three numerical examples. The first example considers brittle fracture where a crack initiates at a notch and propagates brutally, corresponding to a mode I fracture. In the second example, the same specimen is loaded horizontally resulting in a crack that propagates less brutally through a region of high tensile stresses. The third example considers ductile fracture and shows the initiation of two cracks, followed by their propagation within a shear band of localized plastic strains. In addition, the results are benchmarked against the ones obtained from a staggered scheme, where irreversibility of the crack phase-field is treated using the history field of Miehe~et~al.~\cite{mieh10b} or the augmented Lagrangian approach of Wheeler et al.~\cite{whee14a}, and the plastic evolution problem is solved using the interior-point method presented in section~\ref{sec:stagplastic}. Table~\ref{tab:labels} summarizes the labels of the different solution schemes used in this section. All computations are performed on an Intel(R) Xeon(R) Gold 6140 2.30~GHz processor.
\vspace{-0.2em}
\begin{table}[h]
\centering
\caption{Solution schemes used in section~\ref{sec:examples}.}
\vspace{-0.4em}
\begin{tabular}{llll}
\toprule
Label & Scheme & Plastic problem & Damage problem \\
\midrule
Mono IP & monolithic & \mbox{\hspace{9mm}}- & interior-point \\
Mono IP-IP & monolithic & interior-point & interior-point \\
Stag IP & staggered & \mbox{\hspace{9mm}}- & interior-point \\
Stag HF & staggered & \mbox{\hspace{9mm}}- & history field \\
Stag AL & staggered & \mbox{\hspace{9mm}}- & augmented Lagrangian \\
Stag IP-HF & staggered & interior-point & history field \\
Stag IP-IP & staggered & interior-point & interior-point \\
\bottomrule
\end{tabular}
\label{tab:labels}
\end{table}

\vspace{-1.3em}
\subsection{Single-edge notched specimen subject to tension}
\begin{figure}[!h]
\centering
\includegraphics[width=0.38\textwidth]{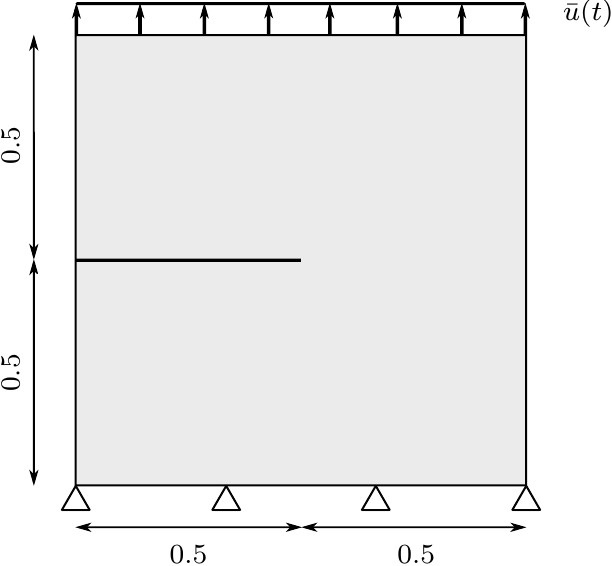}
\vspace{-1em}
\caption{Geometry and boundary conditions of the single-edge notched specimen subject to tension. The dimensions are given in mm.}
\label{fig:sen_tensile}
\end{figure}

The example of the single-edge notched specimen was originally proposed by Bourdin et~al.~\cite{bour00a}, and was later adopted as a benchmark problem by many authors~\cite{mieh10a,mieh10b,heis15a,amba15b,rodr18a,brun20a}. In this example, a square specimen with sides of 1~mm and a pre-existing notch is considered (Fig.~\ref{fig:sen_tensile}). Vertical displacements are imposed on the top edge in steps of $1\times10^{-4}$~mm, while the bottom edge remains completely fixed. The specimen behaves in a brittle way, governed by an AT-2 model, and has a bulk modulus ${K=175000}$~MPa, a shear modulus $G=80769.23$~MPa, a local damage dissipation constant $w_0=101.25$~MPa and a damage length scale ${\etad=0.142}$~N\textsuperscript{1/2}, corresponding to the material properties in Gerasimov and De Lorenzis~\cite{gera19a} (i.e.~a fracture toughness $G_\mathrm{c}=2.7$~N/mm and a damage characteristic length $\ell_\mathrm{d}=0.01$~mm). No plastic effects are considered. In addition, plane strain conditions are assumed. A finite element discretization with 17234 bilinear quadrilateral elements and 17421 nodes is considered, resulting in an element size $h_\mathrm{c}=0.0025$~mm in the critical zone where fracture is expected to occur.

\begin{figure}[!h]
\centering
\subfigure[]{\includegraphics[width=0.325\linewidth]{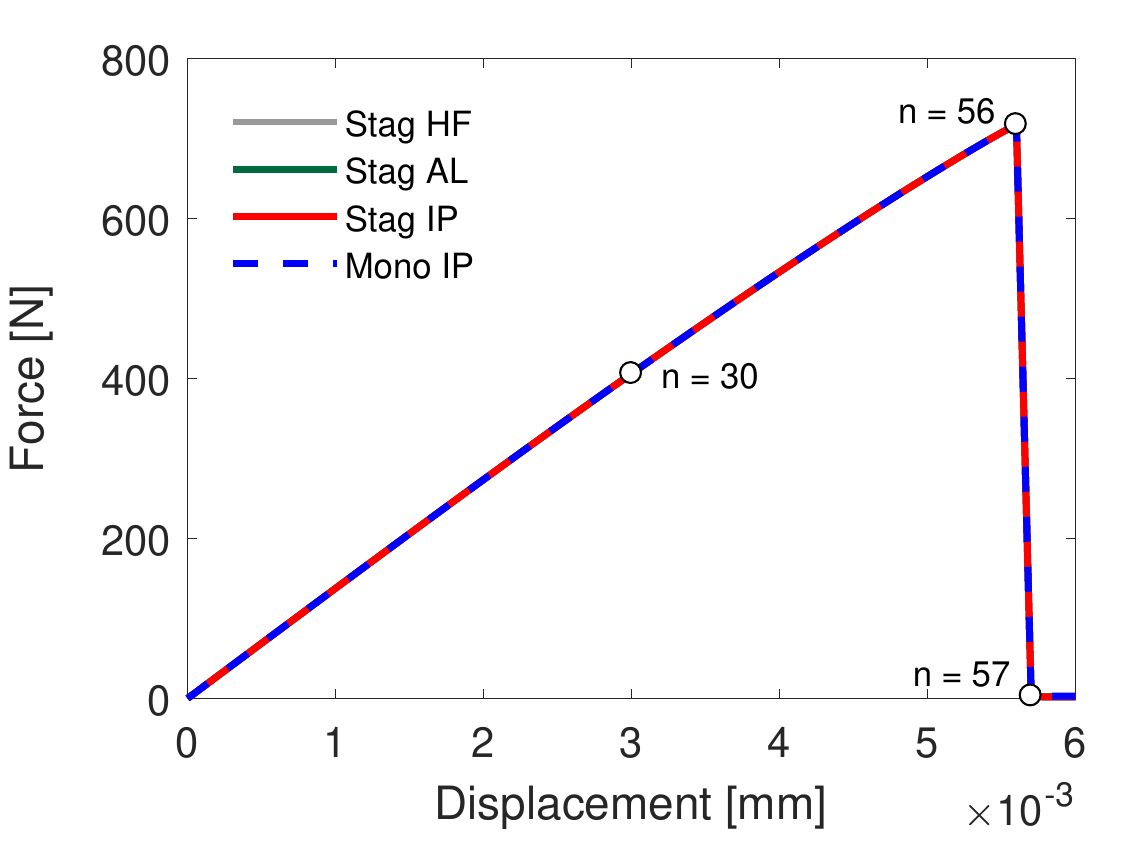}\label{fig:sen_tensile_force}}
\subfigure[]{\includegraphics[width=0.325\linewidth]{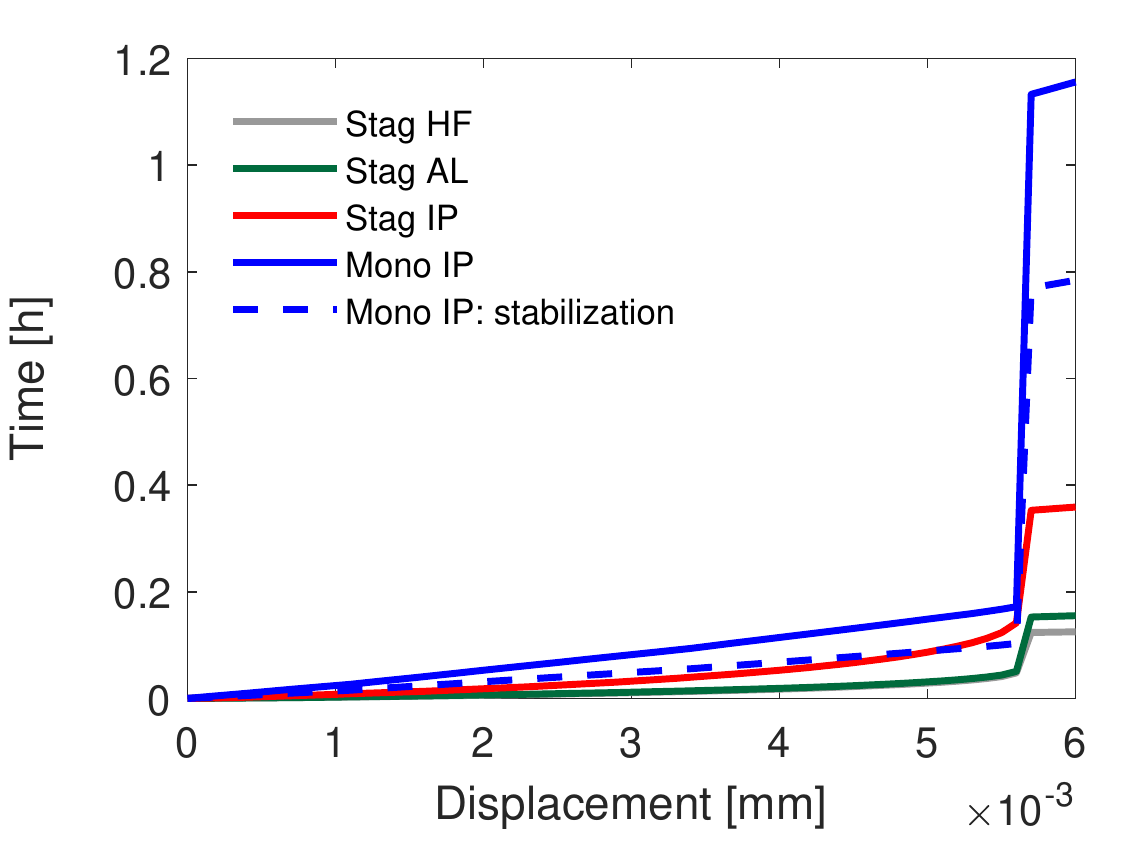}\label{fig:sen_tensile_time}}
\subfigure[]{\includegraphics[width=0.325\linewidth]{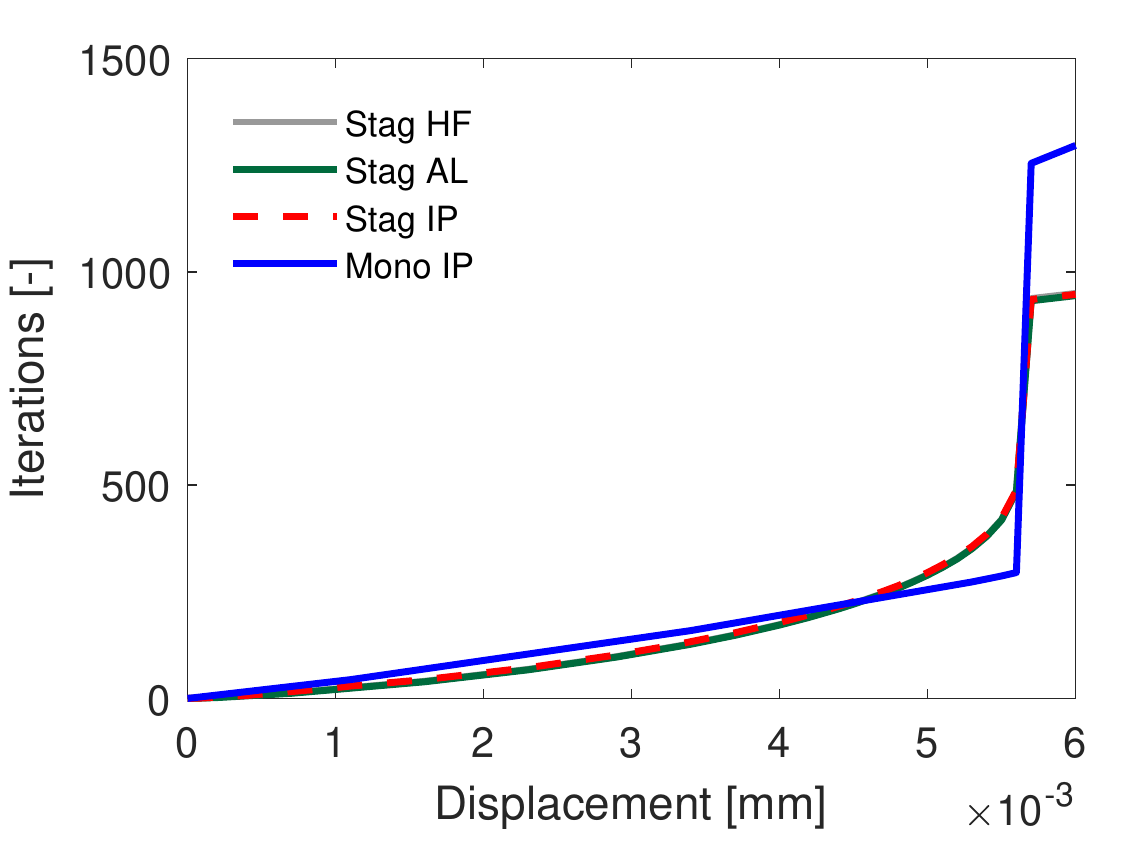}\label{fig:sen_tensile_iter}}
\caption{(a) Load-displacement curve, (b) cumulative computation time, and (c) cumulative number of iterations for the single-edge notched specimen subject to tension. The markers correspond to the time steps at which the energy evolution and residual are shown in Fig.~\ref{fig:sen_tensile_convergence}.}
\end{figure}

\vspace{-1em}
Fig.~\ref{fig:sen_tensile_force} shows the load-displacement curves computed using the staggered scheme and the monolithic scheme. In addition, the results are compared to the results obtained from a staggered scheme where irreversibility is imposed using the history field and an augmented Lagrangian. A penalty constant ${\gamma=2.6\times10^{6}}$~MPa, determined by Gerasimov and De Lorenzis~\cite{gera19a}, is used to obtain the results of the augmented Lagrangian approach. The solutions obtained using all schemes match for all time steps. The curves show that the specimen behaves quasi-linearly until time step $n=56$, after which the force abruptly decreases.
\vspace{-1em}

\vspace{-1em}
Fig.~\ref{fig:sen_tensile_time} shows the computational time required for all schemes. The computational time for the staggered and monolithic schemes where the interior-point method is used is larger than the computational time required for the staggered scheme using the history field or the augmented Lagrangian. In case the history field is used in combination with a staggered scheme, the damage evolution problem is cast in an unconstrained linear equation and requires no iterative solution procedure, and therefore has a smaller computational time. In addition, Fig.~\ref{fig:sen_tensile_time} shows that the majority of the computational time in the monolithic scheme is used for the symmetric indefinite factorization of the Hessian matrix, which is required to compute the stabilization factor.

Fig.~\ref{fig:sen_tensile_iter} shows the cumulative number of iterations for all schemes. The number of factorizations required to obtain the stabilization factor is not included in the number of iterations of the monolithic scheme. In the quasi-linear regime, all schemes require a low number of iterations to converge. However, during time step $n=57$ where the force abruptly decreases, the number of iterations increases significantly. In the monolithic scheme, the stabilization of the Hessian matrix is only active during this time step. In all other time steps, the Hessian matrix is positive definite for all iterations. Note that without the stabilization of the Hessian matrix, the presented monolithic scheme does not converge at this time step. The cumulative number of iterations required in the monolithic scheme is larger than the number of staggered iterations, which is almost independent of the different treatments of irreversibility.

Table~\ref{tab:sen_tensile_stats} summarizes the computational time and number of iterations for the computations of time step $n=57$. Table~\ref{tab:sen_tensile_stats}(a) shows that the number of staggered iterations and the number of equilibrium sub-iterations is very similar for the different treatments of irreversibility. Compared to the augmented Lagrangian approach, the number of sub-iterations required to solve the constrained damage evolution problem using the interior-point method is 40\% larger for this time step. Table~\ref{tab:sen_tensile_stats}(b) summarizes the results for the monolithic scheme, where it can be observed that the factorizations of the Hessian matrix required to compute the stabilization factor correspond to the majority of the computational time.
\begin{table}[h]
\centering
\caption{Computational time and number of iterations for time step $n=57$ at which the crack propagates brutally in the single-edge notched specimen subject to tension.}
\vspace{-1em}
\begin{tabular}{llll}
\multicolumn{4}{c}{\small{(a) Staggered scheme}} \\
\toprule
& Interior-point & History field & Augmented Lagrangian \\
\midrule
Computational time [s] & 763 & 273 & 367 \\
Staggered iterations & 447 & 450 & 447 \\
Equilibrium sub-iterations & 1308 & 1305 & 1306 \\
Damage sub-iterations & 2929 & 450 & 2091 \\
\bottomrule
\end{tabular}
\begin{tabular}{ll}
\vspace{-1em}
\\
\multicolumn{2}{c}{\small{(b) Monolithic scheme}} \\
\toprule
& Interior-point \\
\midrule
Computational time [s] & 3459 \\
Monolithic iterations & 959 \\
Hessian factorizations & 1910 \\
Factorization time [s] & 2399 \\
\bottomrule
\end{tabular}
\label{tab:sen_tensile_stats}
\end{table}

\begin{figure}[!h]
\centering
\subfigure[]{\includegraphics[width=0.325\linewidth]{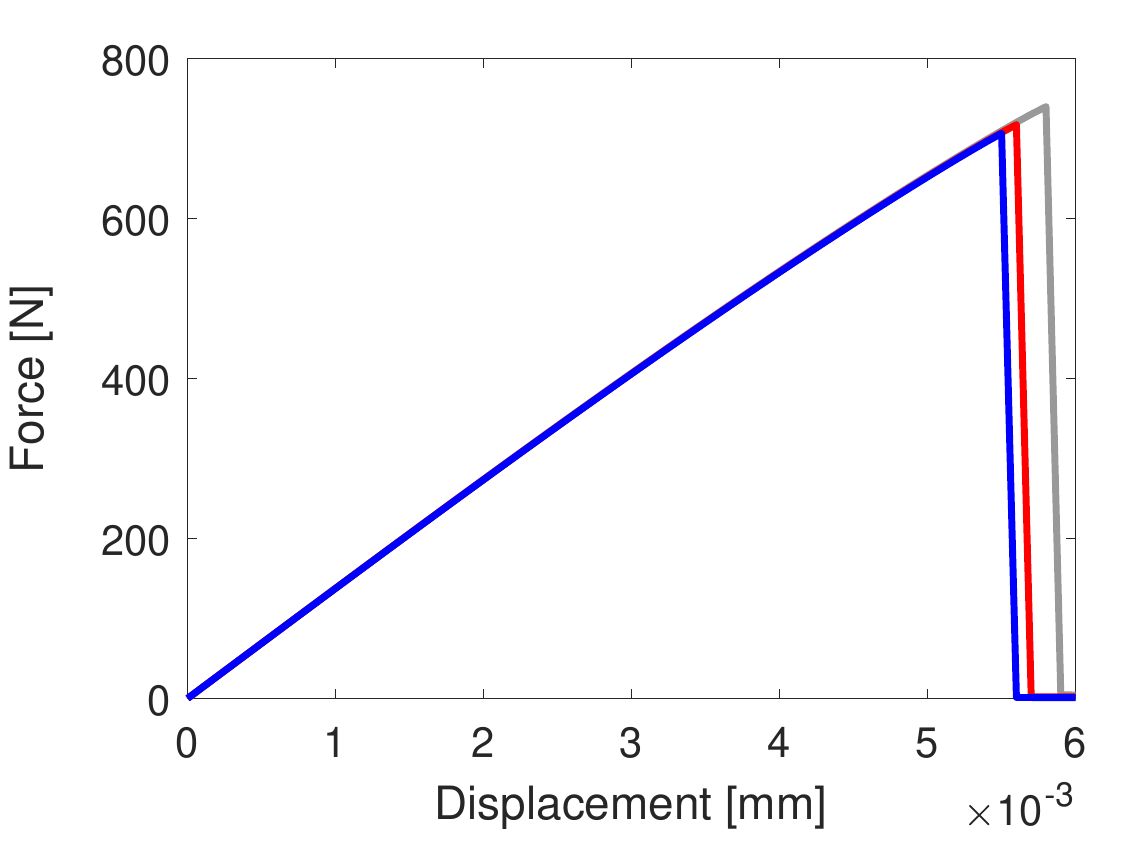}\label{fig:sen_tensile_mesh_force}}
\subfigure[]{\includegraphics[width=0.325\linewidth]{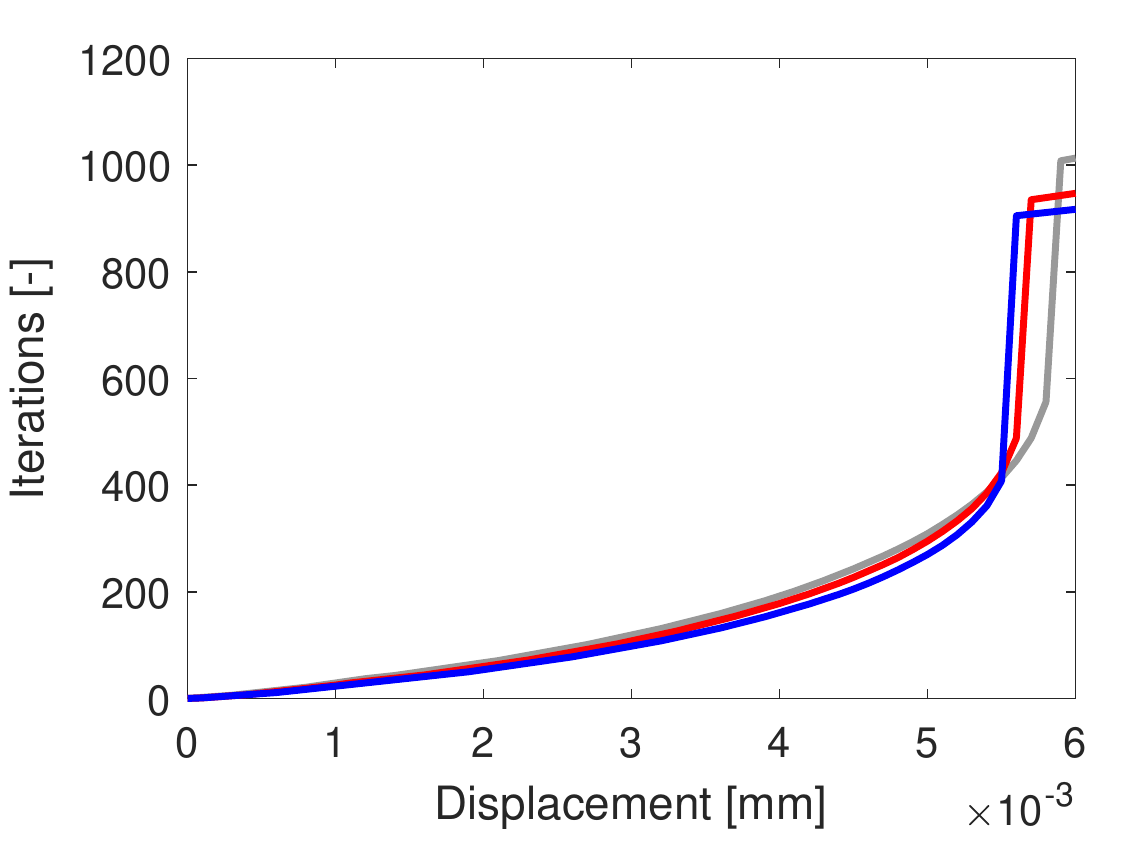}\label{fig:sen_tensile_mesh_stag}}
\subfigure[]{\includegraphics[width=0.325\linewidth]{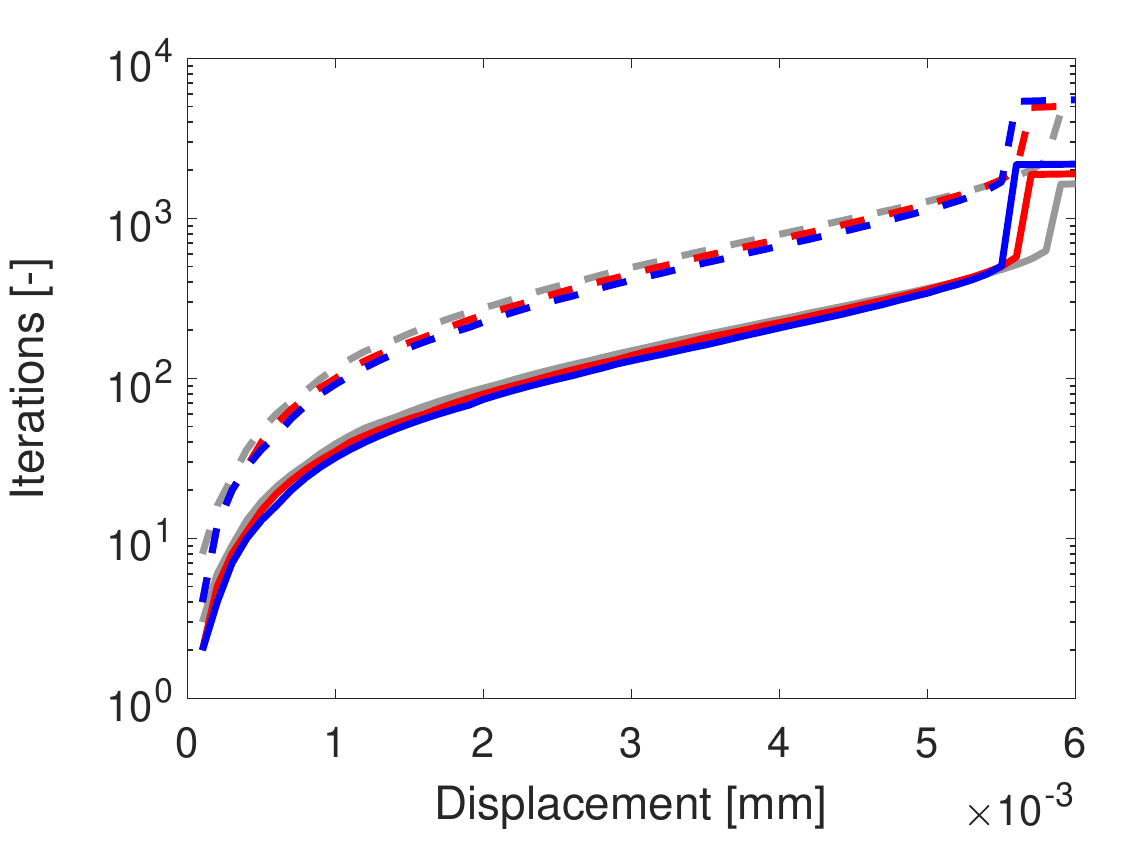}\label{fig:sen_tensile_mesh_sub}}
\caption{(a) Load-displacement curve, (b) cumulative number of staggered iterations, and (c) cumulative number of sub-iterations for the single-edge notched specimen subject to tension computed using a mesh with $h_\mathrm{c}=0.5\ell_\mathrm{d}$ (gray), $h_\mathrm{c}=0.25\ell_\mathrm{d}$ (red), and $h_\mathrm{c}=0.125\ell_\mathrm{d}$ (blue). The solid and dashed lines correspond to equilibrium and damage sub-iterations, respectively.}
\end{figure}

Fig.~\ref{fig:sen_tensile_mesh_force} shows the load-displacement curve obtained using the interior-point method in a staggered scheme for different mesh sizes. Mesh-convergence can be observed. In addition, Figs.~\ref{fig:sen_tensile_mesh_stag} and~\ref{fig:sen_tensile_mesh_sub} show the cumulative number of staggered iterations and sub-iterations, respectively. Note that in Fig.~\ref{fig:sen_tensile_mesh_sub}, the number of sub-iterations is cumulated over both pseudo-time and the staggered iteration loop, as will be done hereafter for all sub-iterations. In this example, the total number of staggered iterations is not very sensitive to the mesh size. In addition, the average number of damage sub-iterations per staggered iteration (i.e.~4.19, 4.41 and 4.67 for the coarse, medium and fine mesh, respectively) does not change significantly.

\begin{figure}[!h]
\centering
\hspace{12mm}
\addtolength{\subfigcapmargin}{-6mm}
\subfigure[]{\includegraphics[height=4.5cm]{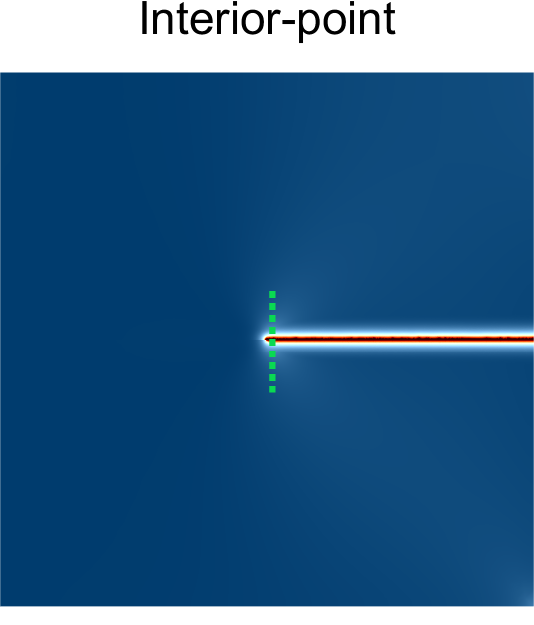}\label{fig:sen_tensile_dam1}}
\hspace{7mm}
\subfigure[]{\includegraphics[height=4.5cm]{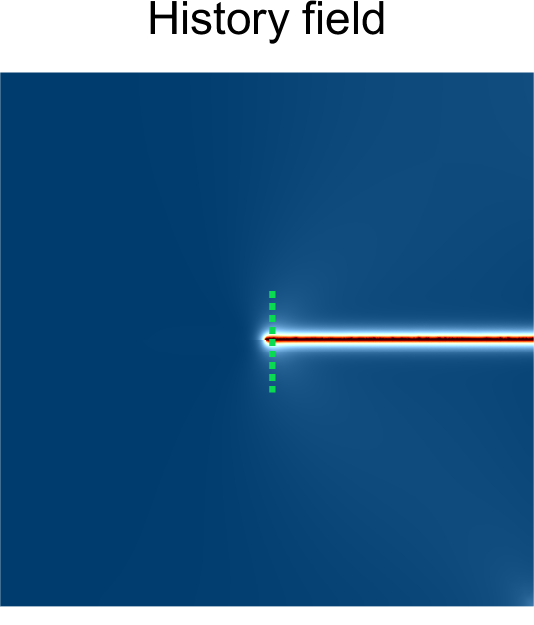}}
\hspace{7mm}
\subfigure[]{\includegraphics[height=4.5cm]{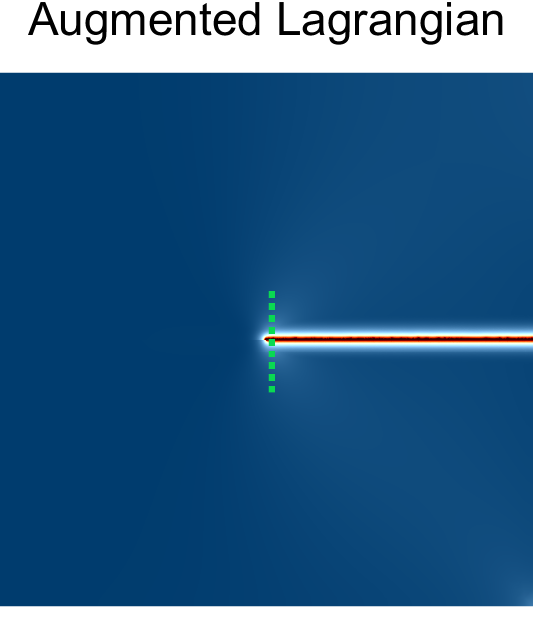}\hspace{0.6cm}\includegraphics[height=4.5cm]{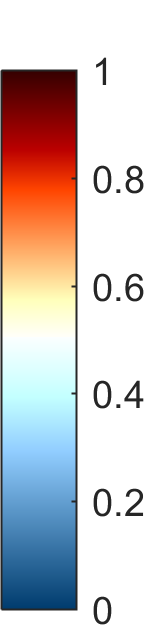}\label{fig:sen_tensile_dam3}}
\addtolength{\subfigcapmargin}{6mm}
\addtolength{\subfigcapskip}{-4mm}
\subfigure[]{\hspace{-11mm}\includegraphics[width=0.949\linewidth]{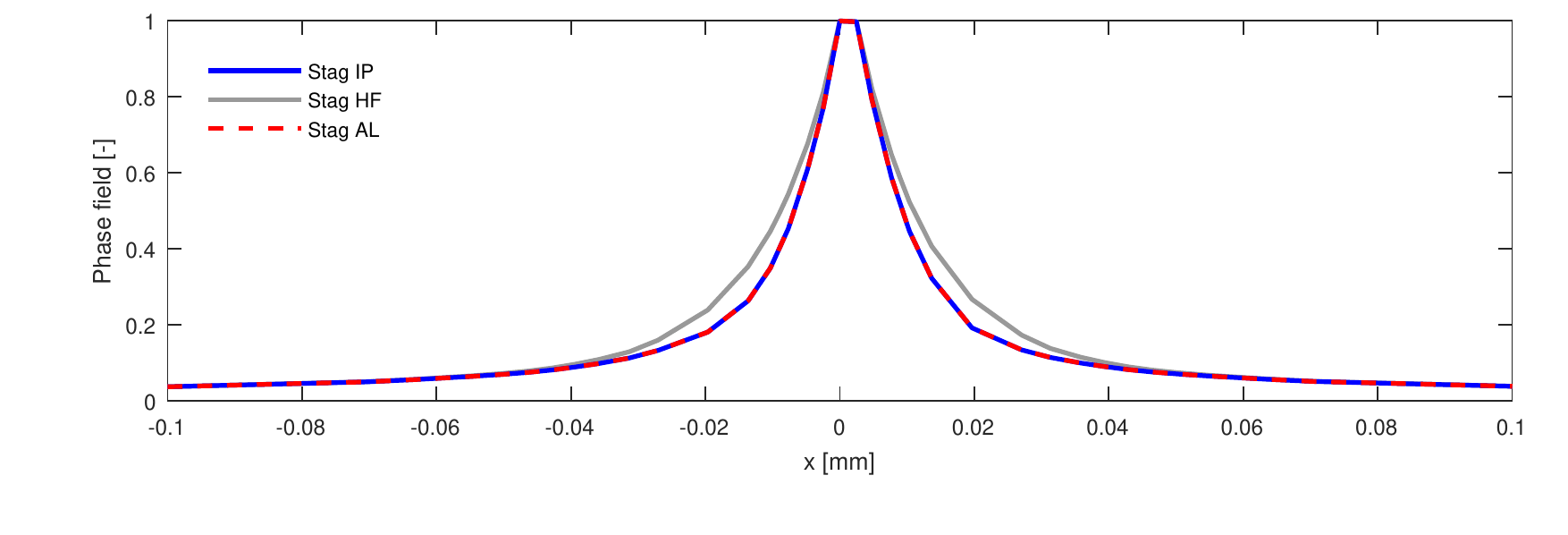}\label{fig:sen_tensile_profile}}
\addtolength{\subfigcapskip}{4mm}
\vspace{-2.1em}
\caption{Crack phase-field profile of the single-edge notched specimen subject to tension at the final time step obtained using (a) the interior-point method, (b) the history field, and (c) the augmented Lagrangian approach. The crack phase-field profile orthogonal to the crack, in a cross-section denoted with a dashed line in (a) to (c), is shown in (d).}
\vspace{-1.6em}
\end{figure}

Figs.~\ref{fig:sen_tensile_dam1} to \ref{fig:sen_tensile_dam3} show the crack phase-field profile at the final time step computed using three staggered schemes with different treatments of irreversibility. Fig.~\ref{fig:sen_tensile_profile} shows that the crack phase-field profile obtained using the history field is similar to that of the other approaches, but shows a slightly wider crack. This is attributed to the overestimation of the fracture surface energy that results from the use of the history field~\cite{gera19a}.
\begin{figure}[!h]
\centering
\subfigure{\includegraphics[width=0.325\linewidth]{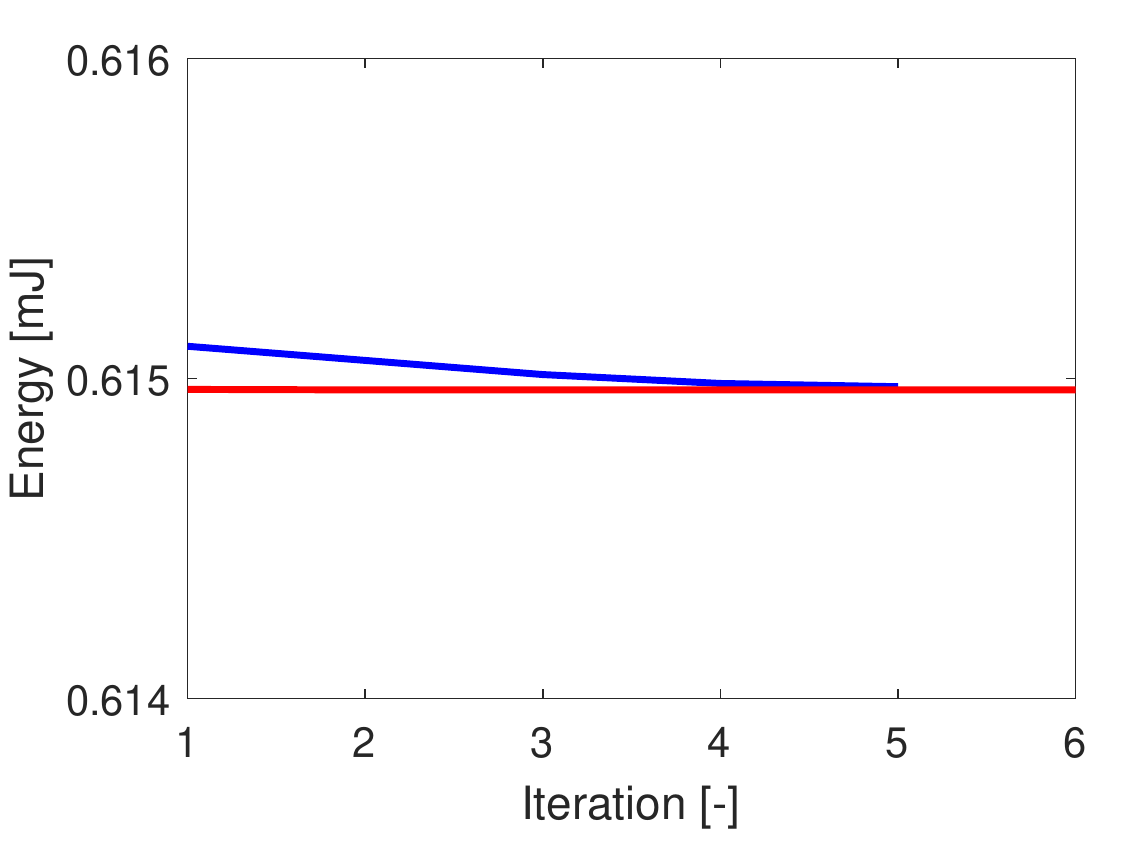}\label{fig:sen_tensile_nrg_30}}
\subfigure{\includegraphics[width=0.325\linewidth]{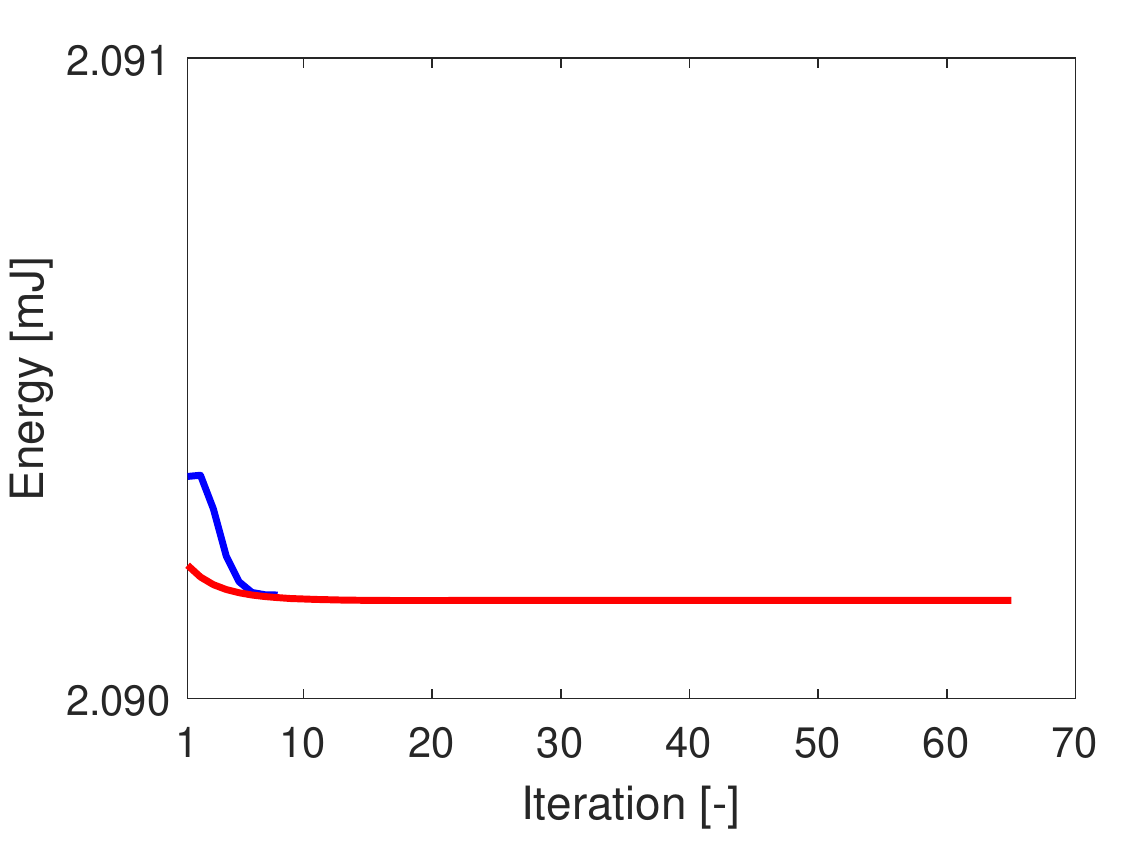}\label{fig:sen_tensile_nrg_56}}
\subfigure{\includegraphics[width=0.325\linewidth]{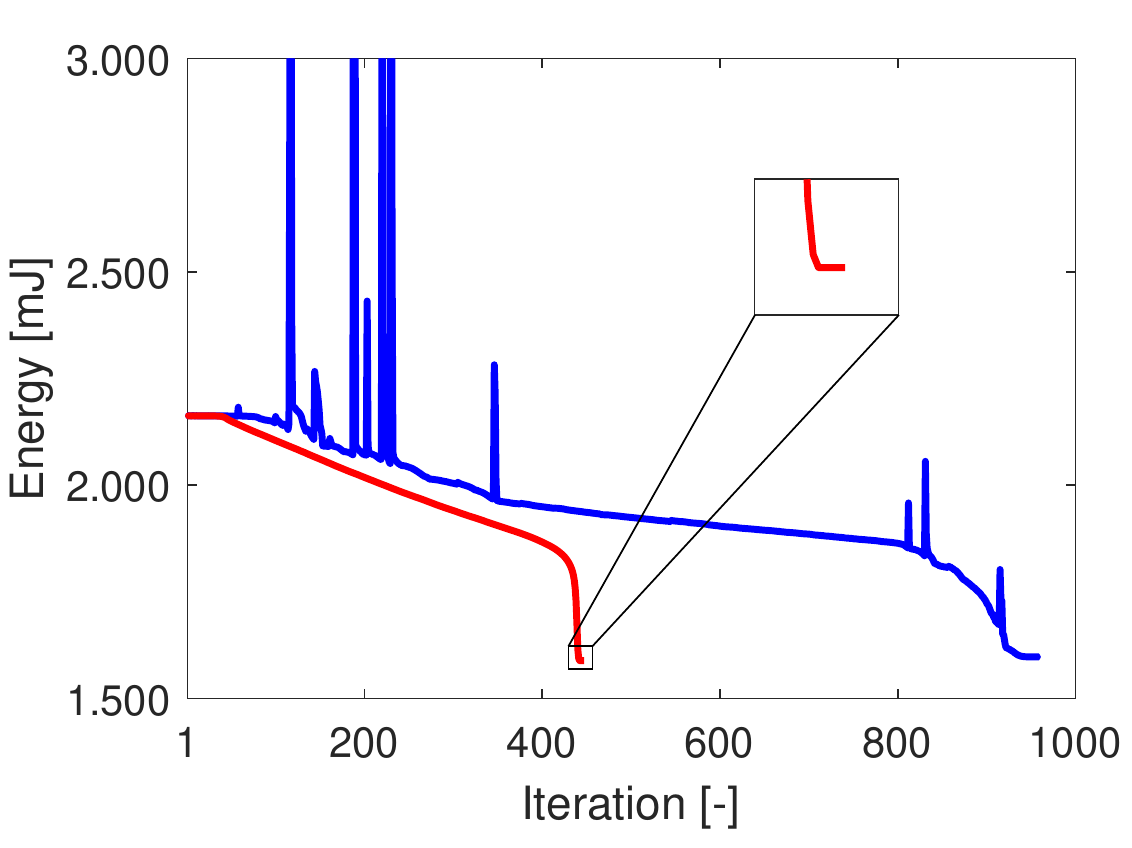}\label{fig:sen_tensile_nrg_57}}
\addtocounter{subfigure}{-3}
\subfigure[]{\includegraphics[width=0.325\linewidth]{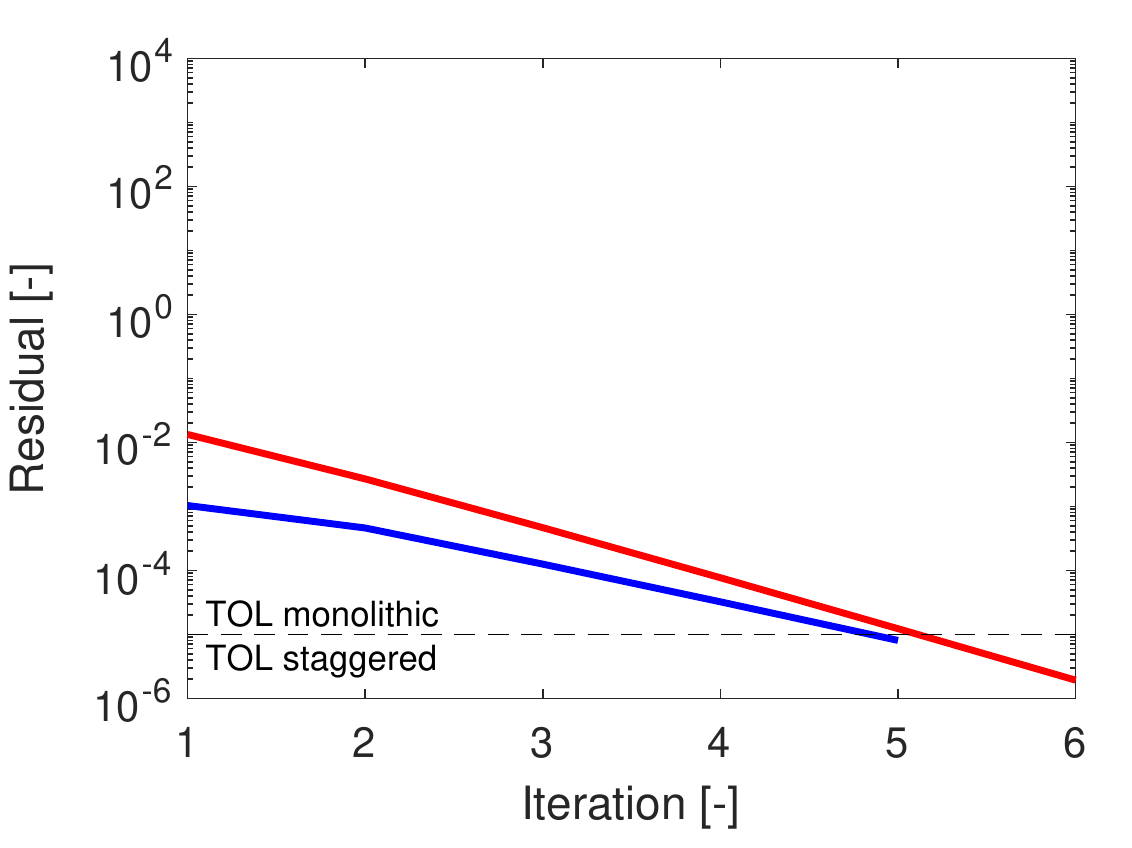}}
\subfigure[]{\includegraphics[width=0.325\linewidth]{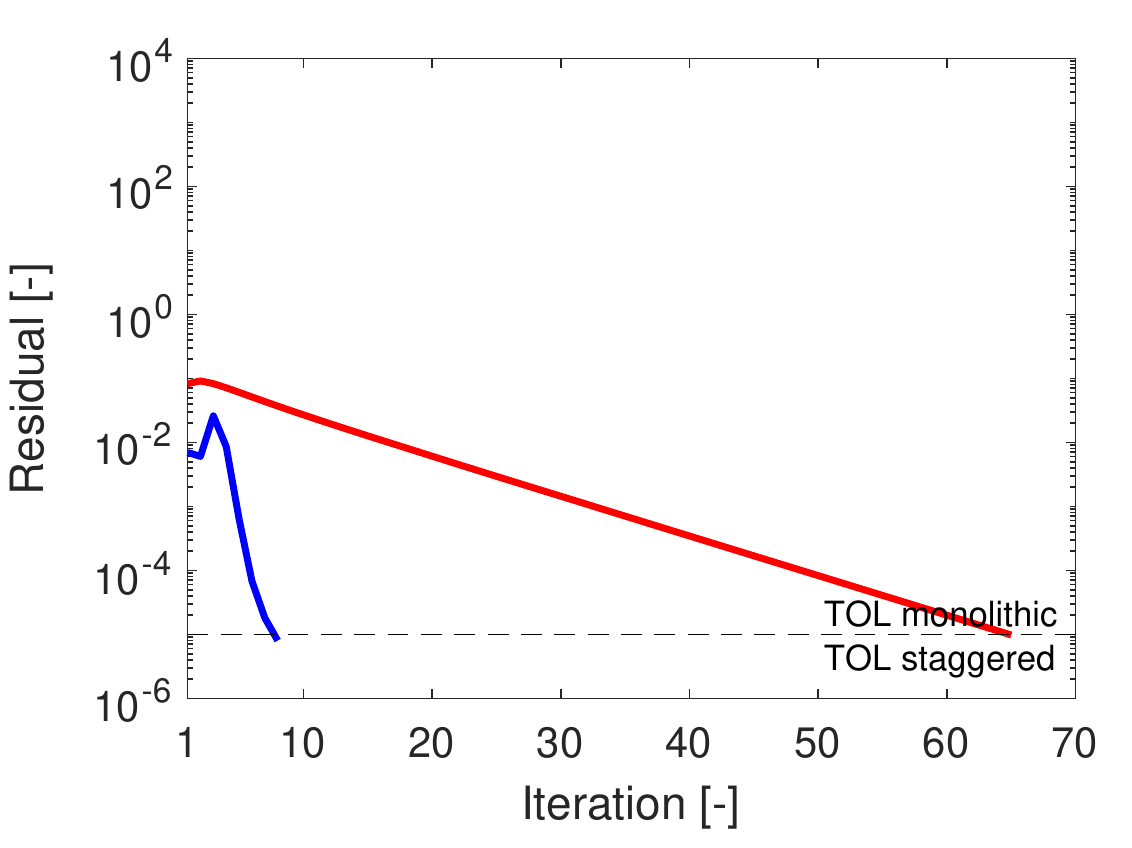}}
\subfigure[]{\includegraphics[width=0.325\linewidth]{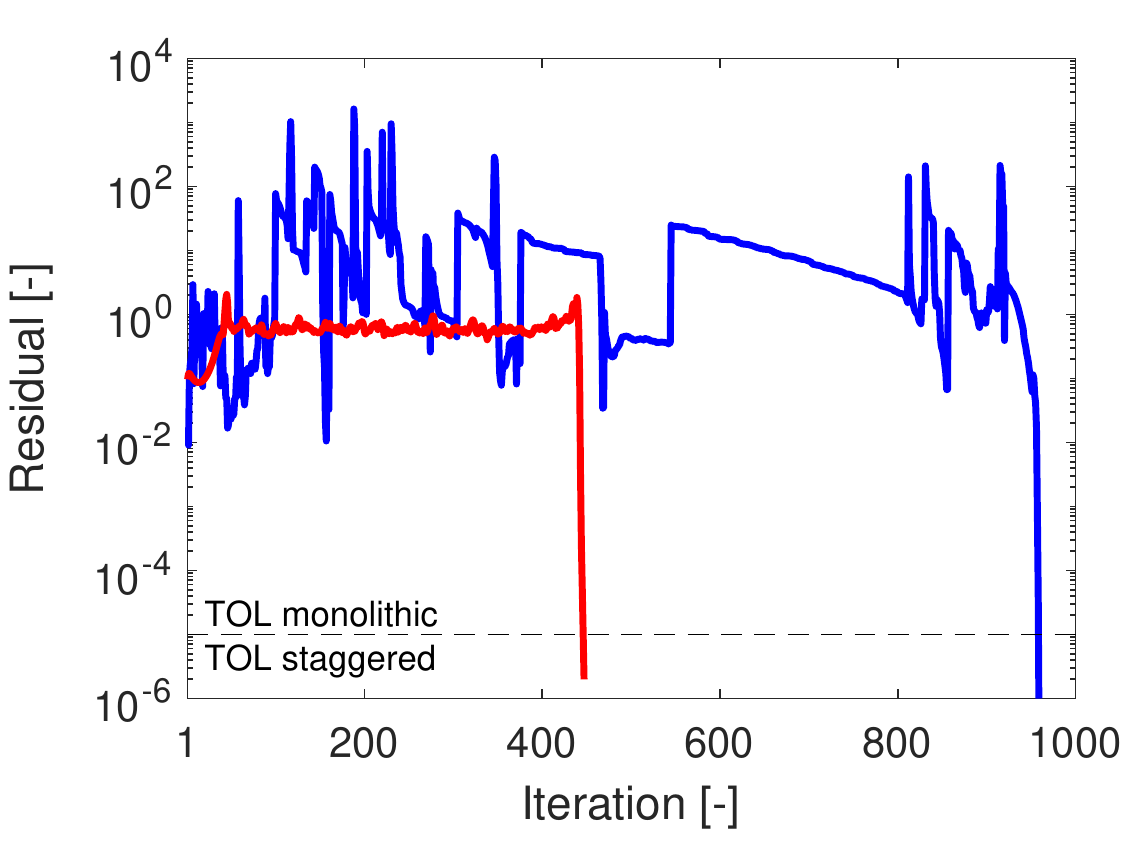}}
\vspace{-0.8em}
\caption{Evolution of the energy functional (top) and the residual (bottom) for the solution of time step (a) $n=30$, (b) $n=56$, and (c) $n=57$ of the single-edge notched specimen subject to tension. The energy and residual of the staggered scheme and monolithic scheme are shown in red and blue, respectively.}
\label{fig:sen_tensile_convergence}
\end{figure}

Fig.~\ref{fig:sen_tensile_convergence} shows the evolution of the energy functional and the residual during the minimization process. Figs.~\ref{fig:sen_tensile_nrg_30}~and~\ref{fig:sen_tensile_nrg_56} show that in the quasi-linear regime, the number of iterations required to converge to a stationary state is smaller in the monolithic scheme. However, in time step $n=57$, the number of iterations is larger in the monolithic scheme, resulting in the increase of computational time in Fig.~\ref{fig:sen_tensile_time}. The staggered scheme solves a sequence of convex sub-problems. As such, the energy functional decreases monotonically for all iterations. From Figs.~\ref{fig:sen_tensile_nrg_30} to \ref{fig:sen_tensile_nrg_57} it can be observed that, for this example, a tolerance ${\tt TOL}=10^{-5}$ is sufficient to obtain a solution where the energy functional is stationary.

In this example, all schemes yield similar results. The main difference is observed in the width of the crack phase-field, which is larger in case the history field is used. The computational time is the lowest when the staggered scheme is used in combination with the history field, which is a consequence of the modification to the original constrained damage evolution problem. The computational time is the largest for the monolithic scheme in combination with the interior-point method, attributed to the large number of factorizations of the Hessian matrix performed during brutal crack propagation. This could possibly be improved by changing the starting values and multiplication factors in algorithm~\ref{alg:inertiacorrection}. However, hereafter, the values given in algorithm~\ref{alg:inertiacorrection} will be used.

\subsection{Single-edge notched specimen subject to shear}
\begin{figure}[!h]
\centering
\includegraphics[width=0.38\textwidth]{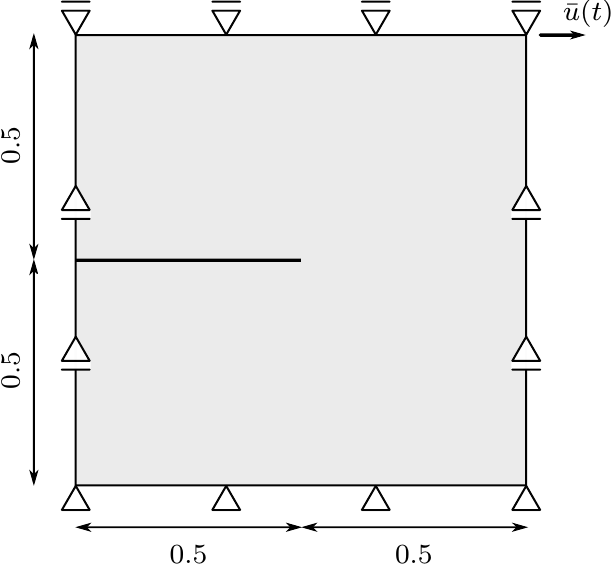}
\caption{Geometry and boundary conditions of the single-edge notched specimen subject to shear. The dimensions are given in~mm.}
\label{fig:sen_shear}
\end{figure}

In this example, the specimen considered in the previous section is loaded horizontally (Fig.~\ref{fig:sen_shear}). The top edge and side edges are restrained from moving in the vertical direction, while the bottom edge remains completely fixed. Horizontal displacements are imposed on the top edge in steps of $3\times10^{-4}$~mm. The loading conditions are chosen such that, contrary to the previous example, a crack will propagate downwards and in a less brutal way. A finite element discretization with 21788 bilinear quadrilateral elements and 21965 nodes is considered, resulting in an element size $h_\mathrm{c}=0.0027$~mm in the critical zone where fracture is expected to occur. In order to prevent volumetric locking in compressed regions of the specimen, a selective reduced integration technique is used to integrate the volumetric term of the elastic energy density in~\eqref{eq:directionalderivativesdamage} and the corresponding volumetric part of the stress tensor in~\eqref{eq:directionalderivativesincremental1}. A detailed overview and motivation of this procedure can be found in Alessi et al.~\cite{ales20a}. Because the history field approach is incompatible with selective reduced integration, no results are presented thereof.

\begin{figure}[!h]
\centering
\subfigure[]{\includegraphics[width=0.325\linewidth]{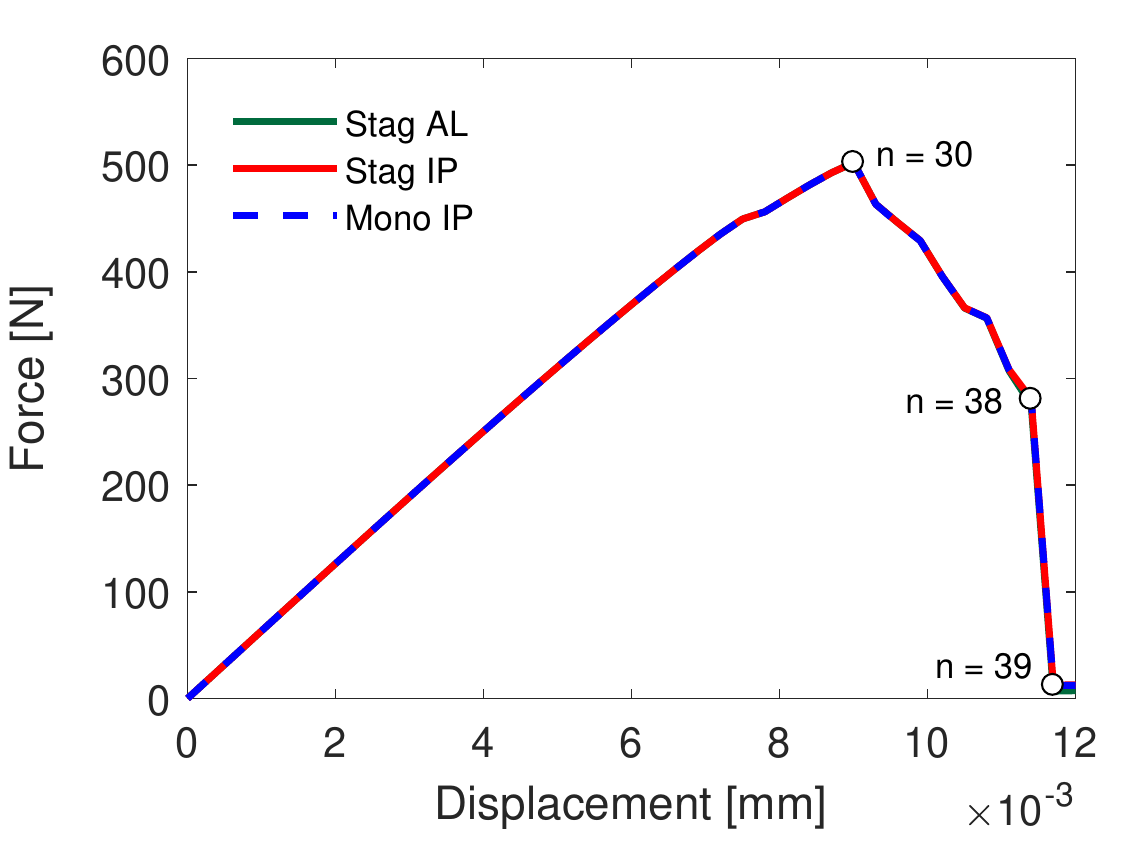}\label{fig:sen_shear_force}}
\subfigure[]{\includegraphics[width=0.325\linewidth]{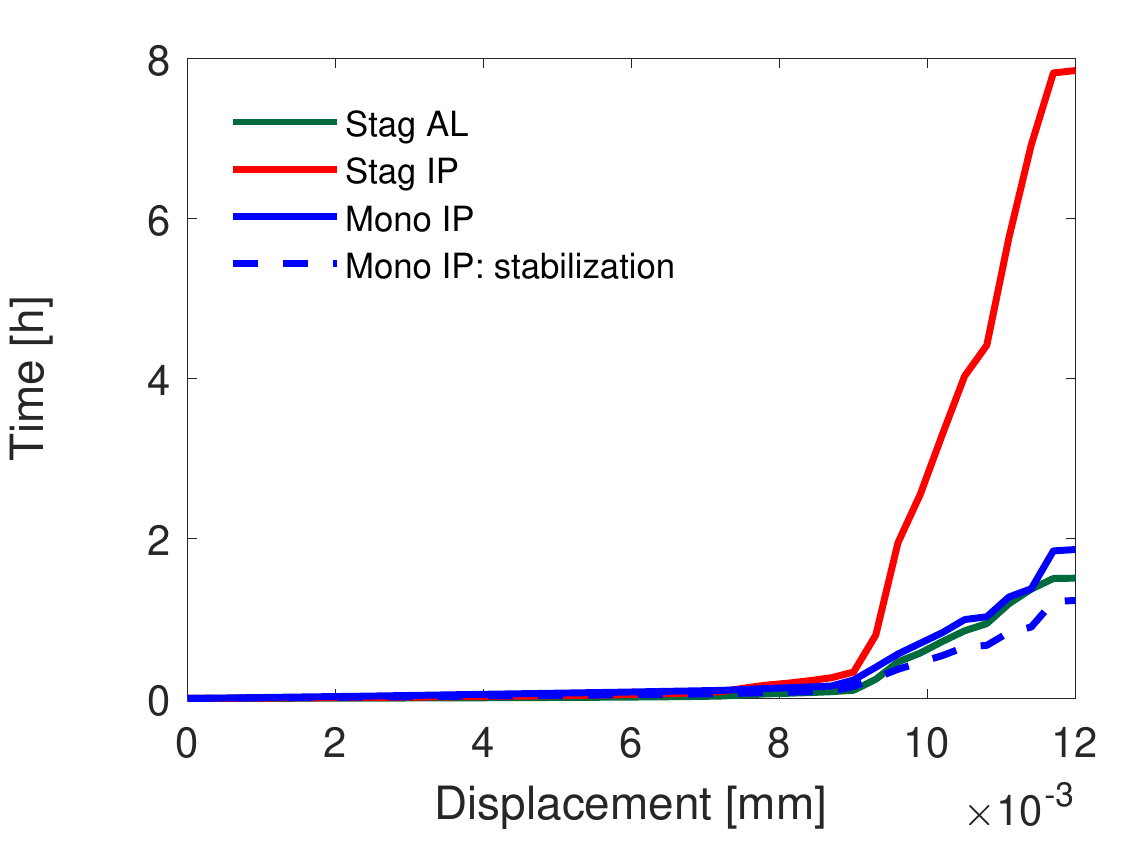}\label{fig:sen_shear_time}}
\subfigure[]{\includegraphics[width=0.325\linewidth]{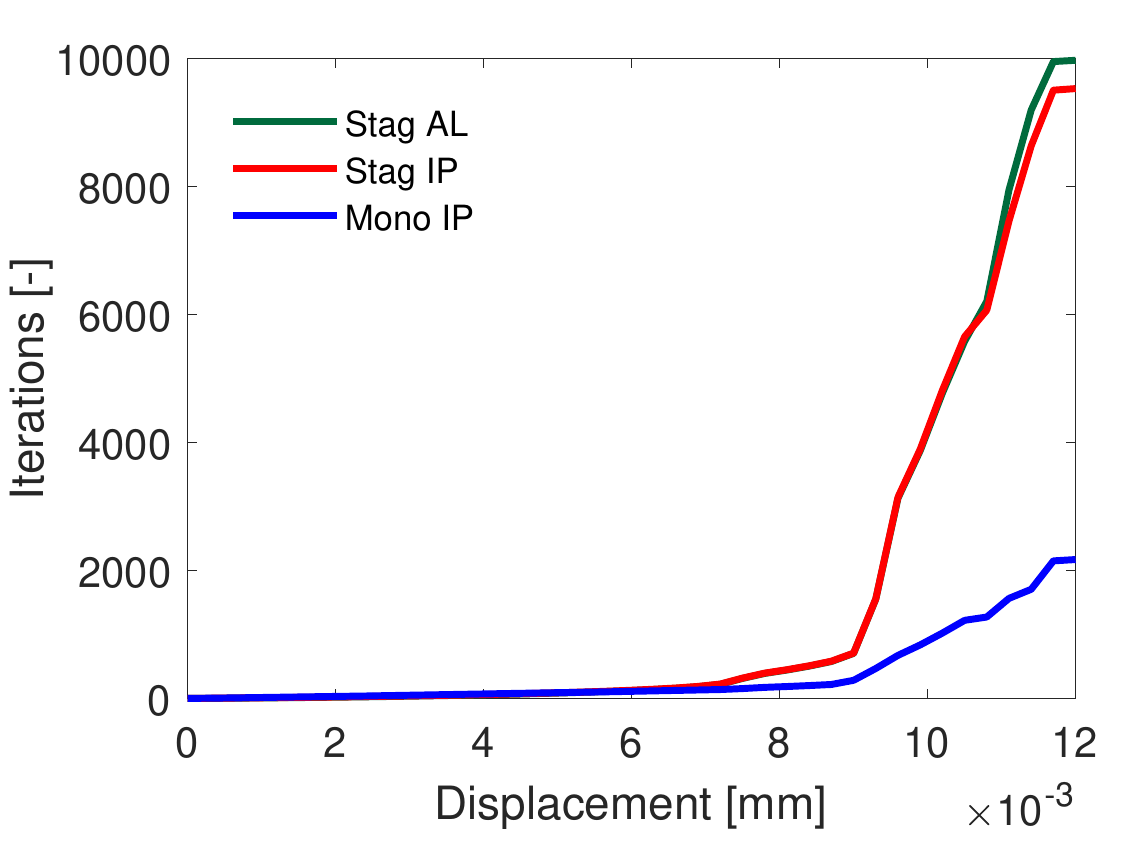}\label{fig:sen_shear_iter}}
\caption{(a) Load-displacement curve, (b) cumulative computation time, and (c) cumulative number of iterations for the single-edge notched specimen subject to shear. The markers correspond to the time steps at which the energy evolution and residual are shown in Fig.~\ref{fig:sen_shear_convergence}.}
\end{figure}

Fig.~\ref{fig:sen_shear_force} shows the load-displacement curve computed using the staggered scheme and the monolithic scheme. First, the specimen behaves quasi-linearly, followed by a softening regime during which a crack propagates. For a penalty constant ${\gamma=2.6\times10^{6}}$~MPa, the results of the augmented Lagrangian approach and the interior-point methods match well. In the monolithic scheme, the stabilization of the Hessian matrix is active between imposed displacements of $9\times10^{-3}$~mm and $12\times10^{-3}$~mm, during which softening is significant.

Figs.~\ref{fig:sen_shear_time} and~\ref{fig:sen_shear_iter} show the computational time and the number of iterations required during the computations, respectively. Similar to the previous example, all schemes require a small number of iterations to obtain a solution in the quasi-linear regime. However, during crack propagation, the number of iterations and the computational time increase. Contrary to the previous example, the computational time of the monolithic scheme is smaller than the computational time of the staggered scheme where an interior-point method is used to impose irreversibility. This is attributed to the less brutal crack propagation. In addition, when the interior-point method is adopted in a staggered scheme, the number of staggered iterations is similar to that of the augmented Lagrangian approach, but the computational time is significantly larger.

\begin{figure}[!h]
\centering
\subfigure[]{\includegraphics[width=0.325\linewidth]{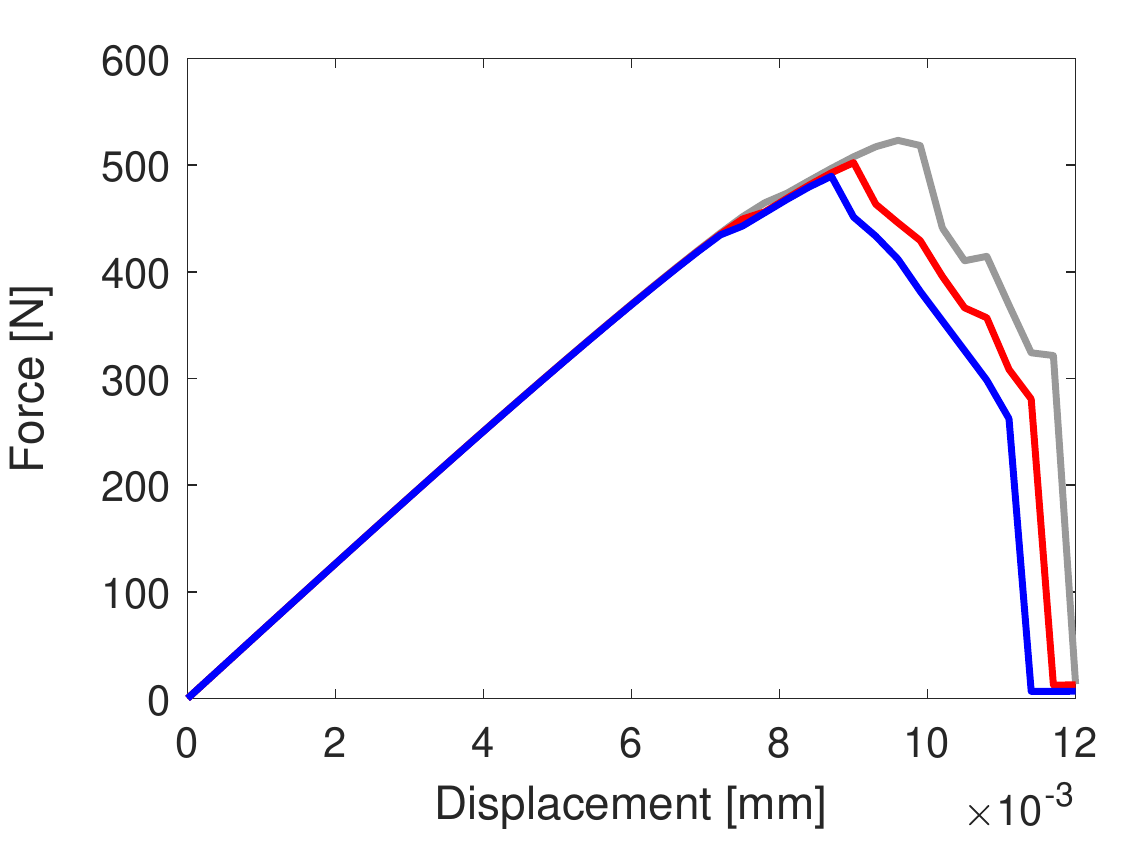}\label{fig:sen_shear_mesh_force}}
\subfigure[]{\includegraphics[width=0.325\linewidth]{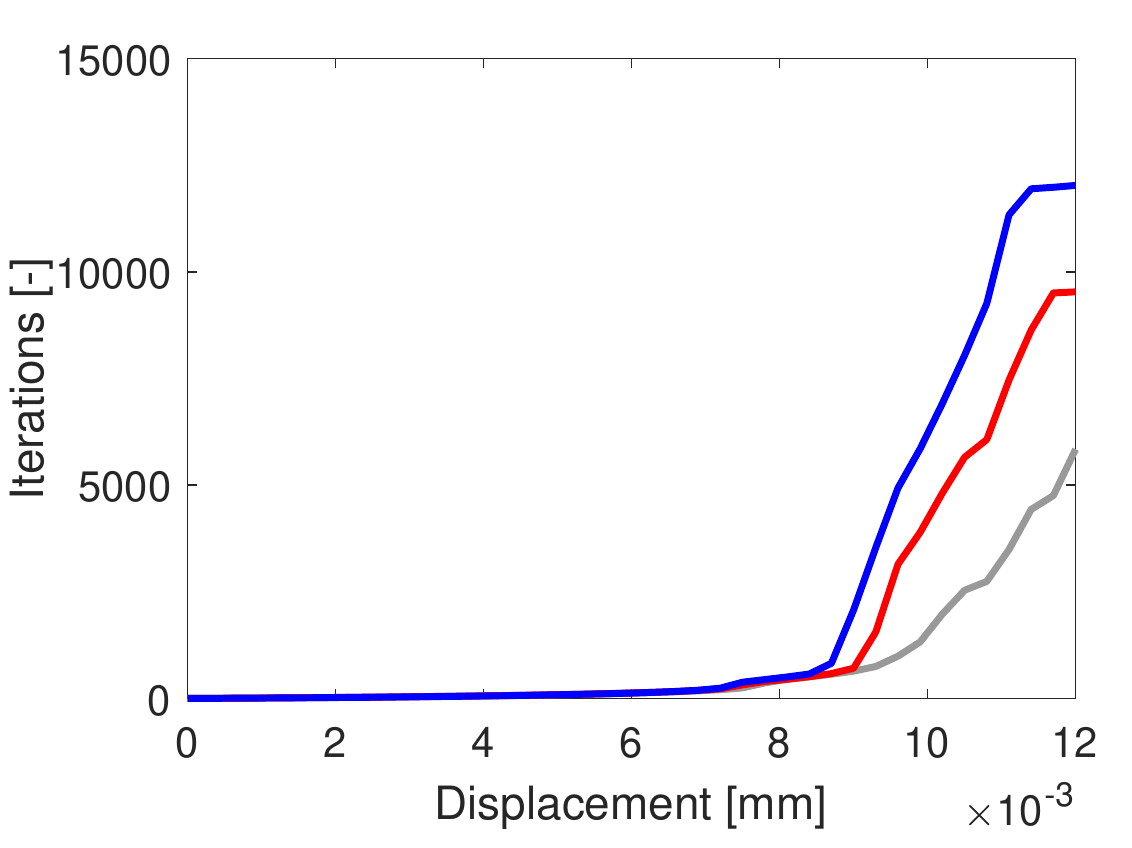}\label{fig:sen_shear_mesh_stag}}
\subfigure[]{\includegraphics[width=0.325\linewidth]{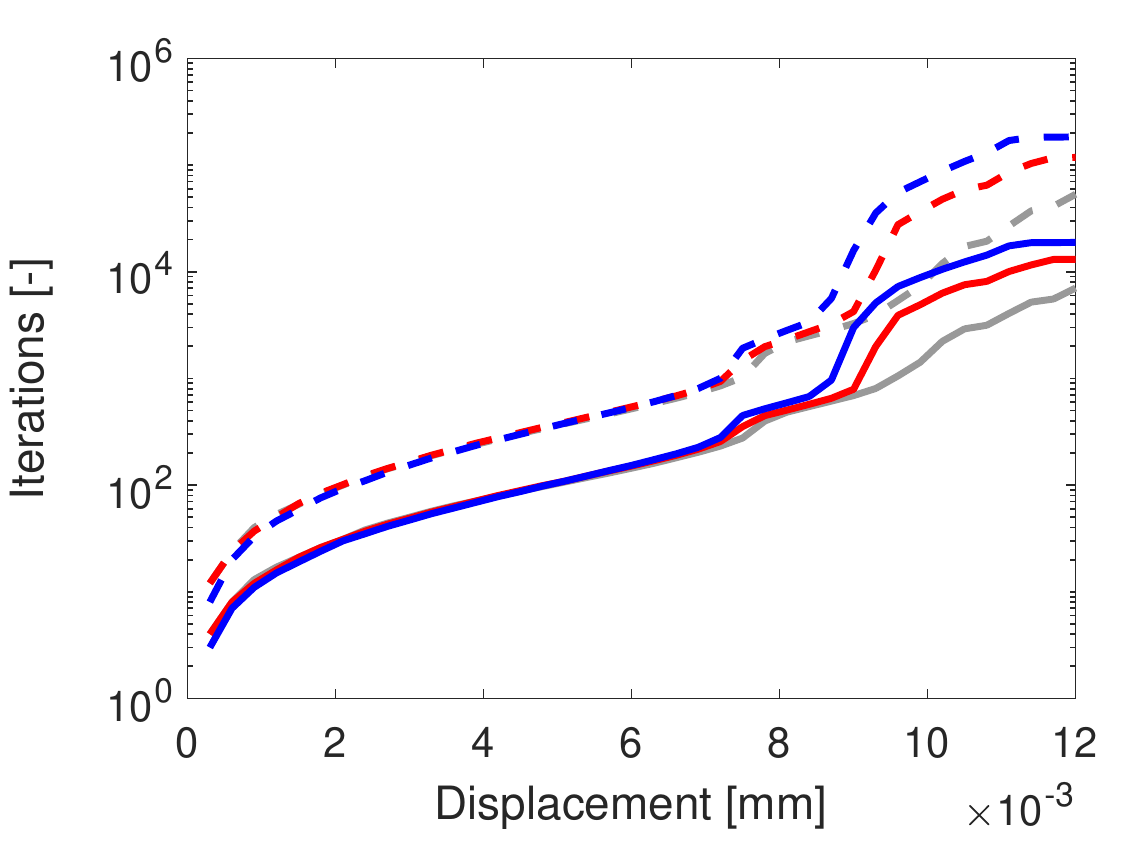}\label{fig:sen_shear_mesh_sub}}
\caption{(a) Load-displacement curve, (b) cumulative number of staggered iterations, and (c) cumulative number of sub-iterations for the single-edge notched specimen subject to shear computed using a mesh with $h_\mathrm{c}=0.546\ell_\mathrm{d}$ (gray), $h_\mathrm{c}=0.273\ell_\mathrm{d}$ (red), and $h_\mathrm{c}=0.136\ell_\mathrm{d}$ (blue). The solid and dashed lines correspond to equilibrium and damage sub-iterations, respectively.}
\end{figure}

Fig.~\ref{fig:sen_shear_mesh_force} shows the load-displacement curve obtained using the interior-point method in a staggered scheme for different mesh sizes. In addition, Figs.~\ref{fig:sen_shear_mesh_stag} and~\ref{fig:sen_shear_mesh_sub} show the cumulative number of staggered iterations and sub-iterations, respectively. In contrast to the previous example, it can be observed that the total number of staggered iterations considerably depends on the mesh size. The average number of damage sub-iterations per staggered iteration (5.54, 6.91 and 8.18 for the coarse, medium and fine mesh, respectively), however, is still of the same order.

\begin{figure}[!h]
\centering
\hspace{15mm}
\addtolength{\subfigcapmargin}{-6mm}
\subfigure[]{\includegraphics[height=4.5cm]{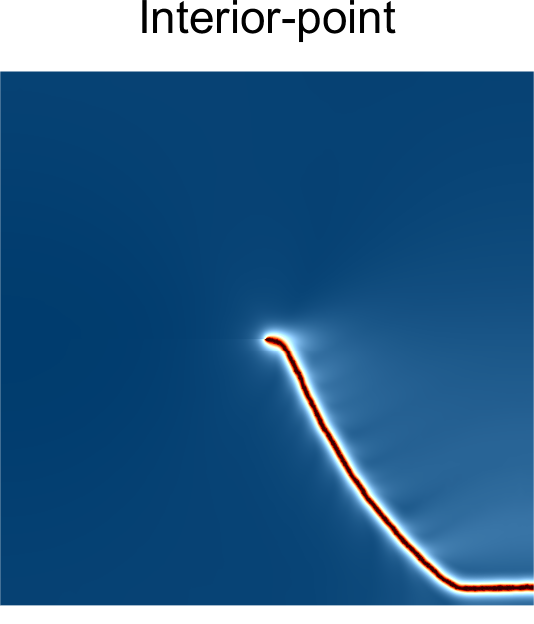}}
\hspace{7mm}
\subfigure[]{\includegraphics[height=4.5cm]{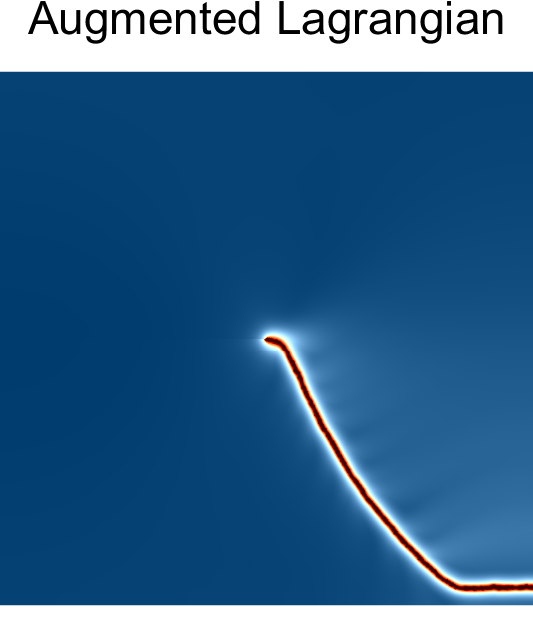}\hspace{0.6cm}\includegraphics[height=4.5cm]{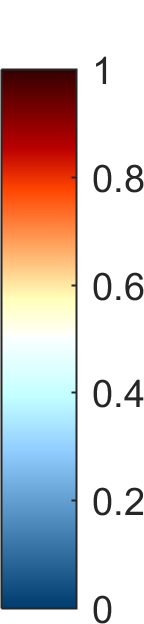}}
\addtolength{\subfigcapmargin}{6mm}
\vspace{-0.8em}
\caption{Crack phase-field profile of the single-edge notched specimen subject to shear at the final time step obtained using (a) the interior-point method and (b) the augmented Lagrangian approach.}
\label{fig:sen_shear_dam}
\vspace{-1.2em}
\end{figure}

Fig.~\ref{fig:sen_shear_dam} shows the crack phase-field profile computed using two staggered schemes at the final time step. There are no significant differences between the crack phase-field profile obtained using an interior-point method and the augmented Lagrangian approach.
\vspace{-1em}

\begin{figure}[!h]
\centering
\subfigure{\includegraphics[width=0.325\linewidth]{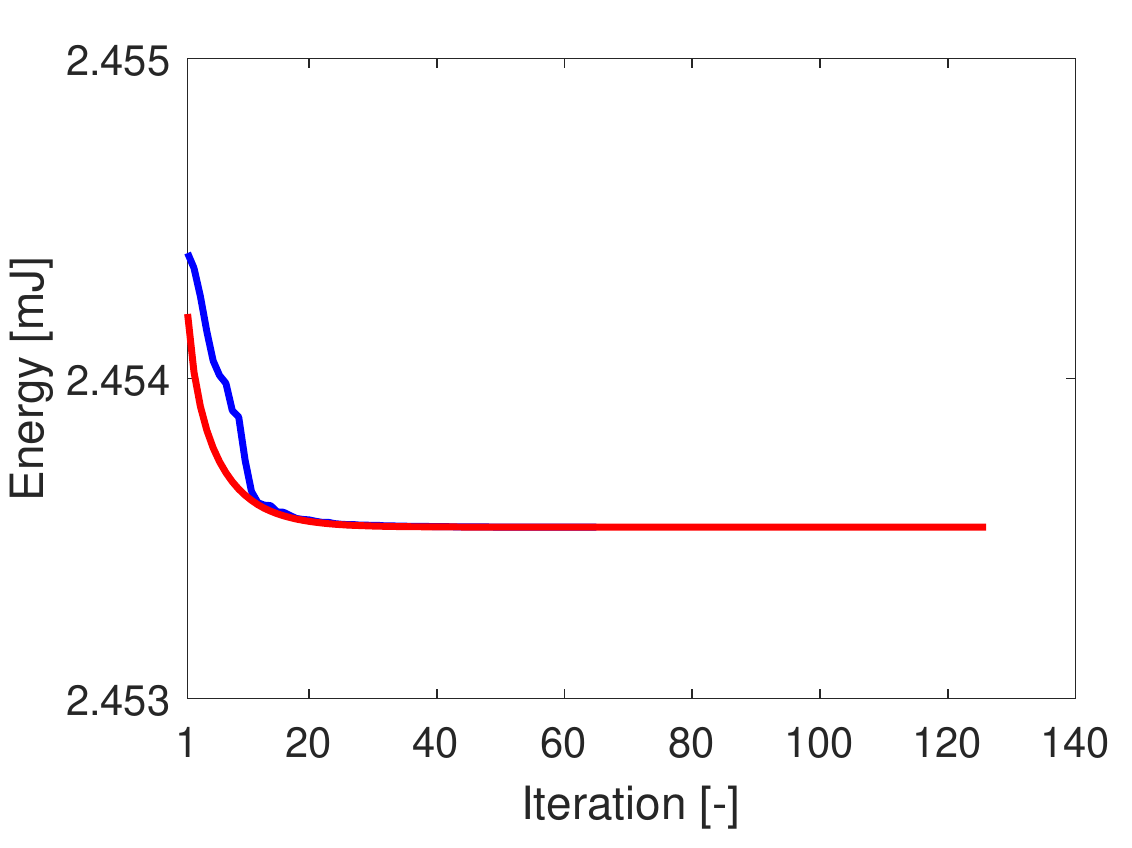}}
\subfigure{\includegraphics[width=0.325\linewidth]{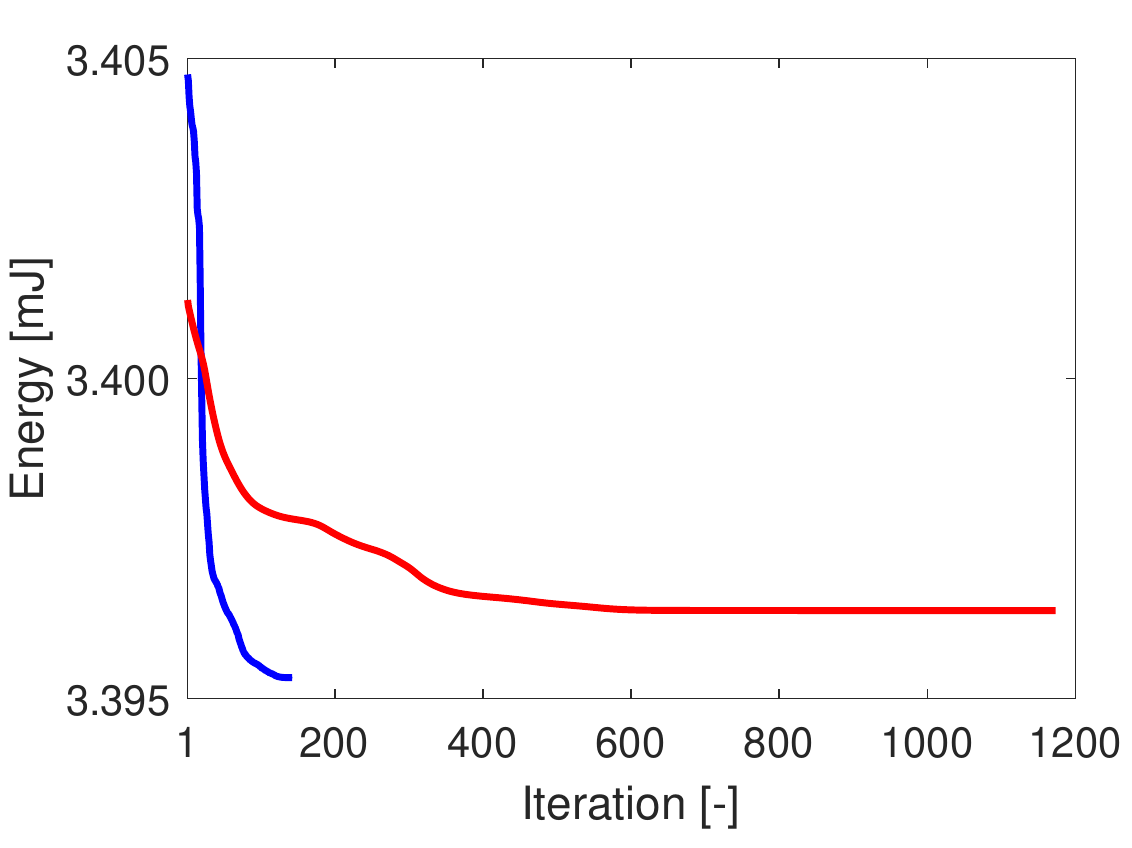}\label{fig:sen_shear_nrg_38}}
\subfigure{\includegraphics[width=0.325\linewidth]{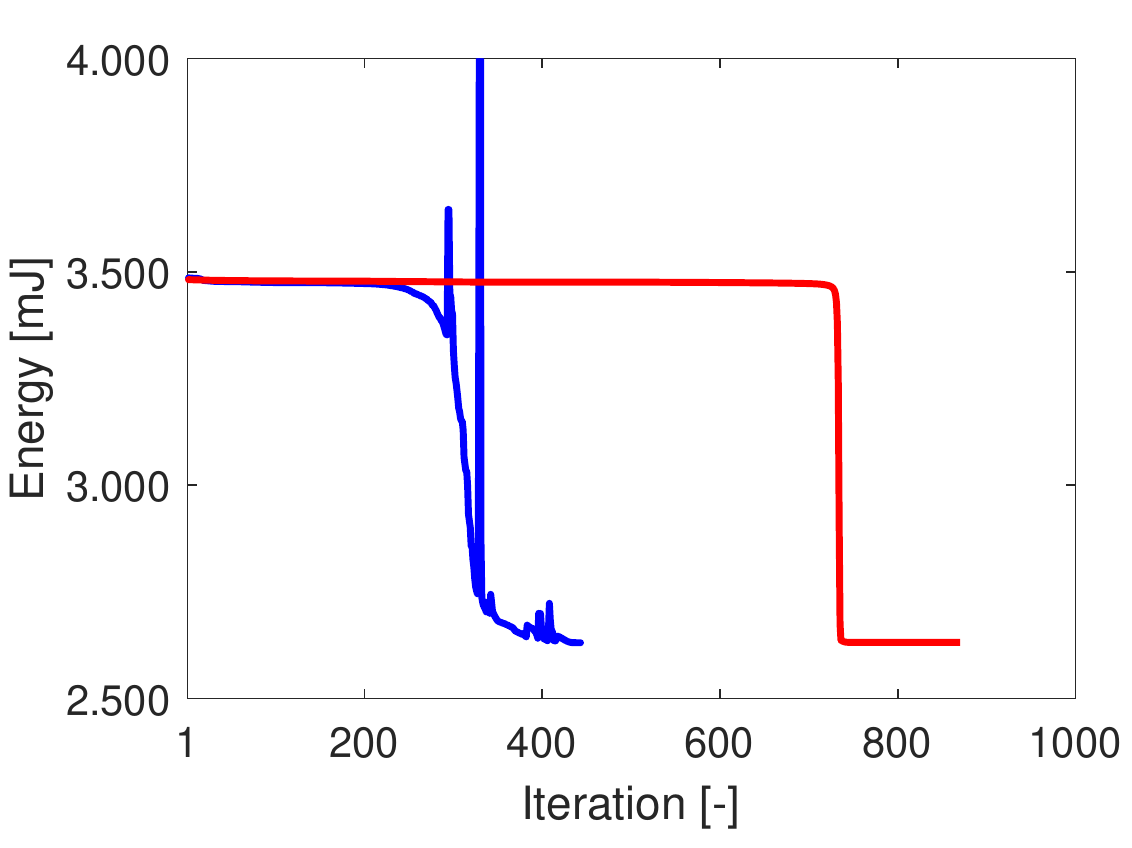}}
\addtocounter{subfigure}{-3}
\subfigure[]{\includegraphics[width=0.325\linewidth]{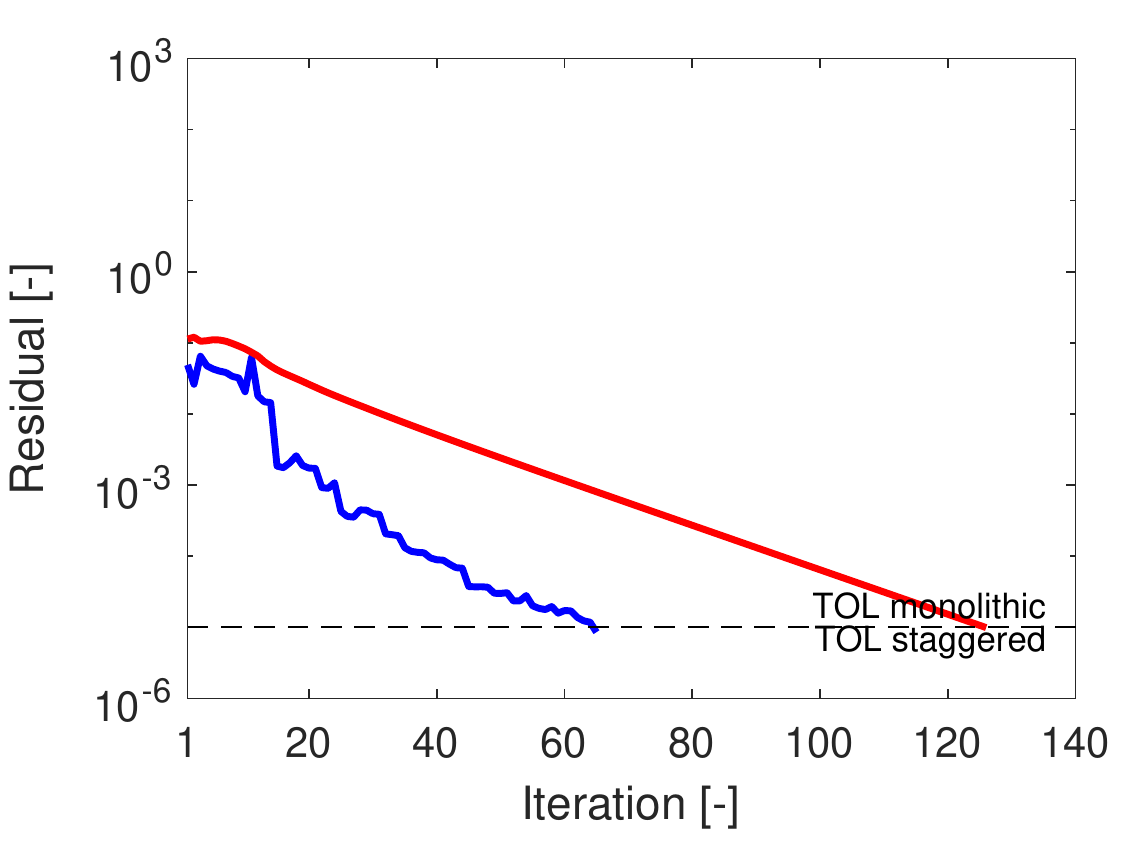}}
\subfigure[]{\includegraphics[width=0.325\linewidth]{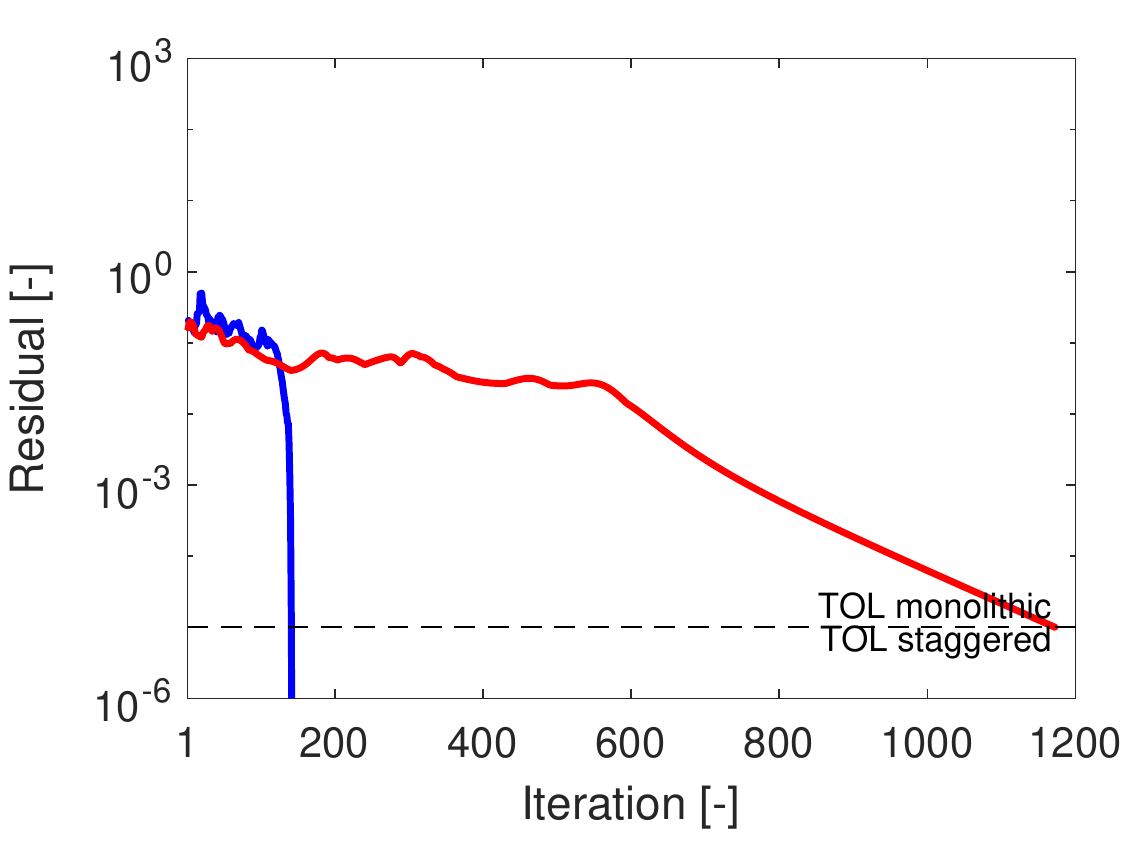}}
\subfigure[]{\includegraphics[width=0.325\linewidth]{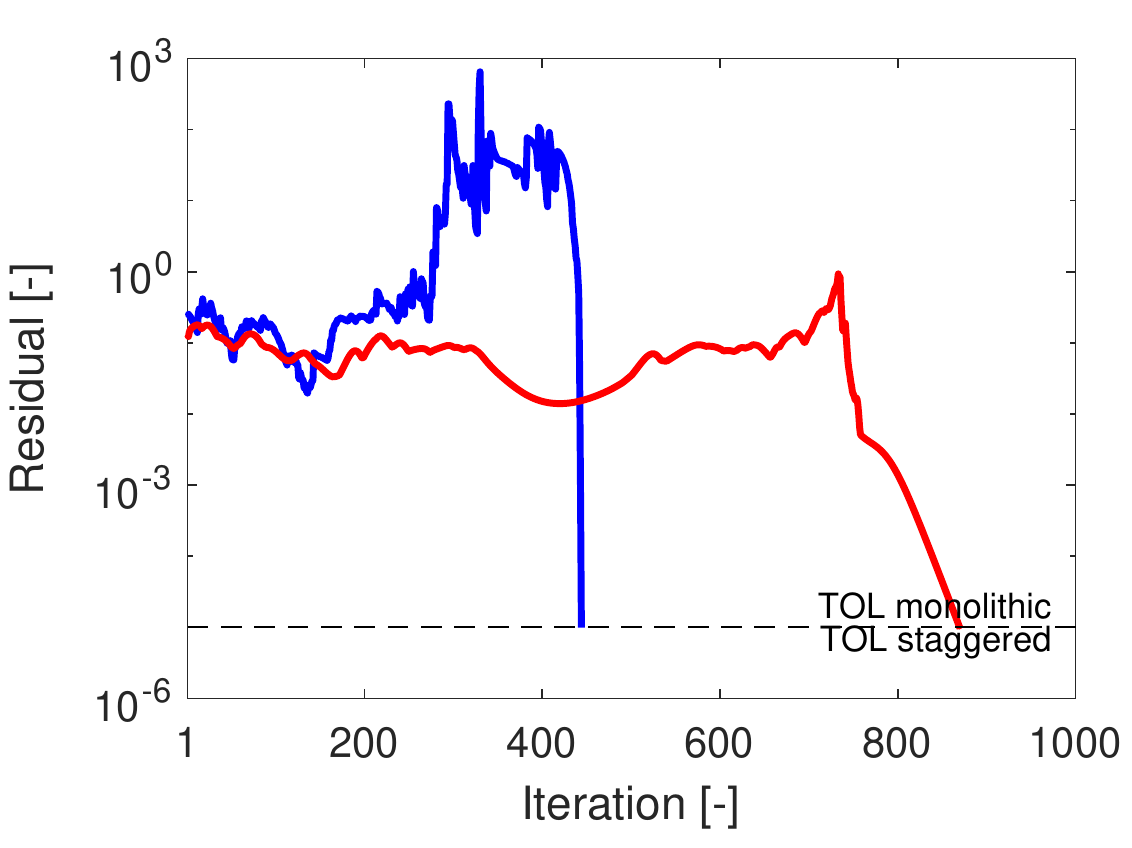}}
\vspace{-0.8em}
\caption{Evolution of the energy functional (top) and the residual (bottom) for the solution of time step (a) $n=30$, (b)~$n=38$, and (c) $n=39$ of the single-edge notched specimen subject to shear. The energy and residual of the staggered scheme and monolithic scheme are shown in red and blue, respectively.}
\label{fig:sen_shear_convergence}
\end{figure}

Fig.~\ref{fig:sen_shear_convergence} shows the evolution of the energy functional and the residual during the minimization process. In the time steps shown, the number of staggered iterations is larger than the number of iterations in the monolithic scheme. Fig.~\ref{fig:sen_shear_nrg_38} shows that the staggered and monolithic schemes result in a stationary point with a slightly different energy level. In this example, the two different solutions do not lead to noticeable differences in the load-displacement curve (Fig.~\ref{fig:sen_shear_force}) or crack phase-field profile, and the difference is attributed to numerical errors. It can be concluded that, for this example, a tolerance ${\tt TOL}=10^{-5}$ is sufficient to obtain a solution where the energy functional is stationary.
\vspace{-1.3em}

Contrary to the previous example, the difference in efficiency of the considered solution schemes is larger. The monolithic scheme requires less computational time compared to the staggered scheme that imposes irreversibility using an interior-point method. This is attributed to the less brutal crack propagation. Like the previous example, the number of staggered iterations required in the interior-point method and the augmented Lagrangian approach is similar. However, when the interior-point method is adopted, the number of sub-iterations is larger than when the augmented Lagrangian approach is used. Nevertheless, the interior-point method has the benefit of not requiring a penalty parameter.

\subsection{Asymmetrically notched specimen subject to tension}
\begin{figure}[!h]
\centering
\includegraphics[width=0.20\textwidth]{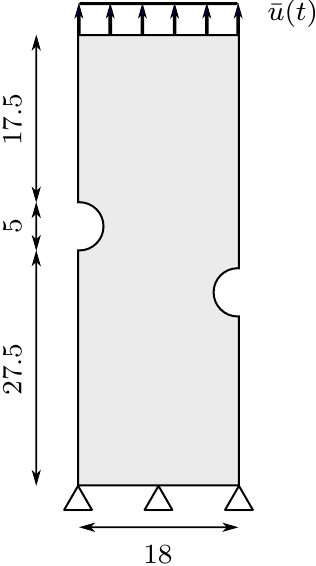}
\caption{Geometry and boundary conditions of the asymmetrically notched specimen subject to tension. The dimensions are given in~mm.}
\label{fig:asn}
\end{figure}

The example of the asymmetrically notched specimen has been studied using the phase-field approach to ductile fracture~\cite{amba15a}, and was later experimentally verified by Ambati et~al.~\cite{amba16a}. In this example, vertical displacements are imposed at the top of a notched specimen in steps of $1\times10^{-2}$~mm, while the bottom part remains fixed (Fig.~\ref{fig:asn}). Contrary to the previous examples, this specimen behaves in a ductile way, governed by an AT-1 model coupled to plasticity, and has a bulk modulus $K=71659.46$~MPa, a shear modulus $G=27296.28$~MPa, a plastic yield strength $\sigmap=345$~MPa, a plastic length scale $\etap=4$~N\textsuperscript{1/2}, a local dissipation constant $w_0=79.35$~MPa, and a damage length scale $\etad=2.217$~N\textsuperscript{1/2}. As such, the plastic characteristic length and the damage characteristic length equal $\ell_\mathrm{p}=0.215$~mm and $\ell_\mathrm{d}=0.176$~mm, respectively. These material properties correspond to the ones used by Rodr{\'i}guez~et~al.~\cite{rodr18a}, except that no plastic hardening is included in this example. In addition, plane strain conditions are assumed. A finite element discretization with 19428 bilinear quadrilateral elements and 19541 nodes is considered, resulting in an element size $h_\mathrm{c}=0.0533$~mm in the critical zone where fracture is expected to occur.

\clearpage

\begin{figure}[!h]
\centering
\subfigure[]{\includegraphics[width=0.325\linewidth]{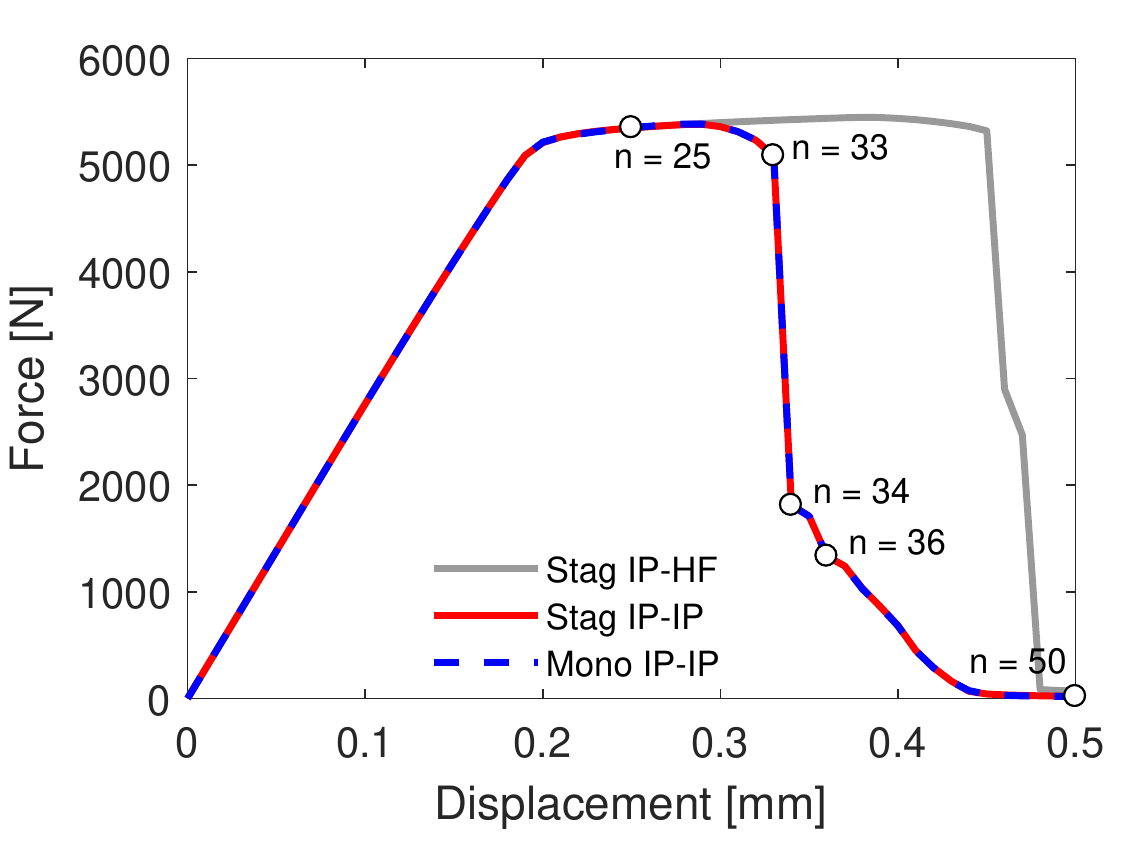}\label{fig:asn_force}}
\subfigure[]{\includegraphics[width=0.325\linewidth]{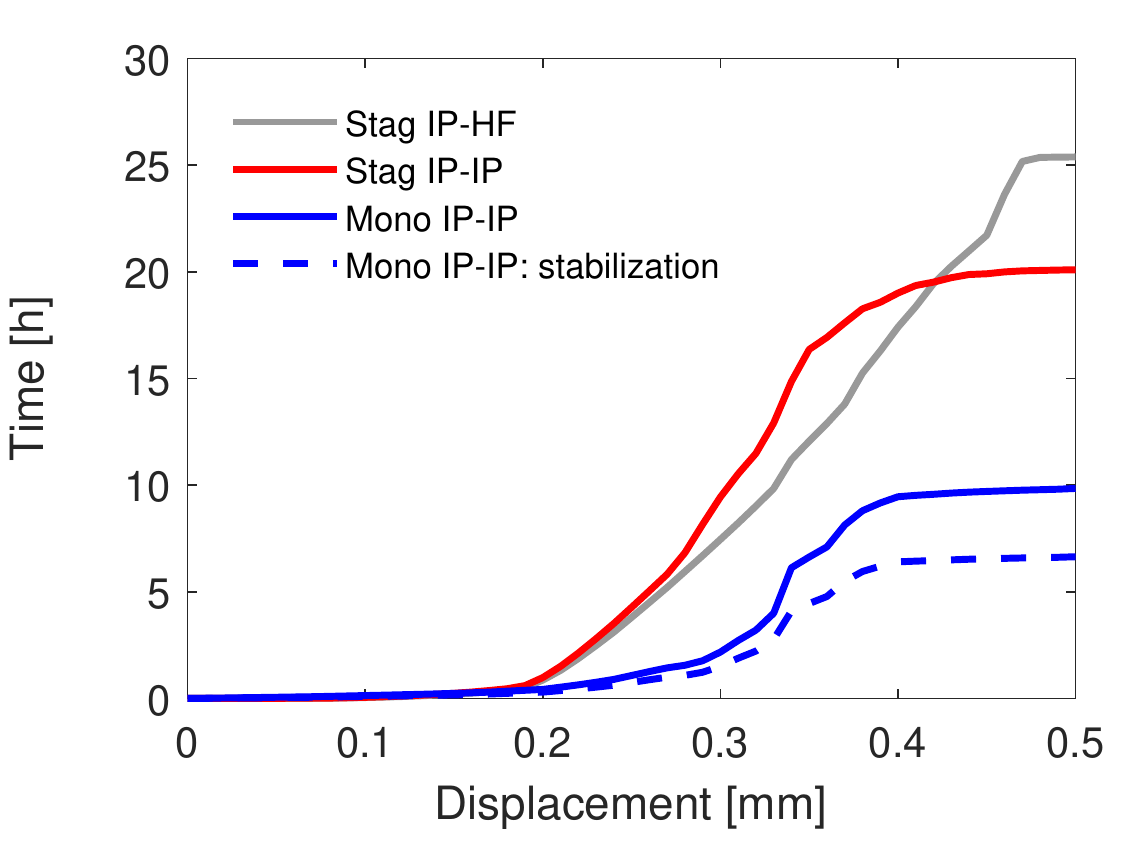}\label{fig:asn_time}}
\subfigure[]{\includegraphics[width=0.325\linewidth]{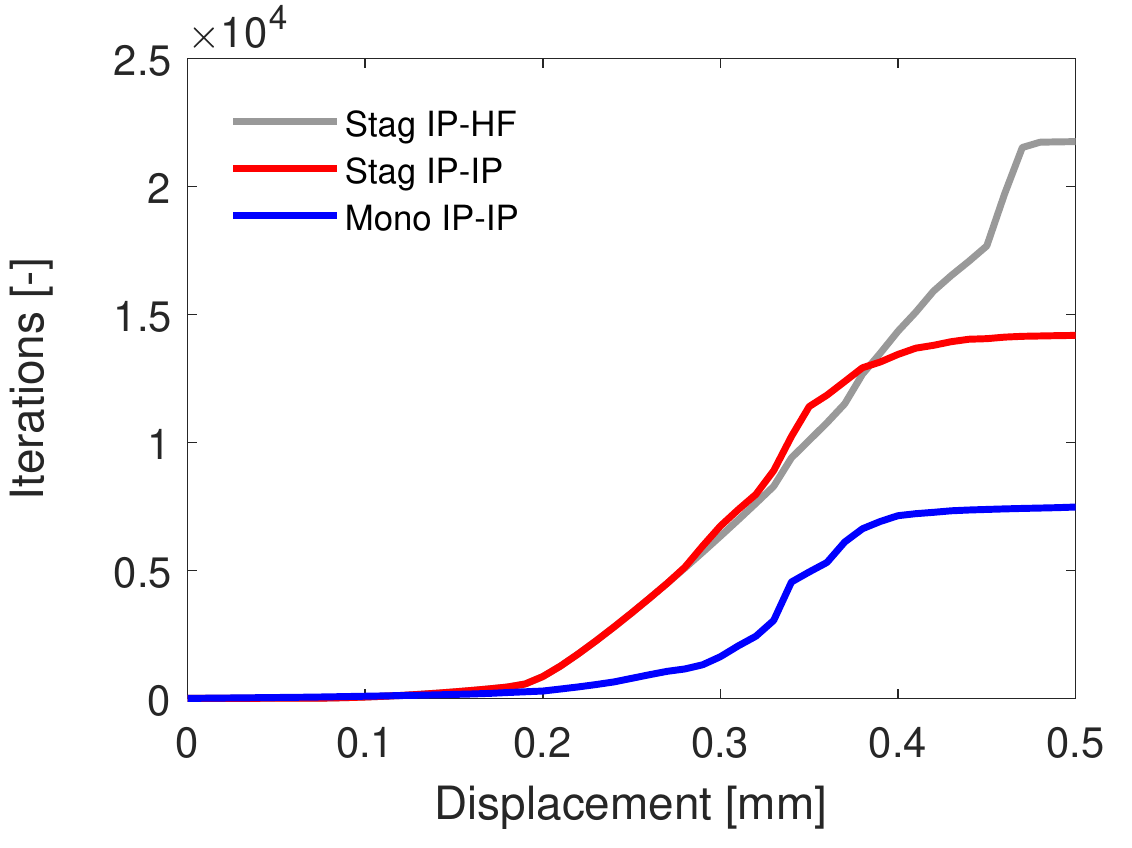}\label{fig:asn_iter}}
\caption{(a) Load-displacement curve, (b) cumulative computation time, and (c) cumulative number of iterations for the asymmetrically notched specimen subject to tension. The markers correspond to the time steps at which the profiles of the equivalent plastic strain and crack phase-field are shown in Fig.~\ref{fig:asn_dam_kappa}, and the energy evolution and residual are shown in Fig.~\ref{fig:asn_convergence}.}
\end{figure}

Fig.~\ref{fig:asn_force} shows the load-displacement curve computed using the staggered scheme and the monolithic scheme. First, the specimen shows an elastic regime as a result of the AT-1 model. This regime is followed by a perfectly plastic regime, and a subsequent fracture stage where damage initiates and the force decreases. This example shows that the use of a history field for the evolution of the crack phase-field can lead to significant deviations of the load-displacement curve. Considering that the interior-point method solves the original constrained system of variational inequalities, it can be concluded that, for this example, the history field results in a delay of crack propagation. No results of the augmented Lagrangian approach are included since no lower bound of a penalty parameter which sufficiently enforces irreversibility in ductile phase-field models is known.

The staggered computations are performed with a tolerance ${\tt TOL}=10^{-3}$ for all time steps. In time step $n=35$, the monolithic scheme progresses slowly and cannot reduce the residual of the governing equations below this value of the tolerance. Therefore, this condition is relaxed and a tolerance of ${\tt TOL}=10^{-3}$ is used for time steps $n\in[1,34]$ which is further increased to a tolerance of ${\tt TOL}=11\times10^{-3}$ in the remaining time steps of the monolithic scheme. Alternatively, a convergence criterion based on the step size of the displacement field, equivalent plastic strain, crack phase-field and dual variables can be adopted. Compared to the staggered scheme, the use of a less strict tolerance in the remaining steps does not lead to a noticeably different load-displacement curve or profile of the internal variables. Since the tolerance used for the staggered and monolithic schemes is only the same for the first 34 time steps, the computational time and cumulative number of iterations can only be objectively compared up to this point. 

Fig.~\ref{fig:asn_time} shows that the computational time of the monolithic scheme is smaller than the computational time of both staggered schemes. This also holds for the first 34 time steps where the same tolerance is adopted for all schemes. In addition, Fig.~\ref{fig:asn_iter} shows that the number of iterations required to converge to a stationary state in the monolithic scheme is smaller than the number of iterations required in the staggered scheme. The computational time for the staggered scheme where the history field is used in combination with the interior-point method is the largest, attributed to a more brutal propagation of the two cracks. Note, however, that the results obtained using the history field deviate significantly from the ones obtained using an interior-point method and, consequently, should be considered when comparing computational times or iteration numbers.

\begin{figure}[!h]
\centering
\subfigure[]{\includegraphics[width=0.325\linewidth]{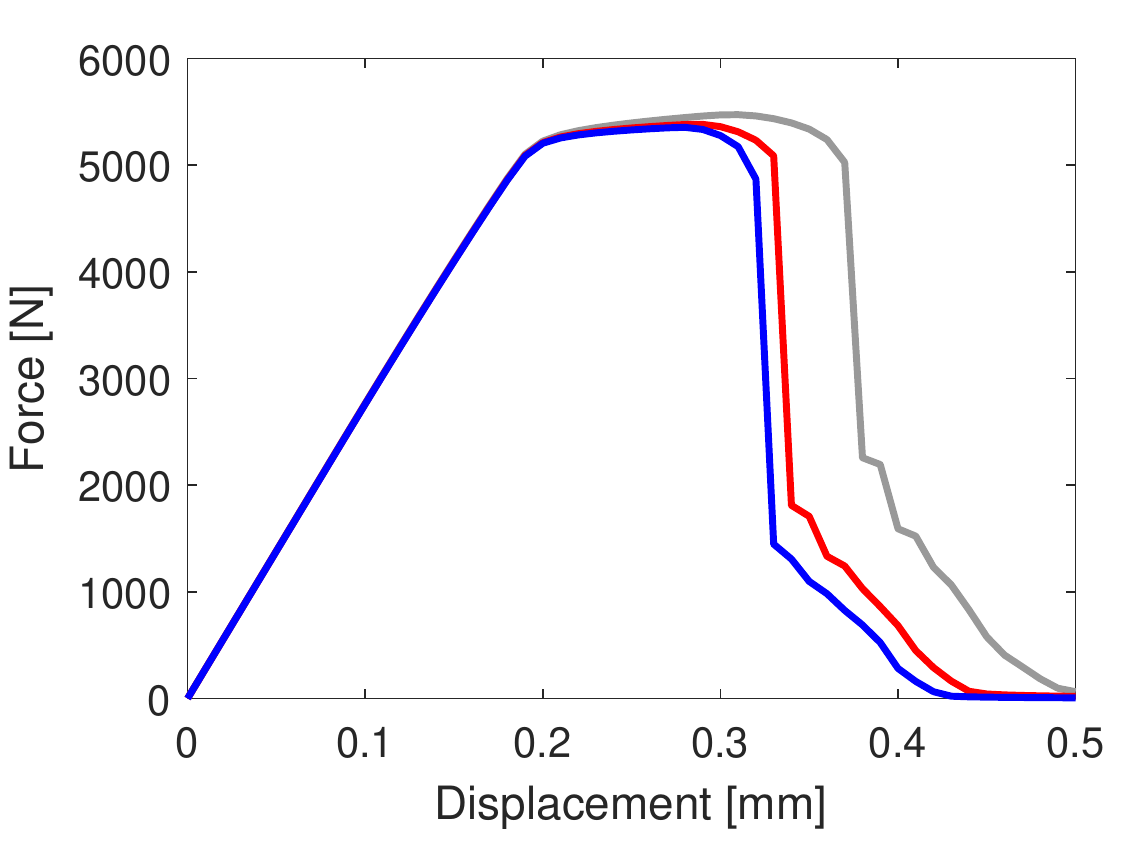}\label{fig:asn_mesh_force}}
\subfigure[]{\includegraphics[width=0.325\linewidth]{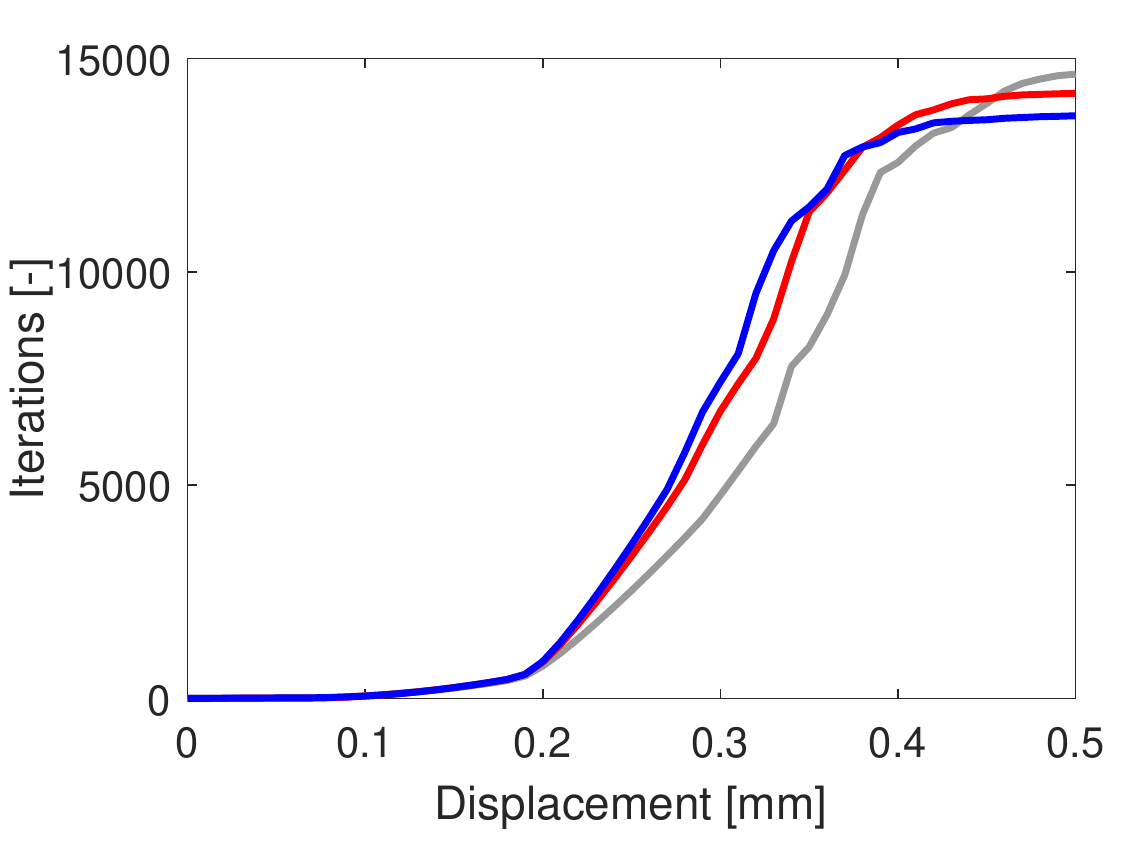}\label{fig:asn_mesh_stag}}
\subfigure[]{\includegraphics[width=0.325\linewidth]{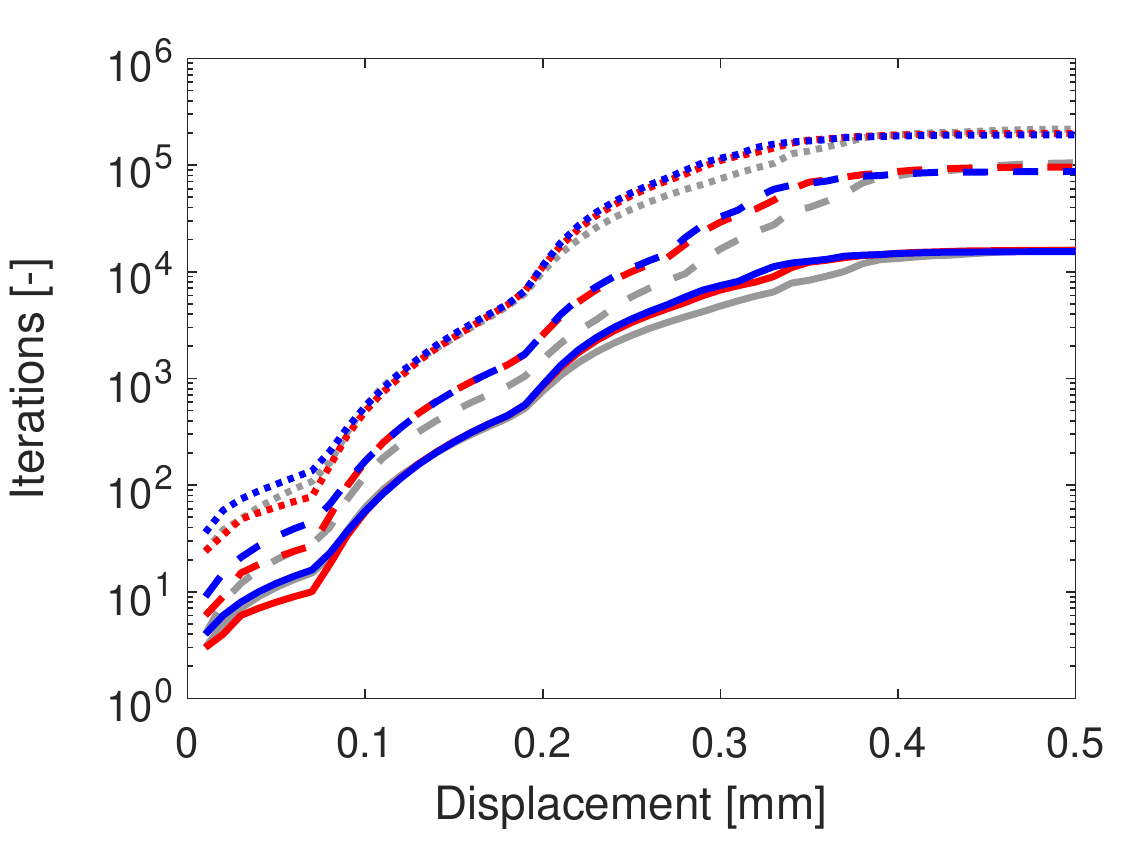}\label{fig:asn_mesh_sub}}
\caption{(a) Load-displacement curve, (b) cumulative number of staggered iterations, and (c) cumulative number of sub-iterations for the asymmetrically notched specimen subject to tension using a mesh with $h_\mathrm{c}=0.606\ell_\mathrm{d}$ (gray), $h_\mathrm{c}=0.303\ell_\mathrm{d}$ (red), and $h_\mathrm{c}=0.152\ell_\mathrm{d}$ (blue). The solid, dotted and dashed lines correspond to equilibrium, plastic, and damage sub-iterations, respectively.}
\end{figure}

Fig.~\ref{fig:asn_mesh_force} shows the load-displacement curve obtained using the interior-point method in a staggered scheme for different mesh sizes. Mesh-convergence can be observed. In addition, Figs.~\ref{fig:asn_mesh_stag} and~\ref{fig:asn_mesh_sub} show the cumulative number of staggered iterations and sub-iterations, respectively. The total number of staggered iterations and sub-iterations is not very sensitive to the mesh size. In addition, the average number of damage and plastic sub-iterations per staggered iteration does not change significantly.

\begin{figure}[!h]
\centering
\hspace{1.5cm}\includegraphics[width=0.8\linewidth]{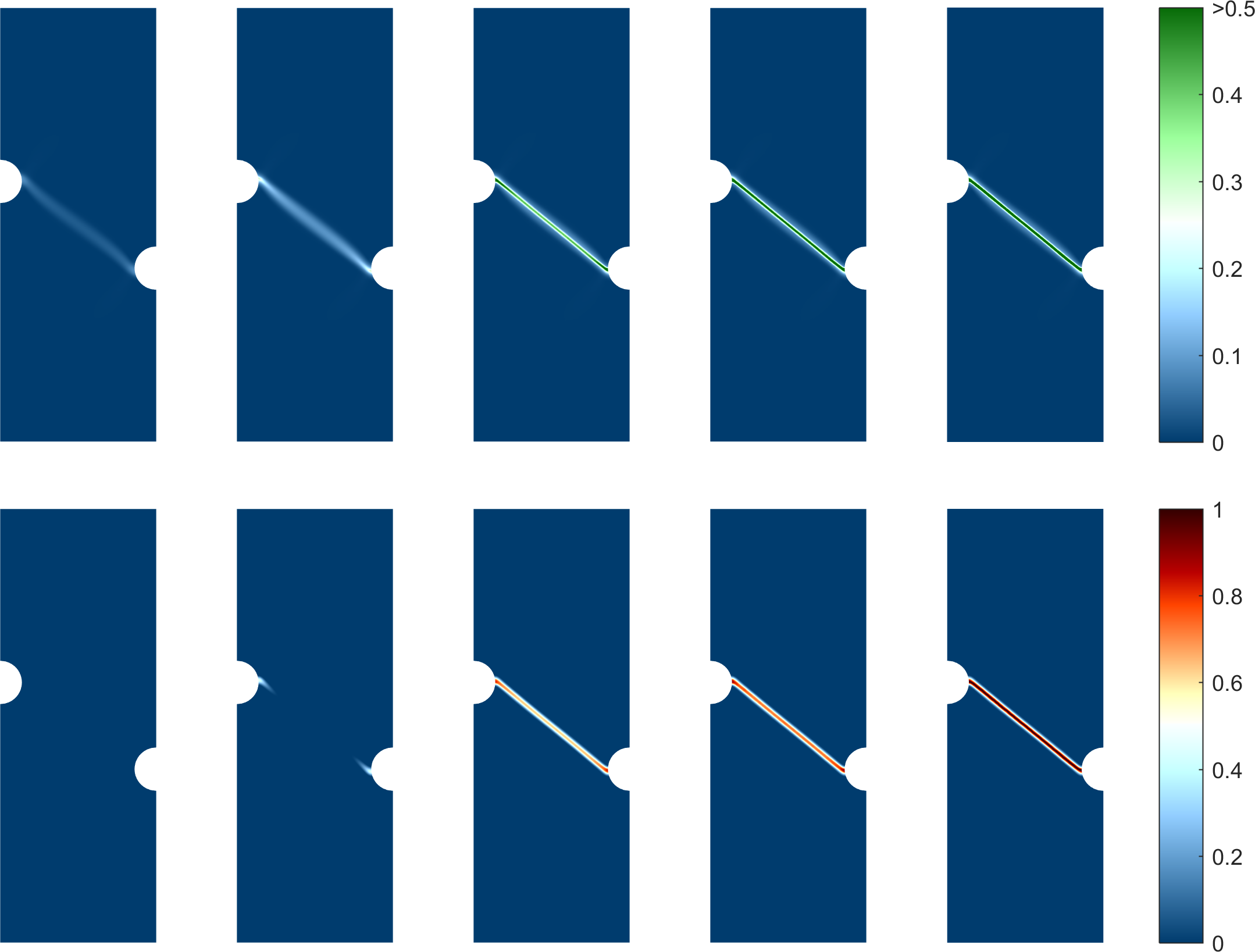}
\put(-390,156){(a)}
\put(-320,156){(b)}
\put(-250,156){(c)}
\put(-180,156){(d)}
\put(-107,156){(e)}
\put(-390,7){(f)}
\put(-320,7){(g)}
\put(-250,7){(h)}
\put(-177,7){(i)}
\put(-106,7){(j)}
\caption{Profiles of the equivalent plastic strain (top) and crack phase-field (bottom) for the asymmetrically notched specimen subject to tension, computed using the staggered scheme at time step $n=25$, $n=33$, $n=34$, $n=36$, and $n=50$ (left to right).}
\label{fig:asn_dam_kappa}
\end{figure}
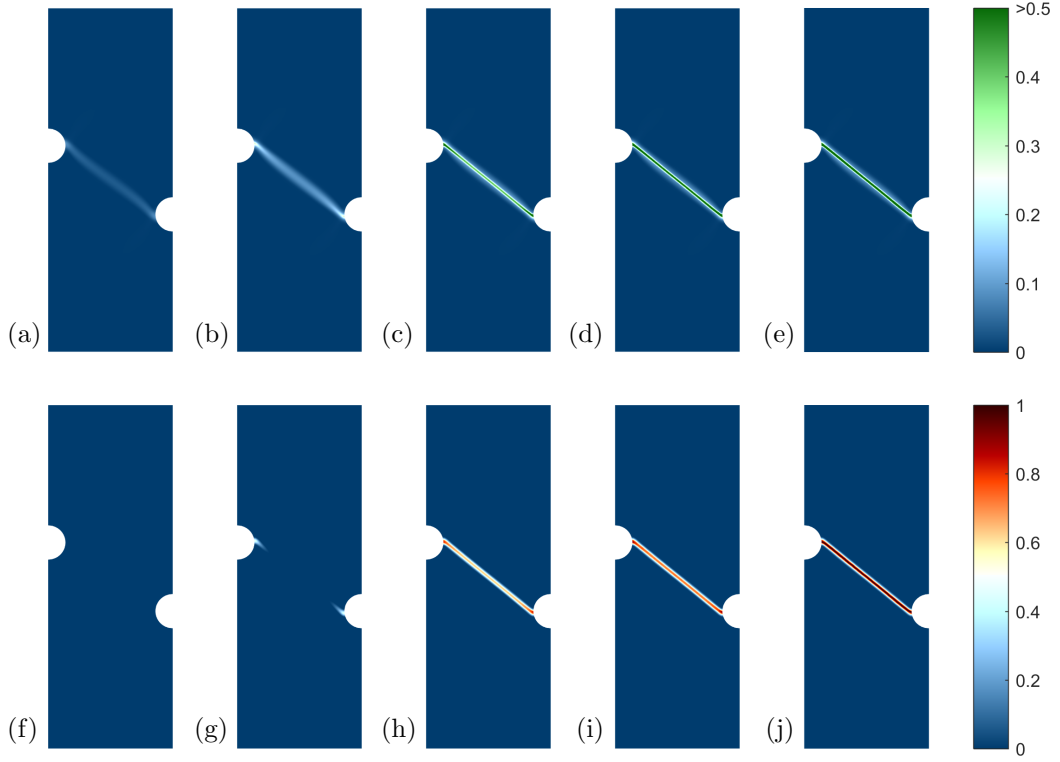

Fig.~\ref{fig:asn_dam_kappa} shows the profiles of the equivalent plastic strain and the crack phase-field, computed using the interior-point method in a staggered scheme at different time steps. Fig.~\ref{fig:asn_dam_kappa}(a) shows that in the beginning of the plastic stage, a shear band starts to develop between the notches. The thickness of this band is governed by the plastic length scale $\etap$. At this stage, the damage threshold $w_0$ has not been reached and the crack has not initiated yet (Fig.~\ref{fig:asn_dam_kappa}(f)). In a later stage, damage initiate near the notches (Fig.~\ref{fig:asn_dam_kappa}(g)), and two cracks eventually coalesce (Figs.~\ref{fig:asn_dam_kappa}(h) and \ref{fig:asn_dam_kappa}(i)). The plastic length scale vanishes in the region of the crack as a result of its degradation in the energy functional~\eqref{eq:functional}, and consequently, the plastic strain evolution localizes within the crack region (Figs.~\ref{fig:asn_dam_kappa}(c) and~\ref{fig:asn_dam_kappa}(d)). There are no significant differences between the profiles of the equivalent plastic strain or the crack phase-field obtained using the staggered scheme and the monolithic scheme in case an interior-point method is used to impose both plastic and damage irreversibility. When the history field is used to impose damage irreversibility, the crack phase-field profile is slightly narrower. The corresponding figures are not shown for brevity.

\clearpage

\begin{figure}[!h]
\centering
\subfigure{\includegraphics[width=0.325\linewidth]{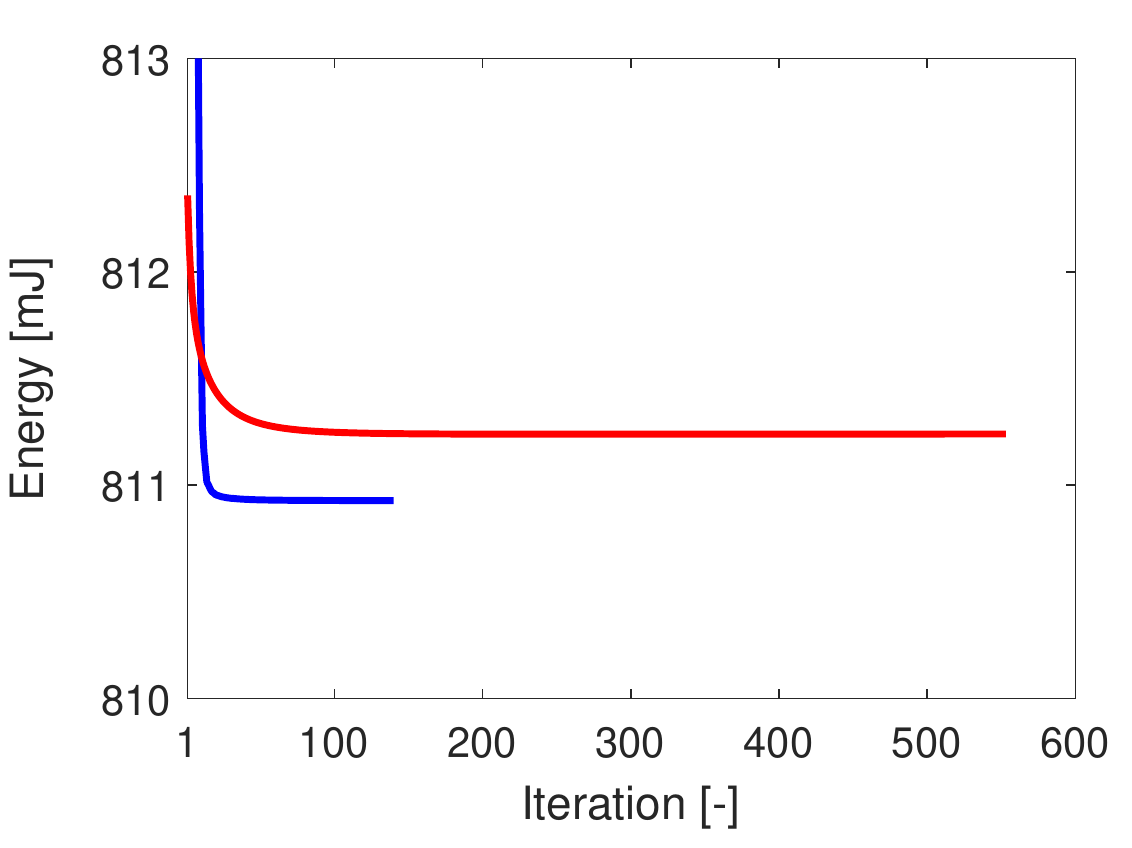}}
\subfigure{\includegraphics[width=0.325\linewidth]{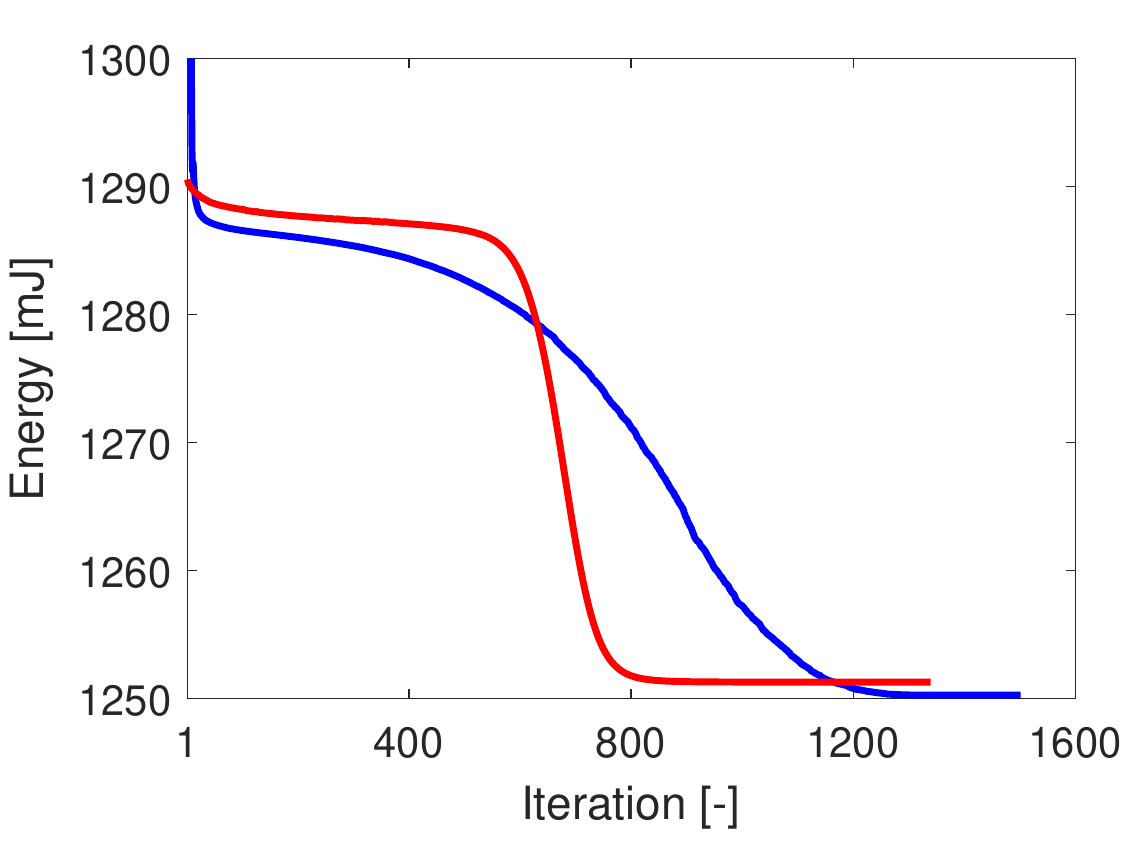}\label{fig:asn_nrg_34}}
\subfigure{\includegraphics[width=0.325\linewidth]{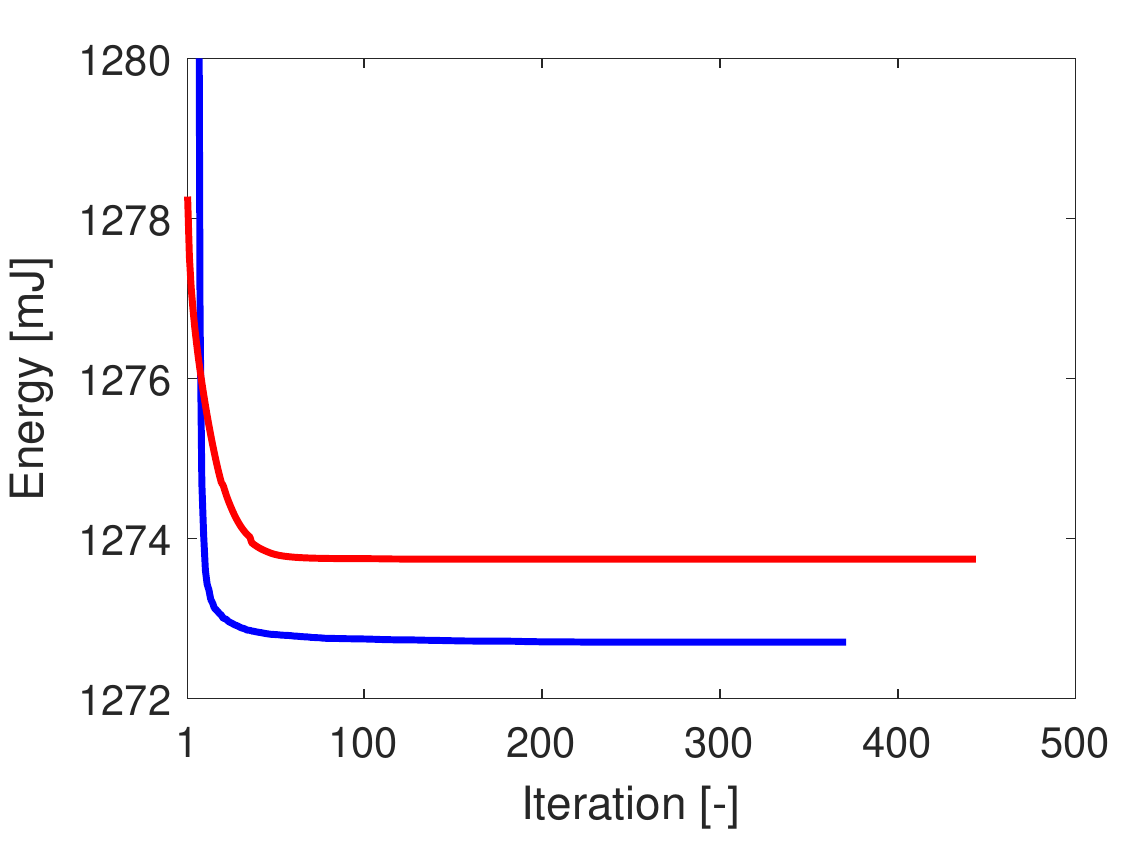}\label{fig:asn_nrg_36}}
\addtocounter{subfigure}{-3}
\subfigure[]{\includegraphics[width=0.325\linewidth]{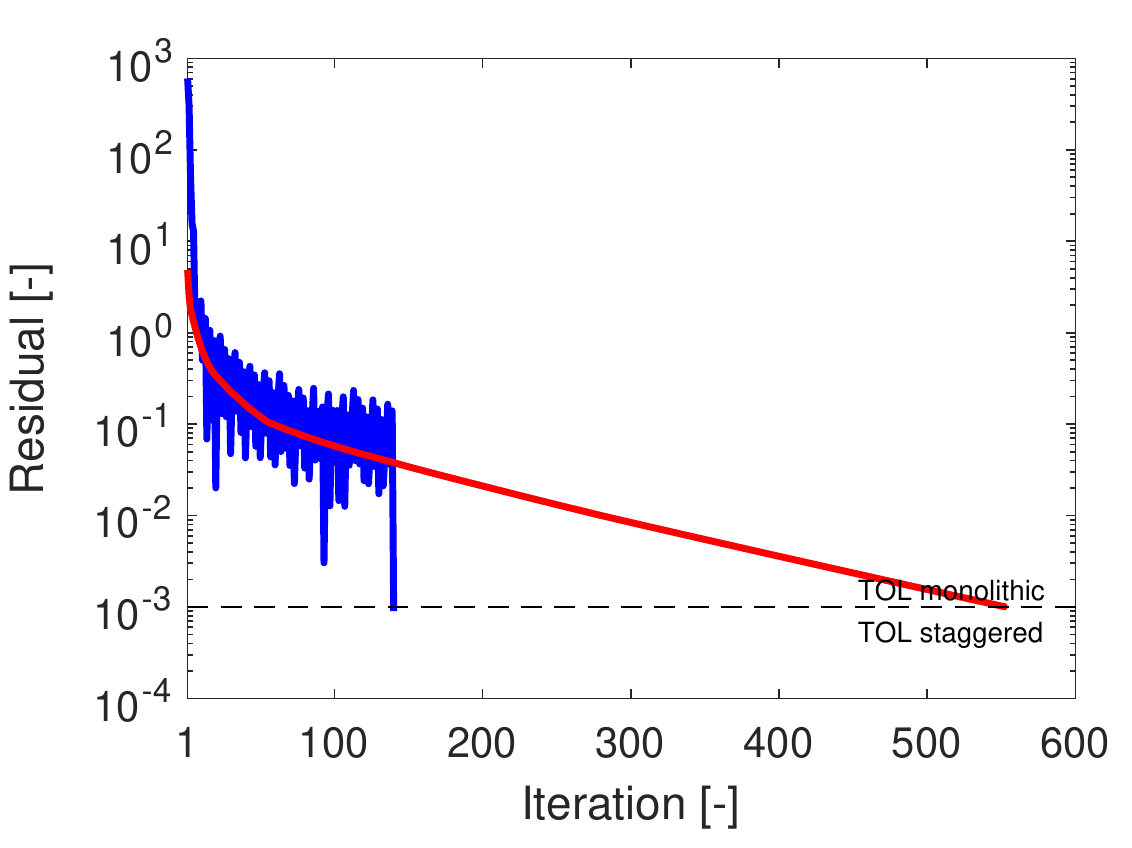}}
\subfigure[]{\includegraphics[width=0.325\linewidth]{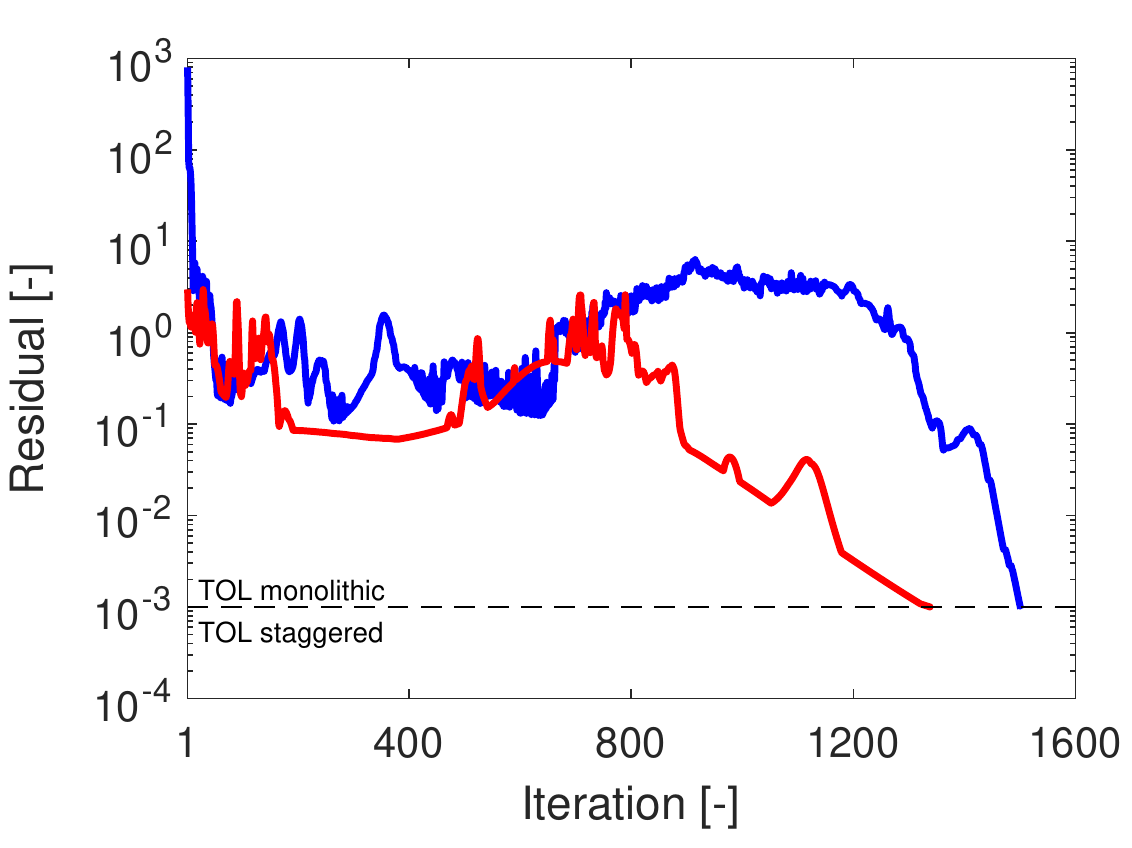}}
\subfigure[]{\includegraphics[width=0.325\linewidth]{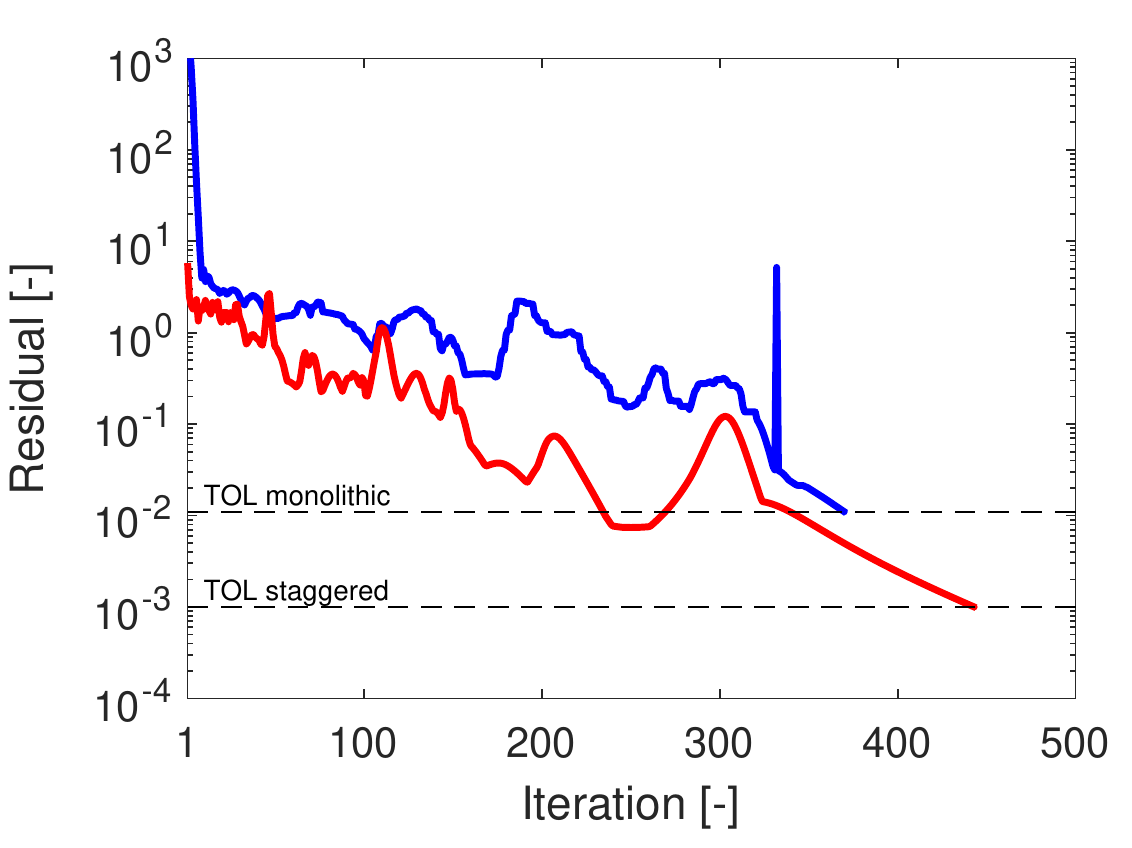}}
\caption{Evolution of the energy functional (top) and the residual (bottom) for the solution of time step (a) $n=25$, (b)~$n=34$, and (c) $n=36$ of the asymmetrically notched specimen subject to tension. The energy and residual of the staggered scheme and monolithic scheme are shown in red and blue, respectively.}
\label{fig:asn_convergence}
\end{figure}

Fig.~\ref{fig:asn_convergence} shows the evolution of the energy functional and the residual during the minimization process. Both schemes based on interior-point methods result in a stationary point with a different energy level. Similar to the second example, the two different solutions do not lead to noticeable differences in the load-displacement curve (Fig.~\ref{fig:asn_force}), or the profiles of the equivalent plastic strain and the crack phase-field. It can be concluded that, for this example, a tolerance ${\tt TOL}=10^{-3}$ is sufficient to obtain a solution where the energy functional is stationary. In addition, when the tolerance is relaxed to $11\times10^{-3}$, it can be observed that the energy is stationary in time step $n=36$ (Fig.~\ref{fig:asn_nrg_36}).

\begin{figure}[!h]
\centering
\subfigure[]{\includegraphics[width=0.325\linewidth]{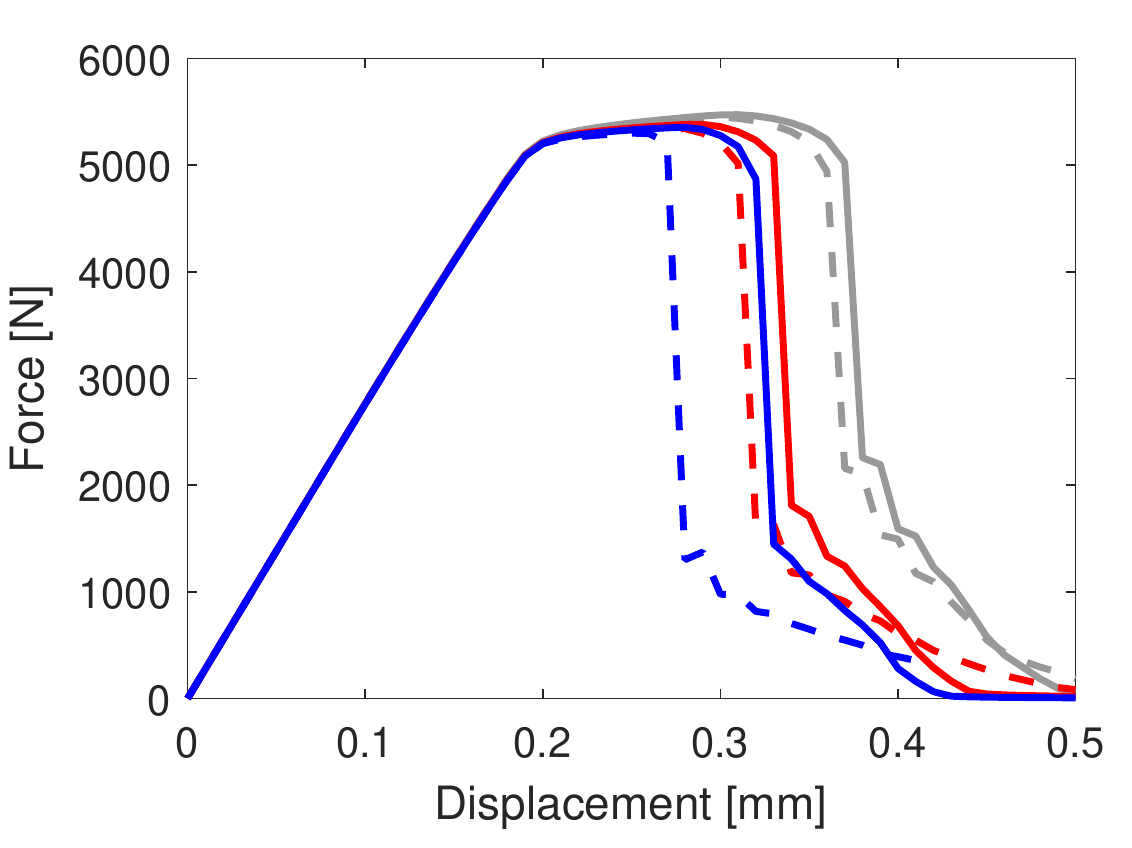}\label{fig:asn_lengthscale}}
\subfigure[]{\includegraphics[width=0.325\linewidth]{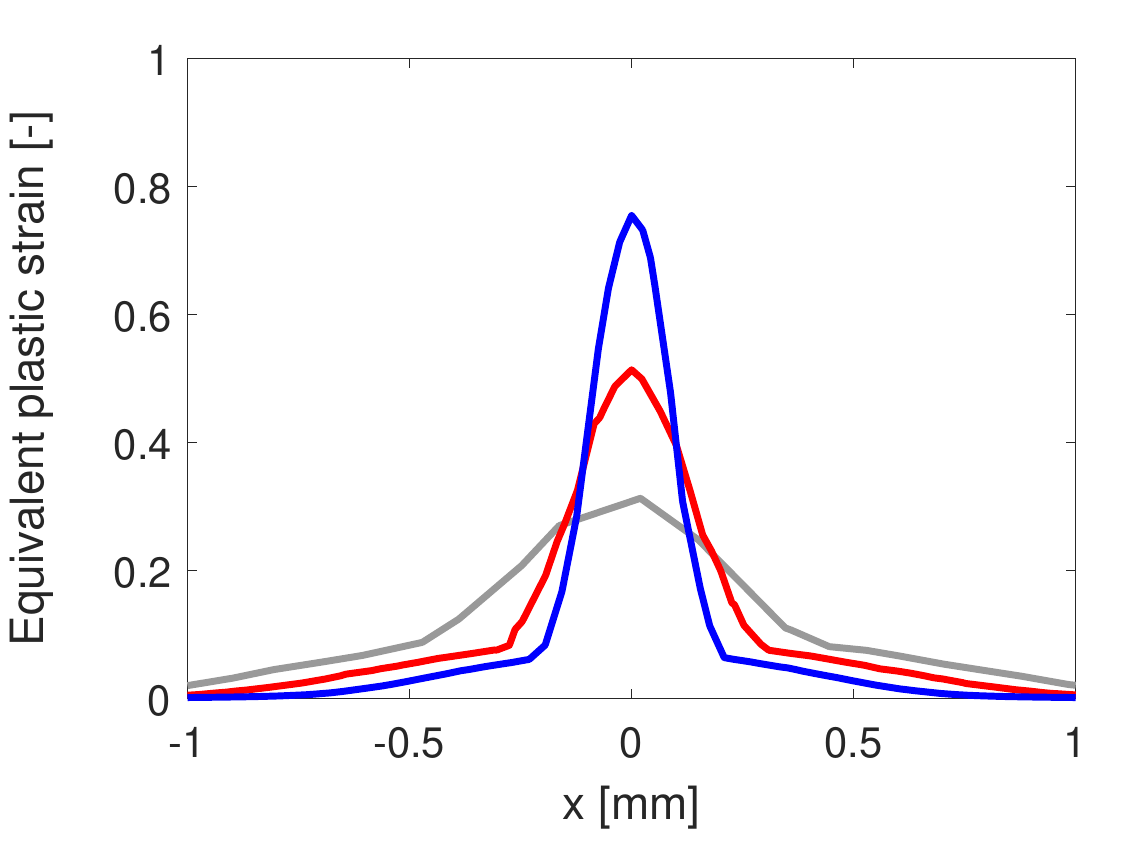}\label{fig:asn_local}}
\subfigure[]{\includegraphics[width=0.325\linewidth]{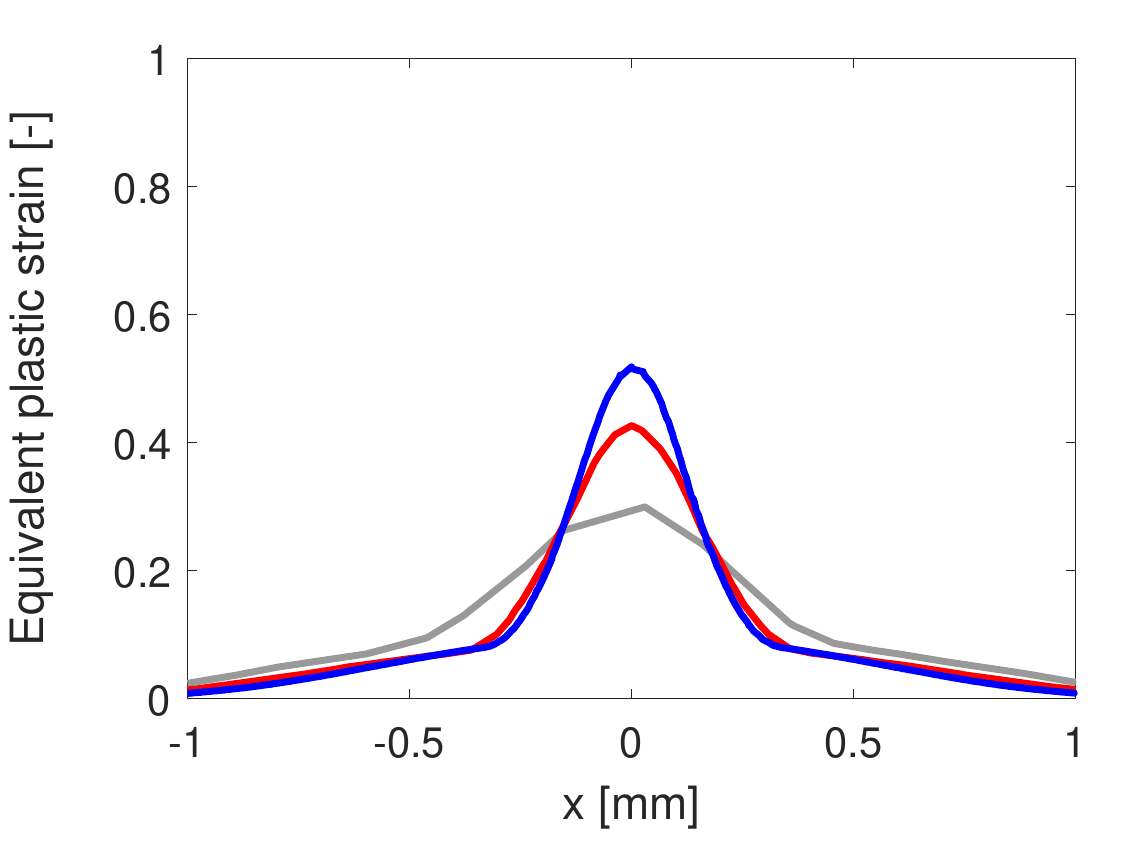}\label{fig:asn_nonlocal}}
\caption{(a) Load-displacement curve, (b) profile of the equivalent plastic strain for $\etap=0$~N\textsuperscript{1/2}, and (c) profile of the equivalent plastic strain for $\etap=4$~N\textsuperscript{1/2} for the asymmetrically notched specimen subject to tension using a mesh with $h_\mathrm{c}=0.606\ell_\mathrm{d}$ (gray), $h_\mathrm{c}=0.303\ell_\mathrm{d}$ (red), and $h_\mathrm{c}=0.152\ell_\mathrm{d}$ (blue). The solid and dashed lines correspond to non-local and local plasticity, respectively.}
\label{fig:asn_plasticstudy}
\end{figure}

Fig.~\ref{fig:asn_plasticstudy} compares the results for different meshes of the asymmetrically notched specimen subject to tension for both local and non-local plasticity. The load-displacement curve does not show any convergence upon mesh refinement in the local case (Fig.~\ref{fig:asn_lengthscale}). In addition, in the local case, numerical issues hamper the computation of the load-displacement curve in time step 43 for the finest mesh. Figs.~\ref{fig:asn_local} and~\ref{fig:asn_nonlocal} show the profile of the equivalent plastic strains at the time step that corresponds to brutal crack propagation, for the local plasticity and non-local plasticity model, respectively. In case of local plasticity, the equivalent plastic strains increase upon mesh refinement. The increase is smaller for the non-local counterpart, where the localization is regularized by the plastic length scale $\etap$.

In summary, a significant delay in crack propagation is observed when a history field is used to impose irreversibility of the crack phase-field. This is attributed to the modifications of the governing equations. For this example, the monolithic scheme in combination with the interior-point method appears to be the most efficient, but cannot decrease the residual to the same value as the staggered scheme. In addition, the present method enables the use of non-local plasticity models, which lead to mesh-insensitive results and a more robust solution scheme.

\section{Conclusions} \label{sec:conclusions}
In this paper, the system of variational inequalities that governs brittle and ductile fracture processes in phase-field models is rigorously solved using interior-point methods. With this method, a sequence of perturbed constraints is considered, which, in the limit, recovers the unperturbed constraints and solves the exact constrained problem. As such, no penalty parameters or modifications of the governing equations are required. The interior-point method is applied in both a staggered and a monolithic scheme, and systematically compared to a benchmark staggered scheme where a history field or an augmented Lagrangian is used to impose irreversibility of the crack phase-field. In order to stabilize the monolithic scheme, a stabilization of the Hessian matrix is suggested.

The presented algorithms are tested on three benchmark problems. It is observed that the computational efficiency of the staggered scheme and monolithic scheme is problem-dependent. When an interior-point method is used to impose irreversibility, the computational time of the monolithic scheme is similar to or smaller than that of the staggered scheme for the time steps where cracks do not propagate brutally. In the time steps corresponding to brutal crack propagation, the computational efficiency appears problem-dependent. In the first example of a brittle specimen, the use of a history field leads to a wider crack phase-field profile. For both examples with brittle specimens, the interior-point method in a staggered scheme is slower than the irreversibility approach that uses an augmented Lagrangian. In the example of a ductile specimen, the use of a history field, in combination with an interior-point method for plastic irreversibility, modifies the results significantly, resulting in a delay of crack propagation and a narrower crack phase-field profile. Alternatively, adopting the interior-point method for both damage and plastic irreversibility allows to rigorously solve the equations of the non-local plasticity model for ductile fracture. As such, a more robust solution scheme is obtained and mesh-sensitivity is avoided. Nevertheless, the use of interior-point methods does require a large number of sub-iterations and, consequently, a higher computational cost than penalty methods.

\section*{References}
\bibliographystyle{model1-num-names}
\bibliography{../../../../../literature/abbrev,../../../../../literature/literature_no_caps} 

\end{document}